\documentclass[a4paper,12pt,twoside]{article}
\usepackage[francais]{babel}
\usepackage{amssymb, amsmath, eucal}
\usepackage[all]{xy}
\usepackage{mathrsfs}
\usepackage[dvips]{graphics}
\begin{document}
\textwidth15.5cm
\textheight22.5cm
\voffset=-13mm
\newtheorem{The}{Theorem}[section]
\newtheorem{Lem}[The]{Lemma}
\newtheorem{Prop}[The]{Proposition}
\newtheorem{Cor}[The]{Corollary}
\newtheorem{Rem}[The]{Remark}
\newtheorem{Obs}[The]{Observation}
\newtheorem{SConj}[The]{Standard Conjecture}
\newtheorem{Titre}[The]{\!\!\!\! }
\newtheorem{Conj}[The]{Conjecture}
\newtheorem{Question}[The]{Question}
\newtheorem{Prob}[The]{Problem}
\newtheorem{Def}[The]{Definition}
\newtheorem{Not}[The]{Notation}
\newtheorem{Claim}[The]{Claim}
\newtheorem{Conc}[The]{Conclusion}
\newtheorem{Ex}[The]{Example}
\newtheorem{Fact}[The]{Fact}
\newcommand{\C}{\mathbb{C}}
\newcommand{\R}{\mathbb{R}}
\newcommand{\N}{\mathbb{N}}
\newcommand{\Z}{\mathbb{Z}}
\newcommand{\Q}{\mathbb{Q}}
\newcommand{\Proj}{\mathbb{P}}

\begin{center}

{\Large\bf Transcendental K\"ahler Cohomology Classes}

\end{center}

\begin{center}

 {\large Dan Popovici}

\end{center}

\vspace{1ex}

\noindent {\small {\bf Abstract.} Associated with a smooth, $d$-closed $(1, \, 1)$-form $\alpha$ of possibly non-rational De Rham cohomology class on a compact complex manifold $X$ is a sequence of asymptotically holomorphic complex line bundles $L_k$ on $X$ equipped with $(0, \, 1)$-connections $\bar\partial_k$ for which $\bar\partial_k^2\neq 0$. Their study was begun in the thesis of L. Laeng. We propose in this non-integrable context a substitute for H\"ormander's familiar $L^2$-estimates of the $\bar\partial$-equation of the integrable case that is based on analysing the spectra of the Laplace-Beltrami operators $\Delta_k''$ associated with $\bar\partial_k$. Global approximately holomorphic peak sections of $L_k$ are constructed as a counterpart to Tian's holomorphic peak sections of the integral-class case. Two applications are then obtained when $\alpha$ is strictly positive\!: a Kodaira-type approximately holomorphic projective embedding theorem  and a Tian-type almost-isometry theorem for compact K\"ahler, possibly non-projective, manifolds. Unlike in similar results in the literature for symplectic forms of integral classes, the peculiarity of $\alpha$ lies in its transcendental class. This approach will be hopefully continued in future work by relaxing the positivity assumption on $\alpha$.}

\vspace{3ex}

\section{Introduction}\label{section:introduction}

 Let $X$ be a compact complex manifold, $\mbox{dim}_{\C}X=n$. Fix an arbitrary Hermitian metric $\omega$ on $X$ (which will be identified throughout with the corresponding $C^{\infty}$ positive-definite $(1, \, 1)$-form on $X$). Let $\alpha$ be any $C^{\infty}$ $d$-closed $(1, \, 1)$-form on $X$. Thus its De Rham cohomology $2$-class $\{\alpha\}\in H^2_{DR}(X, \, \R)$ may be a {\it transcendental} (i.e. non-rational) class.

 For every $q=0, \dots , n$, denote as in [Dem85a] by $X(\alpha, \, q)\subset X$ the open subset of points $z\in X$ such that $\alpha$ has $q$ negative and $n-q$ positive eigenvalues (counted with multiplicities) at $z$. Set $X(\alpha, \, \leq 1)\!\!:=X(\alpha, \, 0)\cup X(\alpha, \, 1)$. When $\{\alpha\}$ is an integral class (i.e. $\{\alpha\}$ is the first Chern class $c_1(L)$ of a holomorphic line bundle $L$ on $X$), it follows from Demailly's holomorphic Morse inequalities [Dem85a] that the (a priori very weak) positivity assumption $\int_{X(\alpha, \, \leq 1)}\alpha^n>0$ suffices to guarantee that $L$ is big. Moreover, it is well known that the first Chern class of any {\it big} holomorphic line bundle $L\rightarrow X$ contains a K\"ahler current $T$ (i.e. a $d$-closed $(1, \, 1)$-current $T$ such that $T\geq\delta\,\omega$ on $X$ for some constant $\delta>0$). 

 On the other hand, the existence of a K\"ahler current (possibly of transcendental class) on $X$ is equivalent, thanks to [DP04], to $X$ being a {\it class} ${\cal C}$ manifold (i.e. bimeromorphically equivalent to a compact K\"ahler manifold). Thus extending to arbitrary classes $\{\alpha\}$ the (by now) classical results for integral classes alluded to above is of the utmost importance in the study of the geometry of $X$. In this paper we make the first moves towards an eventual resolution of Demailly's conjecture on {\it transcendental Morse inequalities}.

\begin{Conj}(Demailly)\label{Conj:Morse} Let $\alpha$ be any $C^{\infty}$ $d$-closed $(1, \, 1)$-form on $X$ of arbitrary (i.e. possibly non-rational) cohomology class $\{\alpha\}\in H^2_{DR}(X, \, \R)$. If

\begin{equation}\label{eqn:Dem-posint-assumption}\int\limits_{X(\alpha, \, \leq 1)}\alpha^n>0,\end{equation}

\noindent then there exists a K\"ahler current in the class $\{\alpha\}$.

\end{Conj}

 After the results we have obtained in [Pop08] and [Pop09], this conjecture of Demailly is the last missing link in a (hopefully forthcoming, but still elusive) resolution of the following long-conjectured fact.

\begin{Conj}(standard) Suppose that in a complex analytic family of compact complex manifolds $(X_t)_{t\in\Delta}$ over the unit disc $\Delta\subset\C$ the fibre $X_t$ is K\"ahler for every $t\in\Delta\setminus\{0\}$. Then $X_0$ is a class ${\cal C}$ manifold.

\end{Conj}

 The results we obtain in this paper build on earlier work by L. Laeng [Lae02] whose set-up and main result we now summarise. They will serve as the starting point of the present work.

\subsection{Setting considered by L. Laeng in [Lae02]}\label{subsection:setting}

 By [Lae02, Th\'eor\`eme 1.3, p. 57], one can find an infinite subset $S\subset\N^{\star}$ (that will be assumed without loss of generality to be $\N^{\star}$) and a sequence $(\alpha_k)_{k\in S}$ of $C^{\infty}$ $d$-closed $2$-forms (in general not of type $(1, \, 1)$) on $X$ such that

\begin{equation}\label{eqn:prop-def-alpha-k}(i)\,\{\alpha_k\}\in H^2_{DR}(X, \, \Z) \hspace{2ex} \mbox{and} \hspace{2ex} (ii)\, ||\alpha_k - k\alpha||_{C^{\infty}}\leq \frac{C}{k^{\frac{1}{b_2}}} \hspace{2ex} \mbox{for all} \,\, k\in S,\end{equation}

\noindent where $C>0$ is a constant independent of $k$ but depending on $X$ and $b_2:=\mbox{dim}_{\R}H^2_{DR}(X, \, \R)$ denotes the second Betti number of $X$. Here, as throughout the paper, the meaning of the symbol $||\,\,\,||_{C^{\infty}}$ is that the stated estimate holds in every $C^l$-norm with a constant $C=C_l>0$ depending on $l\in\N$ (but not on $k$). Thus the notation used in part $(ii)$ of (\ref{eqn:prop-def-alpha-k}) is shorthand for

$$||\alpha_k - k\alpha||_{C^l}\leq \frac{C_l}{k^{\frac{1}{b_2}}} \hspace{2ex} \mbox{for all} \,\, k\in S \,\,\, \mbox{and all} \,\,\, l\in\N.$$

\noindent It is clear that if $\alpha_k=\alpha_k^{2, \, 0} + \alpha_k^{1, \, 1} + \alpha_k^{0, \, 2}$ denotes the splitting of $\alpha_k$ into pure-type components, we have

\begin{equation}\label{eqn:prop-def-alpha-k-bis}(i)\, ||\alpha_k^{1, \, 1} - k\alpha||_{C^{\infty}}\leq \frac{C}{k^{\frac{1}{b_2}}} \hspace{2ex} \mbox{and} \hspace{2ex} (ii)\, ||\alpha_k^{0, \, 2}||_{C^{\infty}}\leq \frac{C}{k^{\frac{1}{b_2}}} \hspace{2ex} \mbox{for all} \,\, k\in S.\end{equation}

\noindent In particular, $\alpha_k^{1, \, 1}$ (as well as $\alpha_k$) comes arbitrarily close to $k\alpha$, while $\alpha_k^{0, \, 2}$ converges to zero in the $C^{\infty}$-topology when $k\rightarrow +\infty$.

 Since the classes $\{\alpha_k\}$ are integral, there exists for every $k\in S$ a Hermitian $C^{\infty}$ (in general not holomorphic) complex line bundle $(L_k, \, h_k)\rightarrow X$ carrying a Hermitian connection $D_k$ of curvature $\frac{i}{2\pi}D_k^2=\alpha_k$. In particular $c_1(L_k)=\{\alpha_k\}$. The complex structure of $X$ induces a splitting of $D_k$ into components $\partial_k$ and $\bar\partial_k$ of respective types $(1, \, 0)$ and $(0, \, 1)$\!\!: 

$$D_k=\partial_k + \bar\partial_k$$

\noindent for which one clearly has

\begin{equation}\label{eqn:del-squared}\bar\partial_k^2=-2\pi i\,\alpha_k^{0, \, 2} \hspace{2ex} \mbox{and} \hspace{2ex} \partial_k\bar\partial_k + \bar\partial_k\partial_k=-2\pi i\,\alpha_k^{1, \, 1}.\end{equation} 
 
\noindent In particular, $\bar\partial_k^2\neq 0$ (i.e. $\bar\partial_k$ is a {\it non-integrable} connection of type $(0, \, 1)$ on $L_k$) if $L_k$ is not holomorphic. This means that, even locally, $L_k$ may admit no holomorphic sections as $\ker\bar\partial_k$ need not contain any non-trivial elements. However, combined with $(ii)$ of (\ref{eqn:prop-def-alpha-k-bis}), the first part of (\ref{eqn:del-squared}) shows that although the line bundle $L_k$ is non-holomorphic, it comes arbitrarily close to being holomorphic as $k\rightarrow +\infty$. The sequence of {\it asymptotically holomorphic} line bundles $(L_k)_{k\in S}$ will play a major role in what follows. 

 (In the classical case when $\{\alpha\}=c_1(L)$ is an integral class, one can of course choose $\alpha_k=k\alpha$ and $L_k=L^k$ is then a genuine holomorphic line bundle in which $\bar\partial_k=\bar\partial$ is {\it integrable}, i.e. $\bar\partial_k^2=0$. However, the {\it non-integrable} case is of concern to us here.)

\subsection{Results of L. Laeng [Lae02]}\label{subsection:Lresults}

One of the main problems considered in [Lae02] was finding a suitable notion of {\it approximately holomorphic} sections of the {\it approximately holomorphic} line bundles $L_k$. Various such notions have been put forward by a number of authors in similar (though not identical) situations (e.g. Donaldson in [Don90] where sections of some vector bundle $E$ lying in the kernel of $D:=\partial_E^{\star} + \bar\partial_E^{\star}$ are considered; Donaldson in [Don96] where Gaussian sections in the local flat model are used to construct global approximately holomorphic sections; Shiffman and Zelditch in [SZ02] where the Boutet de Monvel-Sj\"ostrand [BS76] and the Boutet de Monvel-Guillemin [BG81] methods are used on a suitable circle bundle to construct a pseudodifferential operator $\overline{D}_0$ replacing the standard $\bar\partial_b$ of the integrable case; Ma and Marinescu in [MM02] where the $\mbox{spin}^c$ Dirac operator is used; Borthwick and Uribe in [BU00] where low-lying eigenvalues of a certain rescaled Laplace operator are considered, etc). 

 However, all the works mentioned above share a common, very strong hypothesis that we would like to dispense with\!\!: the curvature form of the line bundle in whose high tensor powers {\it approximately holomorphic} sections are constructed is supposed to be {\it non-degenerate} (i.e. a symplectic form). The Boutet de Monvel-Sj\"ostrand theory [BS76] relies heavily on the non-degeneracy assumption and no version of it in the degenerate case is known. 

 As far as we are aware, the only attempt at tackling the {\it non-integrable, degenerate} case (i.e. where the initial closed $2$-form $\alpha$ does not have rational class and may degenerate at certain points of $X$) was made in the thesis of L. Laeng [Lae02] whose main result we now recall.

 In the setting described in $\S.$\ref{subsection:setting}, the anti-holomorphic Laplace-Beltrami operator acting in bi-degree $(0, \, 0)$ (i.e. on $C^{\infty}$ sections of $L_k$)

$$\Delta_k''=\bar\partial_k^{\star}\bar\partial_k:C^{\infty}(X, \, L_k)\rightarrow C^{\infty}(X, \, L_k)$$

\noindent may have trivial kernel, but the direct sum of its eigenspaces corresponding to small eigenvalues is a natural substitute thereof. Thus Laeng put forward the following space of sections (cf. [Lae02, Propri\'et\'e 4.5, p. 92]).

\begin{Def}\label{Def:space-ah-sections} For every $k\in S\, (=\N^{\star})$, let

$${\cal H}_k:=\bigoplus\limits_{\mu\leq\frac{C}{k^{1+\varepsilon}}}E_{\Delta_k''}^{0, \, 0}(\mu)\subset C^{\infty}(X, \, L_k),$$

\noindent where $E_{\Delta_k''}^{0, \, 0}(\mu)$ stands for the eigenspace of $\Delta_k''$ in bi-degree $(0, \, 0)$ corresponding to the eigenvalue $\mu$, $\varepsilon$ is any constant such that $0<\varepsilon <\frac{2}{b_2}$ (where $b_2=b_2(X)=\mbox{dim}_{\R}H^2(X, \, \R)$ is the second Betti number of $X$) and $C>0$ is an arbitrary constant.

\end{Def}

 The spaces ${\cal H}_k$ are not uniquely or even canonically associated with $\{\alpha\}$ since there is no privileged choice of rational classes $1/k\,\{\alpha_k\}$ in $H^2_{DR}(X, \, \R)$ approximating $\{\alpha\}$. The $C^{\infty}$ sections of $L_k$ belonging to the space ${\cal H}_k$ will be termed {\it approximately holomorphic} sections of $L_k$ by virtue of their satisfying the following obvious property (cf. [Lae02, p.92]).

\begin{Lem}\label{Lem:L2-ah} For every $k\in S$ and every section $s\in{\cal H}_k$, we have

\begin{equation}\label{eqn:L2-ah}||\bar\partial_ks||^2\leq\varepsilon_k||s||^2,\end{equation}

\noindent where $\varepsilon_k:=\frac{C}{k^{1+\varepsilon}}$ and $||\,\cdot\,||$ denotes the $L^2$-norm defined by $\omega$ and $h_k$ on any space $C^{\infty}_{p, \, q}(X, \, L_k)$.

\end{Lem}

\noindent {\it Proof.} If $\langle\langle\,\cdot\, , \,\cdot\,\rangle\rangle$ denotes the $L^2$ scalar product induced by $\omega$ and $h_k$ on any space $C^{\infty}_{p, \, q}(X, \, L_k)$, for every $s\in{\cal H}_k$ we have

$$||\bar\partial_ks||^2 = \langle\langle\Delta_k''s, \, s\rangle\rangle \leq \varepsilon_k||s||^2$$

\noindent by the definition of ${\cal H}_k$.   \hfill $\Box$

\vspace{2ex}

 The main result of Laeng is the following asymptotic growth estimate of the dimension of ${\cal H}_k$ which provides the {\it non-integrable} analogue of the key estimate in Demailly's holomorphic Morse inequalities [Dem85a].

\begin{The}([Th\'eor\`eme 4.4. in [Lae02])\label{The:Hk-asymp-est} We have\!\!:

\begin{equation}\label{eqn:Hk-asymp-est}\liminf\limits_{k\rightarrow +\infty,\, k\in S}\frac{n!}{k^n}\,\mbox{dim}_{\C}{\cal H}_k\geq \int\limits_{X(\alpha, \, \leq 1)}\alpha^n.\end{equation}

\noindent In particular, if we assume $\int_{X(\alpha, \, \leq 1)}\alpha^n>0$ (Demailly's hypothesis in [Dem85a]), then $\mbox{dim}_{\C}{\cal H}_k$ has maximal growth rate (i.e. $O(k^n)$) as $k\rightarrow +\infty$.

\end{The}

 Thus Theorem \ref{The:Hk-asymp-est} shows that the {\it asymptotically holomorphic} line bundles $L_k$ for which one has singled out spaces of {\it approximately holomorphic} sections ${\cal H}_k$ display in the {\it non-integrable} context a property analogous to the familiar notion of {\it big} holomorphic line bundle of the {\it integrable} context. 

 The underlying idea in this approach to Demailly's conjecture on {\it transcendental Morse inequalities} is to manufacture the desired K\"ahler current in the class $\{\alpha\}$ by modifying in the same class a positive current obtained as a limit of currents explicitly constructed from {\it approximately holomorphic} sections of the {\it approximately holomorphic} line bundles $L_k$. If $(\sigma_{k,\, l})_{0\leq l\leq N_k}$ (where $N_k + 1:=\mbox{dim}_{\C}{\cal H}_k$) is an orthonormal basis of ${\cal H}_k$, it is natural to consider the closed $(1, \, 1)$-currents $T_k$ (cohomologous to $\alpha$) on $X$ defined by

\begin{equation}\label{eqn:def1-currents}T_k=\alpha + \frac{i}{2\pi k}\partial\bar\partial\log\sum\limits_{l=0}^{N_k}|\sigma_{k, \, l}|^2_{h_k}, \hspace{3ex}  k\in S\, (=\N^{\star}),\end{equation}

\noindent where at every point $z\in X$ we denote by $|\sigma_{k, \, l}(z)|_{h_k}$ the $h_k$-norm of $\sigma_{k, \, l}(z)\in (L_k)_z$ (see [Lae02, p. 92-100] where an extra $\varepsilon_k$ is inserted to enable the calculations). Since the sections $\sigma_{k, \, l}$ are not holomorphic, each current $T_k$ may have a negative part. Moreover, the $L^2$-estimate (\ref{eqn:L2-ah}) which makes precise the sense in which sections in ${\cal H}_k$ are {\it approximately holomorphic} falls far short of what is needed to control the negative parts of these currents. To obtain such a control, pointwise estimates of the sections $\sigma_{k, \, l}$ and their derivatives of order $\leq 2$ are needed. We need to be able to either produce sections of $L_k$ that are {\it approximately holomorphic} in a sense much stronger than $L^2$, or to get a far better grip on the existing sections making up the space ${\cal H}_k$. This is where the work of [Lae02] comes to an end and new ingredients are needed.

\subsection{Results obtained in this paper}\label{subsection:results-this-paper}

 We propose in this paper a method of constructing global $C^{\infty}$ {\it approximately holomorphic peak sections} of $L_k$ by a careful analysis of the spectrum and (lack of) commutation properties of the anti-holomorphic Laplace-Beltrami operator $\Delta_k''$ of $L_k$. The familiar $L^2$ techniques of the integrable case based on resolutions of the $\bar\partial$-operator are inapplicable in our case where $\bar\partial_k^2\neq 0$ (hence $\bar\partial_k$-exact forms need not even be $\bar\partial_k$-closed). Thus $\bar\partial_k$ is replaced in this approach by $\Delta_k''$ as the main object of study. The starting point is a Weitzenb\"ock-type formula for non-holomorphic vector bundles that essentially appears in [Lae02]. We only very slightly simplify it in section \ref{section:preliminaries}.  

 Having fixed an arbitrary point $x\in X$, the construction of an approximately holomorphic peak section at $x$ will proceed in two stages. First, a local peak section is constructed on a neighbourhood of $x$ as a Gaussian section lying in the kernel of a coupled $\bar\partial$-operator $\bar\partial_{kA}\!\!:=\bar\partial + A_k^{0, \, 1}$ defined by an appropriate $(0, \, 1)$-form $A_k^{0, \, 1}$ coming from the curvature $2$-form $\alpha_k$ of $L_k$. This local construction, performed in $\S.$\ref{subsection:local-model}, has been inspired by Donaldson's approach in [Don96]. Second, we extend the local section $v$ to a global $C^{\infty}$ section $\theta v$ of $L_k$ by multiplying by a cut-off function $\theta$ and then we take the orthogonal projection $s_h$ of this extension onto the space ${\cal H}_k$. This is tantamount to correcting $s\!\!:=\theta v$ to an {\it approximately holomorphic} global section $s_h$ of $L_k$ by subtracting its orthogonal projection $s_{nh}$ onto the orthogonal complement of ${\cal H}_k$ in $C^{\infty}(X, \, L_k)$. We are faced with the challenge of estimating (for example in $L^2$-norm over $X$ and in a stronger norm on a neighbourhood of $x$) the correction $s_{nh}$ in terms of $\bar\partial_ks$. This is done (cf. Proposition \ref{Prop:s_nh-L2-est}) for an arbitrary global section $s\in C^{\infty}(X, L_k)$ in $\S.$\ref{subsection:corrections} which is the heart of the paper. When $\alpha >0$, the global $L^2$-estimate obtained is

$$||s_{nh}||^2\leq\frac{C}{k}\,||\bar\partial_ks||^2, \hspace{3ex} k\gg 1.$$

\noindent A refined local estimate is then obtained in $\S.$\ref{subsection:global-peak} in a neighbourhood of $x$ when $s$ is the extension of the local Gaussian section of $\S.$\ref{subsection:local-model}. 

 We then go on to give two applications of the {\it approximately holomorphic peak sections} in the special case when $\alpha >0$ on $X$ (i.e. $\alpha$ is a K\"ahler metric). This (rather strong positivity) assumption will be removed in future work, but there are already some interesting features in this most basic case. 

 The first application is an {\it approximately holomorphic} Kodaira-type projective embedding theorem for compact K\"ahler manifolds (see Theorem \ref{The:Kod-embed} for a more precise statement).

\begin{The}\label{The:introd-Kod} Let $X$ be a compact K\"ahler (possibly non-projective) manifold. Pick any K\"ahler metric $\alpha$ (possibly of non-rational class) on $X$ and a choice of spaces ${\cal H}_k$ ($k\in S\subset \N^{\star}$) for $\{\alpha\}$. Then the Kodaira-type map 

$$\Phi_k : X\longrightarrow \Proj^{N_k}$$

\noindent associated with ${\cal H}_k$ is everywhere defined and an embedding for $k$ large enough.

\end{The}

 Similar statements have been proved by various authors (e.g. [SZ02]) for symplectic forms of {\it integral} classes. The novelty of our result lies in allowing for the class $\{\alpha\}$ to be {\it transcendental}.

 The other application, an {\it approximately holomorphic} analogue of Tian's {\it almost isometry theorem} [Tia90, Theorem A], can be stated as follows (see Theorem \ref{The:a-isometry} for a more precise statement).

\begin{The}\label{The:introd-Tian} The assumptions are those of Theorem \ref{The:introd-Kod}.

\noindent $(a)$\, The $(1, \, 1)$-current $T_k\!\!:=\alpha + \frac{i}{2\pi k}\partial\bar\partial\log\sum\limits_{l=0}^{N_k}|\sigma_{k, \, l}|^2_{h_k}$ defined in (\ref{eqn:def1-currents}) converges to $\alpha$ in the $C^2$-topology as $k\rightarrow +\infty$.

\noindent $(b)$\, If $\omega_{FS}^{(k)}$ denotes the Fubini-Study metric of $\Proj^{N_k}$ and $\Phi_k$ is the embedding of Theorem \ref{The:introd-Kod}, then $1/k\,\Phi_k^{\star}\omega_{FS}^{(k)}$ converges to $\alpha$ in the $C^2$-topology as $k\rightarrow +\infty$.

\end{The}

 The natural question arising is whether the above $C^2$-norm convergences can be improved to $C^{\infty}$-topology convergences and, moreover, whether there exists an asymptotic expansion for the Bergman kernel function $\sum\limits|\sigma_{k, \, l}|^2_{h_k}$ that would parallel Zelditch's results of [Zel98]. Given the non-degeneracy assumption on $\alpha$, this is likely but the approach would be probably different to the one based on {\it approximately holomorphic peak sections} that we have undertaken here. However, in our view, the present approach has the advantage of lending itself to generalisations when $\alpha$ is allowed to degenerate. This far more general situation that one faces in tackling Demailly's conjecture on {\it transcendental Morse inequalities} will be taken up in future work.

 The global {\it approximately holomorphic peak sections} of $L_k$ constructed in section \ref{section:peak-sections-construction} provide a non-integrable analogue to Tian's {\it holomorphic peak sections} of [Tia90]. While in the case of an ample holomorphic line bundle treated in [Tia90] the Kodaira Embedding Theorem was already available, we show in section \ref{section:proj-embed} that its standard proof can be imitated in the present non-integrable context using our peak sections. Finally, the proof of Theorem \ref{The:introd-Tian} is spelt out in section \ref{section:a-isometry} along the lines of Tian's proof of his holomorphic case result with an emphasis on the handling of the extra derivatives peculiar to the approximately holomorphic case at hand.

\vspace{3ex}

 The author is happy to acknowledge the hospitality of the {\it Instituto de Matem\'aticas de la UNAM} in Cuernavaca (Mexico) where this work was completed and the kind invitation of Professors Alberto Verjovsky and Jos\'e Seade. He is also very grateful to Nick Buchdahl who, from the other side of the world in Australia, had the kindness of reading parts of this paper and making very useful comments.

\section{Preliminaries}\label{section:preliminaries}

 We collect here a few essentially known facts about the Bochner-Kodaira-Nakano and Weitzenb\"ock identities for not necessarily holomorphic vector bundles on possibly non-K\"ahler compact complex manifolds that we arrange in a form that will be useful to us in the subsequent sections. The references are [Gri66], [Dem85b] (the non-K\"ahler case), [Lae02] (the non-holomorphic bundle case) and [Don90] (whose straightforward approach in the almost K\"ahler situation inspired the presentation in \ref{subsection:Weitzenbock} and \ref{subsection:Weitzenbock-special}). This section also fixes the notation for the rest of the paper.

\subsection{Weitzenb\"ock formula\!\!: the non-integrable case}\label{subsection:Weitzenbock}

 Let $(E, \, h_E, \, D_E)\rightarrow (X, \, \omega)$ be a complex Hermitian $C^{\infty}$ vector bundle ($\mbox{rank}_{\C}E=r\geq 1$) equipped with a Hermitian connection $D_E$ over a complex Hermitian manifold ($\mbox{dim}_{\C}X=n$). One denotes

$$D_E=\partial_E + \bar\partial_E  \hspace{2ex} \mbox{and} \hspace{2ex} d=\partial + \bar\partial$$

\noindent the splittings into $(1, \, 0)$ and $(0, \, 1)$-type components of $D_E$ (acting on $E$-valued forms) and respectively $d$ (the Poincar\'e differential operator acting on scalar-valued forms of $X$) w.r.t. the complex structure of $X$. Note that in the general case when $E$ is not holomorphic, we have

$$\bar\partial_E^2= \Theta(E)^{0, \, 2}\neq 0,$$

\noindent where $\Theta(E)^{0, \, 2}$ denotes the $(0, \, 2)$-component of the curvature form of $(E, \, h_E)$. Thus $\bar\partial_E$ is a non-integrable connection of type $(0, \, 1)$ when $E$ is non-holomorphic. 

 For all $p, q = 0, \dots , n$, one considers the Laplace-Beltrami operators

$$\Delta_E'=\partial_E\partial_E^{\star} + \partial_E^{\star}\partial_E:C^{\infty}_{p, \, q}(X, \, E)\rightarrow C^{\infty}_{p, \, q}(X, \, E) \hspace{2ex} (\mbox{hence}\,\, \Delta_E'=[\partial_E, \, \partial_E^{\star}])$$

\noindent and 

$$\Delta_E''=\bar\partial_E\bar\partial_E^{\star} + \bar\partial_E^{\star}\bar\partial_E:C^{\infty}_{p, \, q}(X, \, E)\rightarrow C^{\infty}_{p, \, q}(X, \, E) \hspace{2ex} (\mbox{hence}\,\, \Delta_E''=[\bar\partial_E, \, \bar\partial_E^{\star}])$$

\noindent acting on spaces of $E$-valued  $C^{\infty}$ $(p, \, q)$-forms. (As usual, we use the notation $[A, \, B]:=AB - (-1)^{ab}BA$ for operators of degrees $\deg(A)=a$ and $\deg(B)=b$ on the graded algebra $C^{\infty}_{\cdot, \, \cdot}(X, \, E)$.)

 If $\Lambda=\Lambda_{\omega}$ denotes the formal adjoint (w.r.t. $\omega$) of the multiplication operator $L=L_{\omega}:=\omega\wedge\cdot$, one considers (cf. [Dem85]) the following torsion operator associated with the Hermitian metric $\omega$ on $X$\!\!:

$$\tau:=[\Lambda, \, \partial\omega].$$

\noindent It is clear that $\tau$ is an operator of order zero and bi-degree $(1, \, 0)$. The metric $\omega$ is K\"ahler if and only if $\tau=0$. 

 The Bochner-Kodaira-Nakano identity for holomorphic vector bundles $E$ was extended to the case of a Hermitian (possibly non-K\"ahler) metric $\omega$ by Griffiths in [Gri66] and in a more precise form by Demailly in [Dem85b]. It was later further extended to the case of a possibly non-holomorphic $C^{\infty}$ vector bundle $E$ by Laeng in [Lae02] in the following form. \\

\noindent {\bf Bochner-Kodaira-Nakano identity\!\!:}

\begin{equation}\label{eqn:BKN-nK-nH}\Delta_E''=\Delta_{E, \, \tau}' + [i\Theta(E)^{1, \, 1}, \, \Lambda] + T_{\omega}\end{equation}

\noindent {\it as operators acting on $C^{\infty}_{p, \, q}(X, \, E)$ (for any $p, q=0, \dots , n$), where} 

$$\Delta_{E, \, \tau}':=[\partial_E+\tau, \, \partial_E^{\star} + \tau^{\star}]:C^{\infty}_{p, \, q}(X, \, E)\rightarrow C^{\infty}_{p, \, q}(X, \, E)$$ 

\noindent {\it is the torsion-twisted version of $\Delta_E'$ (and clearly a non-negative formally self-adjoint elliptic operator of order two) and} 

$$T_{\omega}:=[\Lambda, \, [\Lambda, \, \frac{i}{2}\partial\bar\partial\omega]] - [\partial\omega, \, (\partial\omega)^{\star}]:C^{\infty}_{p, \, q}(X, \, E)\rightarrow C^{\infty}_{p, \, q}(X, \, E)$$

\noindent {\it is a zero-order operator that vanishes if $\omega$ is K\"ahler, while $i\Theta(E)^{1, \, 1}$ denotes the $(1, \, 1)$-component of the curvature form of $(E, \, h_E)$ (which is not of type $(1, \, 1)$ when $E$ is not holomorphic).}

\vspace{2ex}

 To derive a Weitzenb\"ock-type formula for $E$-valued $(0, \, q)$-forms in this general setting, one follows the usual route. For every $q=1, \dots , n$, set $\Omega_E^{0, \, q}:=\Lambda^{0, \, q}T^{\star}X\otimes E$. This complex $C^{\infty}$ vector bundle has a natural Hermitian metric induced by $\omega$ and $h_E$. Let

$$\nabla=\nabla' + \nabla''$$

\noindent be the connection on $\Omega_E^{0, \, q}$ induced by $D_E$, $\omega$ and $h_E$, while $\nabla'$, $\nabla''$ denote its respective $(1, \, 0)$ and $(0, \, 1)$-components w.r.t. the complex structure of $X$. Since the vector bundle $\Lambda^{0, \, q}T^{\star}X$ is anti-holomorphic, it has a natural Hermitian connection whose $(1, \, 0)$-component is $\partial$, so $\nabla$ is the unique connection on $\Omega_E^{0, \, q}$ compatible with the metric induced by $h_E$ and $\omega$ and having $\partial_E$ as $(1, \, 0)$-component. Thus

\begin{equation}\label{eqn:nabla'delE}\nabla'=\partial_E.\end{equation}

  It is clear that $E$-valued $(0, \, q)$-forms identify naturally with sections of $\Omega_E^{0, \, q}$ under the obvious isomorphism

\begin{equation}\label{eqn:ident-0q}C^{\infty}_{0, \, q}(X, \, E)\simeq C^{\infty}(X, \, \Omega_E^{0, \, q}).\end{equation} 

 One naturally defines Laplace-Beltrami operators on sections of $\Omega_E^{0, \, q}:$

$$\Box_E'':=\nabla''^{\star}\nabla'':C^{\infty}(X, \, \Omega_E^{0, \, q})\rightarrow C^{\infty}(X, \, \Omega_E^{0, \, q})$$

\noindent and

$$\Box_E':=\nabla'^{\star}\nabla':C^{\infty}(X, \, \Omega_E^{0, \, q})\rightarrow C^{\infty}(X, \, \Omega_E^{0, \, q}),$$

\noindent as well as the torsion-twisted version of the latter by

$$\Box_{E, \, \tau}':=(\nabla'^{\star} + \tau^{\star})(\nabla' + \tau):C^{\infty}(X, \, \Omega_E^{0, \, q})\rightarrow C^{\infty}(X, \, \Omega_E^{0, \, q}).$$

 It is clear that the identifications (\ref{eqn:nabla'delE}) and (\ref{eqn:ident-0q}) give

\begin{equation}\label{eqn:box'delta'}\Box_{E, \, \tau}'=\Delta_{E, \, \tau}'.\end{equation}

\noindent The Weitzenb\"ock formula will follow from a double application of the Bochner-Kodaira-Nakano identity (\ref{eqn:BKN-nK-nH}), first to relate $\Delta_E''$ and $\Delta_{E, \, \tau}'$ acting on $C^{\infty}_{0, \, q}(X, \, E)$ and then to relate $\Box_E''$ and $\Box_{E, \, \tau}'$ acting on $C^{\infty}(X, \, \Omega_E^{0, \, q})$. The link is provided by (\ref{eqn:box'delta'}).

 Indeed, applying (\ref{eqn:BKN-nK-nH}) on sections (i.e. $(0, \, 0)$-forms) of 
$\Omega_E^{0, \, q}$, we get

\begin{equation}\label{eqn:BKN-sections0q}\Box_E''=\Box_{E, \, \tau}' + [i\Theta(\Omega_E^{0, \, q})^{1, \, 1}, \, \Lambda] + T_{\omega} \hspace{3ex} \mbox{on} \hspace{1ex} C^{\infty}(X, \, \Omega_E^{0, \, q}).\end{equation}

\noindent Now using identification (\ref{eqn:ident-0q}) and identity (\ref{eqn:box'delta'}) and putting together the two instances (\ref{eqn:BKN-nK-nH}) and (\ref{eqn:BKN-sections0q}) of the Bochner-Kodaira-Nakano identity, we get\!\!: \\

\noindent {\bf Weitzenb\"ock formula for $E$-valued $(0, \, q)$-forms\!\!:}

\begin{equation}\label{eqn:Weitzenbock-0q}\Delta_E''=\Box_E'' - [i\Theta(\Omega_E^{0, \, q})^{1, \, 1}, \, \Lambda] + [i\Theta(E)^{1, \, 1}, \, \Lambda]\end{equation}

\noindent {\it on} $C^{\infty}_{0, \, q}(X, \, E)\simeq C^{\infty}(X, \, \Omega_E^{0, \, q})$.

\subsection{Weitzenb\"ock formula in Laeng's special setting}\label{subsection:Weitzenbock-special}

 We now specialise the discussion in subsection \ref{subsection:Weitzenbock} to the situation described in subsection \ref{subsection:setting}. This was already done in [Lae02]. We only very slightly simplify the formulae. The complex manifold $X$ is supposed to be compact.

 If we choose $E=(L_k, \, h_k, \, D_k)\rightarrow (X, \, \omega)$ (so $\mbox{rank}_{\C}L_k=1$), the curvature form is $\frac{i}{2\pi}\Theta(L_k)=\alpha_k$ and the Laplace-Beltrami operators $\Delta_E''$, $\Delta_E'$  and $\Delta_{E, \, \tau}'$ become the following operators acting on $C^{\infty}_{p, \, q}(X, \, L_k)$\!:

$$\Delta_k''=[\bar\partial_k, \, \bar\partial_k^{\star}], \hspace{3ex} \Delta_k'=[\partial_k, \, \partial_k^{\star}] \hspace{3ex} \mbox{and} \hspace{3ex} \Delta_{k, \, \tau}'=[\partial_k+\tau, \, \partial_k^{\star} + \tau^{\star}].$$

\noindent If the connection on $\Omega_k^{0, \, q}:=\Lambda^{0, \, q}T^{\star}X\otimes L_k$ induced by $D_k$ and $\omega$ is denoted (when splitting into $(1, \, 0)$ and $(0, \, 1)$-components) 

$$\nabla_k=\nabla_k' + \nabla_k'',$$ 

\noindent we have $\nabla_k'=\partial_k$ (cf. \ref{eqn:nabla'delE}) on $C^{\infty}_{0, \, q}(X, \, L_k)\simeq C^{\infty}(X, \, \Omega_k^{0, \, q})$. The Laplace-Beltrami operators induced by $\nabla_k$ on $C^{\infty}(X, \, \Omega_k^{0, \, q})$ read\!\!:

$$\Box_k'':={\nabla_k''}^{\star}\nabla_k'', \hspace{2ex} \Box_k':={\nabla_k'}^{\star}\nabla_k' \hspace{3ex} \mbox{and} \hspace{3ex} \Box_{k, \, \tau}':=({\nabla_k'}^{\star} + \tau^{\star})(\nabla_k' + \tau).$$

\noindent Thus by (\ref{eqn:box'delta'}) we have $\Box_k'=\Delta_{k, \, \tau}'$ and the Weitzenb\"ock formula (\ref{eqn:Weitzenbock-0q}) for $L_k$-valued $(0, \, q)$-forms reduces to

\begin{equation}\label{eqn:Weitzenbock_k-rough}\Delta_k''=\Box_k'' + \Lambda(i\Theta(\Omega_k^{0, \, q})^{1, \, 1}) - \Lambda(i\Theta(L_k)^{1, \, 1})\end{equation}

\noindent on $C^{\infty}_{0, \, q}(X, \, L_k)\simeq C^{\infty}(X, \, \Omega_k^{0, \, q})$ because $\Lambda$ (which is of type $(-1, \, -1)$) acts trivially on $(0, \, q)$-forms for bi-degree reasons.

 Now $i\Theta(L_k)^{1, \, 1}=\alpha_k^{1, \, 1}$ and $i\Theta(\Omega_k^{0, \, q}) = i\Theta(\Lambda^{0, \, q}T^{\star}X)\otimes \mbox{Id}_{L_k} + \mbox{Id}_{\Lambda^{0, \, q}T^{\star}X}\otimes i\Theta(L_k)$. Since $i\Theta(\Lambda^{0, \, q}T^{\star}X)$ is of type $(1, \, 1)$ (due to $\Lambda^{0, \, q}T^{\star}X$ being anti-holomorphic), passing to $(1, \, 1)$-components in the last identity we are left with $i\Theta(\Omega_k^{0, \, q})^{1, \, 1} = i\Theta(\Lambda^{0, \, q}T^{\star}X)\otimes \mbox{Id}_{L_k} + \mbox{Id}_{\Lambda^{0, \, q}T^{\star}X}\otimes \alpha_k^{1, \, 1}$. The Weitzenb\"ock formula (\ref{eqn:Weitzenbock_k-rough}) translates to

\begin{equation}\label{eqn:Weitzenbock_k-rough-alpha}\Delta_k''=\Box_k'' + \Lambda(\mbox{Id}_{\Lambda^{0, \, q}T^{\star}X}\otimes\alpha_k^{1, \, 1}) - \Lambda(\alpha_k^{1, \, 1}) + \Lambda(i\Theta(\Lambda^{0, \, q}T^{\star}X)\otimes\mbox{Id}_{L_k}).\end{equation}

 Set $R:=\Lambda(i\Theta(\Lambda^{0, \, q}T^{\star}X)\otimes\mbox{Id}_{L_k})$, a zero-order operator independent of $k$. In order to better exploit the fact that $\alpha_k^{1, \, 1}$ is close to $k\alpha$ for $k$ large, we write $\alpha_k^{1, \, 1}=(\alpha_k^{1, \, 1}-k\alpha) + k\alpha$ and (\ref{eqn:Weitzenbock_k-rough-alpha}) translates to the following \\

\noindent {\bf Weitzenb\"ock formula for $L_k$-valued $(0, \, q)$-forms\!\!:}

\begin{equation}\label{eqn:Weitzenbock_k-rough-short}\Delta_k''=\Box_k'' + kV + (R_{\alpha, \, k} + R)\end{equation}

\noindent {\it on $C^{\infty}_{0, \, q}(X, \, L_k)\simeq C^{\infty}(X, \, \Lambda^{0, \, q}T^{\star}X\otimes L_k)$, where we have denoted}

\begin{equation}\label{eqn:Weitzenbock_k-not1}V:=\Lambda(\mbox{Id}_{\Lambda^{0, \, q}T^{\star}X}\otimes\alpha) - \Lambda(\alpha),\end{equation}

\noindent {\it a zero-order operator, independent of $k$,}

\begin{equation}\label{eqn:Weitzenbock_k-not2}R_{\alpha, \, k}:=\Lambda(\mbox{Id}_{\Lambda^{0, \, q}T^{\star}X}\otimes(\alpha_k^{1, \, 1}-k\alpha)) - \Lambda(\alpha_k^{1, \, 1}-k\alpha),\end{equation}

\noindent {\it a zero-order operator that depends on $k$ but tends to zero as $k\rightarrow +\infty$,}

\begin{equation}\label{eqn:Weitzenbock_k-not3}R:=\Lambda(i\Theta(\Lambda^{0, \, q}T^{\star}X)\otimes\mbox{Id}_{L_k}),\end{equation}

\noindent {\it a zero-order operator, independent of $k$.}

\vspace{2ex}

\noindent {\bf Calculation of $V$.} Let $\lambda_1\leq \ldots \leq \lambda_n$ stand for the eigenvalues of $\alpha$ w.r.t. $\omega$ ordered non-decreasingly. Let us fix an arbitrary point $x\in X$. We can find local holomorphic coordinates $z_1, \dots , z_n$ about $x$ such that

$$\omega(x)=i\sum\limits_{j=1}^ndz_j\wedge d\bar{z}_j  \hspace{3ex} \mbox{and} \hspace{3ex} \alpha(x)=i\sum\limits_{j=1}^n\lambda_j(x)dz_j\wedge d\bar{z}_j.$$

 Fix some $0\leq q\leq n$ and let $u\in C^{\infty}_{0, \, q}(X, \, L_k)$ be arbitrary. Then in a neighbourhood of $x$ we can write

$$u=\sum\limits_{|J|=q}u_J\, d\bar{z}_J, \hspace{3ex} \mbox{for some smooth functions} \hspace{1ex} u_J,$$

\noindent where $d\bar{z}_J:=d\bar{z}_{j_1}\wedge\ldots\wedge\bar{z}_{j_q}$ for every $J=(j_1< \ldots < j_q)$. The first identity in the following formula is known to hold at $x$ (cf. e.g. [Dem97, VI-$\S.5.2.$] or [Lae02, II. 2.6, p. 68])\!\!:

\begin{eqnarray}\nonumber\langle[\alpha, \, \Lambda]u, \, u\rangle & = & \sum\limits_{|J|=q}(\sum\limits_{j\in J}\lambda_j - \sum\limits_{l=1}^n\lambda_l)\,|u_J|^2\\
\nonumber  & = & \sum\limits_{|J|=q}(\sum\limits_{j\in J}\lambda_j)\,|u_J|^2 - (\sum\limits_{l=1}^n\lambda_l)\sum\limits_{|J|=q}|u_J|^2.\end{eqnarray}

\noindent Since $[\alpha, \, \Lambda]u = -(\Lambda\alpha)u$ (because $\Lambda u=0$ for bi-degree reasons) and since $\sum\limits_{l=1}^n\lambda_l = \mbox{Tr}_{\omega}\alpha$ is the trace of $\alpha$ w.r.t. $\omega$, we have obtained the formula

\begin{equation}\label{eqn:Lambda-alpha-u}\langle(\Lambda\alpha)u, \, u\rangle = (\mbox{Tr}_{\omega}\alpha)\,|u|^2 - \sum\limits_{|J|=q}(\sum\limits_{j\in J}\lambda_j)\,|u_J|^2\end{equation}

\noindent at $x$ for any $L_k$-valued $(0, \, q)$-form $u$.

 We can now regard $u$ as a $(0, \, 0)$-form with values in $\Omega_k^{0, \, q}$. The above formula (\ref{eqn:Lambda-alpha-u}) applied on $\Omega_k^{0, \, q}$-valued $(0, \, 0)$-forms reads\!\!:

\begin{equation}\label{eqn:Lambda-alpha-u-00}\langle\Lambda(\mbox{Id}_{\Lambda^{0, \, q}T^{\star}X}\otimes\alpha)u, \, u\rangle = (\mbox{Tr}_{\omega}\alpha)\,|u|^2\end{equation}

\noindent at $x$ for any section $u$ of $\Omega_k^{0, \, q}$.

 Combining (\ref{eqn:Lambda-alpha-u}) and (\ref{eqn:Lambda-alpha-u-00}), we obtain the formula

\begin{equation}\label{eqn:V-formula}\langle Vu, \, u\rangle = \sum\limits_{|J|=q}(\sum\limits_{j\in J}\lambda_j)\,|u_J|^2\end{equation}

\noindent at every point $x\in X$ and for every $u=\sum\limits_{|J|=q}u_J\, d\bar{z}_J\in C^{\infty}_{0, \, q}(X, \, L_k).$

 Consequently the Weitzenb\"ock formula (\ref{eqn:Weitzenbock_k-rough-short}) translates to the following \\

\noindent {\bf Explicit Weitzenb\"ock formula for $C^{\infty}_{0, \, q}(X, \, L_k)\simeq C^{\infty}(X,\, \Lambda^{0, \, q}T^{\star}X\otimes L_k)$\!\!:}

\begin{equation}\label{eqn:explicit-Weitzenbock-0q}\langle \Delta_k''u, \, u\rangle = \langle \Box_k''u, \, u\rangle + k \sum\limits_{|J|=q}(\sum\limits_{j\in J}\lambda_j)\,|u_J|^2 + \langle(R_{\alpha, \, k} + R)u, \, u\rangle\end{equation}

\noindent {\it at every point $x\in X$ and for every $u=\sum\limits_{|J|=q}u_J\, d\bar{z}_J\in C^{\infty}_{0, \, q}(X, \, L_k).$}

 When $q=1$, we get\!\!:

\begin{equation}\label{eqn:explicit-Weitzenbock-01}\langle \Delta_k''u, \, u\rangle = \langle \Box_k''u, \, u\rangle + k \sum\limits_{j=1}^n\lambda_j\,|u_j|^2 + \langle(R_{\alpha, \, k} + R)u, \, u\rangle.\end{equation}

 Note that $\Box_k''$ is a non-negative operator (i.e. $\langle\langle \Box_k''u, \, u\rangle\rangle\geq 0$ for every $u$) and $R_{\alpha, \, k} + R$ is an operator of order $0$, bounded independently of $k$ (see (\ref{eqn:Weitzenbock_k-not2}) and (\ref{eqn:Weitzenbock_k-not3})). Hence in the special case when $\alpha$ is supposed to be {\it positive-definite} (i.e. $\lambda_j(z)>0$ for all $j=1, \dots , n$ and all $z\in X$), we immediately get

\begin{Cor}\label{Cor:W-pos-cons} Suppose that $\alpha>0$ at every point of $X$. Then the Laplace-Beltrami operator $\Delta_k'':C^{\infty}_{0, \, 1}(X, \, L_k)\rightarrow C^{\infty}_{0, \, 1}(X, \, L_k)$ satisfies

\begin{equation}\label{eqn:delta''-lbound-01}\Delta_k''\geq\delta_0k>0, \hspace{3ex}  \mbox{for all}\hspace{1ex} k\gg 1,\end{equation}

\noindent (i.e. $\langle\langle\Delta_k''u, \, u\rangle\rangle \geq \delta_0k\,||u||^2$ for all $u\in C^{\infty}_{0, \, 1}(X, \, L_k)$ and all $k\gg 1$), where $\delta_0>0$ is any constant for which $\alpha\geq 2\delta_0\,\omega$ on $X$.

\end{Cor}

\subsection{The spectral gap in bi-degree $(0, \, 0)$ when $\alpha >0$}\label{subsection:spectral-gap}

 Before developing the new arguments, we explain in this subsection how a result of Laeng [Lae02, 4.2.2, p. 90-91] gives additional information on the spectrum of $\Delta_k'':C^{\infty}(X, \, L_k)\rightarrow C^{\infty}(X, \, L_k)$ in a special case. This will be needed in the next section.

 Since $\bar\partial_k^2\neq 0$, $\bar\partial_k$ does not commute with $\Delta_k''$. Indeed, the commutation defect is easily seen to be 

$$(\Delta_k''\bar\partial_k - \bar\partial_k\Delta_k'')\,s = \bar\partial_k^{\star}\bar\partial_k^2s,  \hspace{2ex} s\in C^{\infty}(X, \, L_k),$$

\noindent and it was given the following $L^2$-norm estimate in [Lae02, p.90]\!\!:

\begin{equation}\label{eqn:comm-defect-Laurent}||\bar\partial_k^{\star}\bar\partial_k^2s||^2 \leq \frac{C}{k^{2/b_2}}\,\bigg(k||s||^2 + ||\bar\partial_ks||^2\bigg), \hspace{2ex} s\in C^{\infty}(X, \, L_k).\end{equation}

 The lack of commutation between $\bar\partial_k$ and $\Delta_k''$ means that an eigenvalue $\lambda$ of $\Delta_k''$ in bi-degree $(0, \, 0)$ need not be an eigenvalue of $\Delta_k''$ in bi-degree $(0, \, 1)$. In particular, $\bar\partial_k$ need not define an injection of the eigenspace $E_{\Delta_k''}^{0, \, 0}(\lambda)$ into $E_{\Delta_k''}^{0, \, 1}(\lambda)$ when $\lambda\neq 0$ as is the case when $\bar\partial_k$ is integrable. However, it was shown in [Lae02] that a part of the spectrum in bi-degree $(0, \, 0)$ injects into an appropriate part of the spectrum in bi-degree $(0, \, 1)$. For any $(p, \, q)$, let $E_{\Delta_k''}^{p, \, q}(\mu)$ stand for the eigenspace of $\Delta_k'':C^{\infty}_{p, \, q}(X, \, L_k) \rightarrow C^{\infty}_{p, \, q}(X, \, L_k)$ corresponding to the eigenvalue $\mu$ (with the understanding that $E_{\Delta_k''}^{p, \, q}(\mu) = \{0\}$ if $\mu$ is not an actual eigenvalue). For any $0<\lambda_1<\lambda_2$ and any $\varepsilon_0'>0$, considering the intervals $I=(\lambda_1, \, \lambda_2]$ and $J=[0, \, \lambda_2 + \varepsilon_0']$ of $\R$ and setting

$$E_I^{0, \, 0}:=\bigoplus\limits_{\lambda\in I}E_{\Delta_k''}^{0, \, 0}(\lambda)\subset C^{\infty}(X, \, L_k) \hspace{2ex} \mbox{and} \hspace{2ex} E_J^{0, \, 1}:=\bigoplus\limits_{\mu\in J}E_{\Delta_k''}^{0, \, 1}(\mu)\subset C^{\infty}_{0, \, 1}(X, \, L_k)$$

\noindent (the present notation differs from that of [Lae02] where $\lambda$ and $\mu$ stood for eigenvalues of $\frac{1}{k}\Delta_k''$ rather than $\Delta_k''$), it was shown in [Lae02] that the map

$$\Pi_J\circ\bar\partial_k : E_I^{0, \, 0}\rightarrow E_J^{0, \, 1},$$

\noindent which is the composition of $\bar\partial_k$ with the orthogonal projection $\Pi_J$ of $C^{\infty}_{0, \, 1}(X, \, L_k)$ onto $E_J^{0, \, 1}$, is injective if appropriate choices of $\lambda_1, \lambda_2$ and $\varepsilon_0'$ are made. The reasoning proceeds in [Lae02] by contradiction in the following way\!\!: if for some $s\in E_I^{0, \, 0}\setminus\{0\}$ we had $\Pi_J(\bar\partial_ks)=0$, then $\bar\partial_ks\in\oplus_{\mu >\lambda_2+\varepsilon_0'}E_{\Delta_k''}^{0, \, 1}(\mu)$, hence

\begin{equation}\label{eqn:est1}||\Delta_k''\bar\partial_ks||\geq (\lambda_2 + \varepsilon_0')\,||\bar\partial_ks||,\end{equation}

\noindent while, on the other hand, we have

$$||\bar\partial_ks||^2 = \langle\langle\Delta_k''s, \, s\rangle\rangle\geq \lambda_1||s||^2$$

\noindent which, for this particular $s$, transforms (\ref{eqn:comm-defect-Laurent}) to

\begin{equation}\label{eqn:comm-defect-Laurent-special}||\bar\partial_k^{\star}\bar\partial_k^2s||^2 \leq \frac{C}{k^{2/b_2}}\,\bigg(1 + \frac{k}{\lambda_1}\bigg) ||\bar\partial_ks||^2.\end{equation}

\noindent Writing now $\Delta_k''\bar\partial_ks = \bar\partial_k\Delta_k''s + \bar\partial_k^{\star}\bar\partial_k^2s$, we get

\begin{equation}\label{eqn:est2}||\Delta_k''\bar\partial_ks||\leq ||\bar\partial_k\Delta_k''s|| + ||\bar\partial_k^{\star}\bar\partial_k^2s|| \leq \bigg(\lambda_2 + Ck^{-1/b_2}\frac{\sqrt{k}}{\sqrt{\lambda_1}}\bigg)\,||\bar\partial_ks||\end{equation}

\noindent if we choose $\lambda_1<k$. (Indeed, $\Delta_k''s\in E_I^{0, \, 0}$ for $s\in E_I^{0, \, 0}$, so $||\bar\partial_k\Delta_k''s||\leq \lambda_2 ||\bar\partial_ks||$. Meanwhile, $1 + k/\lambda_1 < 2\,k/\lambda_1$ if we choose $\lambda_1 < k$; the factor $2$ can be absorbed in the constant $C$ in (\ref{eqn:comm-defect-Laurent-special}).) Now putting (\ref{eqn:est1}) and (\ref{eqn:est2}) together and using the fact that $\bar\partial_ks\neq 0$ (because $\bar\partial_ks = 0 \,\, \Leftrightarrow \,\, \langle\langle\Delta_k''s, \, s\rangle\rangle = 0 \,\, \Leftrightarrow \,\, s\in E_{\Delta_k''}^{0, \, 0}(0)$ which is ruled out by the choices made above), we get

$$\varepsilon_0' \leq Ck^{-1/b_2}\frac{\sqrt{k}}{\sqrt{\lambda_1}}.$$

\noindent Now fix an arbitrary $\varepsilon_0''>0$ independent of $k$ and choose $\varepsilon_0'=\varepsilon_0'(k):=\varepsilon_0''\,k$. The above inequality translates to

\begin{equation}\label{eqn:varepsilon_0'}\varepsilon_0'' \leq Ck^{-1-1/b_2}\frac{\sqrt{k}}{\sqrt{\lambda_1}}.\end{equation}

\noindent Choose moreover $\lambda_1=\lambda_1(k) = C\,k^{-1-\varepsilon}$ with an arbitrary $0 < \varepsilon < 2/b_2$. Then (\ref{eqn:varepsilon_0'}) reads

$$\varepsilon_0'' \leq \frac{C}{k^{\frac{1}{2}(2/b_2 - \varepsilon)}}$$

 \noindent which is impossible if $k$ is large enough since $2/b_2 - \varepsilon>0$ with the above choices. Thus no $s$ as above exists, which means that the map $\Pi_J\circ\bar\partial_k$ is injective, when $k\gg 1$. Hence

\begin{equation}\label{eqn:E-J-E-I-est}\mbox{dim}_{\C}E_J^{0, \, 1}\geq \mbox{dim}_{\C}E_I^{0, \, 0}\geq 0  \hspace{3ex} \mbox{if} \hspace{1ex} k\gg 1,\end{equation}

\noindent with the above choices of $\lambda_1=\lambda_1(k)>0$, $\varepsilon_0'=\varepsilon_0'(k)$ (for any fixed $\varepsilon_0''>0$ independent of $k$) and every $\lambda_2>0$ (that may or may not depend on $k$). This inequality (\ref{eqn:E-J-E-I-est}) was used in [Lae02] in the proof of Theorem \ref{The:Hk-asymp-est}.

 Alternatively, we can use it in the following way. If we assume that $\alpha >0$ on $X$, Corollary \ref{Cor:W-pos-cons} implies that all the eigenvalues of $\Delta_k''$ in bi-degree $(0, \, 1)$ are $\geq \delta_0k$ for $k$ large enough. Hence $E_J^{0, \, 1}=\{0\}$ if we choose $\lambda_2 >0$, $\varepsilon_0''>0$ and $\varepsilon_0'=\varepsilon_0''k>0$ such that $\lambda_2 + \varepsilon_0''k < \delta_0k$. (We can choose, for example, $0 < \varepsilon_0'' < \delta_0$ and $\lambda_2=\lambda_2(k):=(\delta_0 - \varepsilon_0)k$ where $0 < \varepsilon_0'' < \varepsilon_0 < \delta_0$ with $\varepsilon_0$ independent of $k$). In this case, (\ref{eqn:E-J-E-I-est}) implies that $E_I^{0, \, 0}=\{0\}$ which means that $\Delta_k''$ acting in bi-degree $(0, \, 0)$ has no eigenvalues in the interval $I=(\lambda_1, \, \lambda_2]$. We thus get the following

\begin{Cor}\label{Cor:spectral-gap} Suppose that $\alpha >0$ at every point of $X$. Let $\mbox{Spec}^{0, \, 0}(\Delta_k'')$ denote the set of eigenvalues of $\Delta_k'':C^{\infty}(X, \, L_k)\rightarrow C^{\infty}(X, \, L_k)$. Then, for any constants $C>0$, $0<\varepsilon < \frac{2}{b_2}$, $\delta_0>0$ such that $\alpha \geq 2\delta_0\omega$ on $X$ and any $0 < \varepsilon_0 < \delta_0$,  we have

\begin{equation}\label{eqn:spectral-gap}\mbox{Spec}^{0, \, 0}(\Delta_k'') \cap \bigg(\frac{C}{k^{1+\varepsilon}}, \, (\delta_0-\varepsilon_0)k\bigg] = \emptyset   \hspace{3ex} \mbox{if} \hspace{1ex} k \hspace{1ex} \mbox{is large enough}.\end{equation}

\end{Cor}

\section{Construction of peak sections}\label{section:peak-sections-construction}

 To control the negative part of $T_k$ (defined in (\ref{eqn:def1-currents})), the first estimate we need is a pointwise lower bound for the Bergman kernel function

\begin{equation}\label{eqn:ak-def}a_k(x):=\sum\limits_{l=1}^{N_k}|\sigma_{k, \, l}(x)|^2=\sup\limits_{\sigma\in\bar{B}_k(1)}|\sigma(x)|^2, \hspace{3ex} x\in X,\end{equation}

\noindent where $\bar{B}_k(1)\subset{\cal H}_k$ denotes the closed unit ball of ${\cal H}_k$ and the latter identity follows by considering the evaluation linear map ${\cal H}_k\ni\sigma\mapsto\sigma(x)\in\C$ (whose squared $L^2$-norm equals both $a_k(x)$ and the right-hand term of (\ref{eqn:ak-def})) at any given point $x\in X$. The function $a_k$ features in the denominators of expressions that make up $T_k$ when the derivatives of $\log a_k$ are calculated.   

 In the integrable case, the only known way of obtaining a pointwise lower bound for expressions similar to $a_k$ (in which the $\sigma_{k, \, l}$'s are genuine holomorphic functions) is Demailly's method of [Dem92, Proposition 3.1] consisting in an application of the Ohsawa-Takegoshi $L^2$-extension theorem at the given point $x$\!\!: a global holomorphic function $\sigma$ exists with prescribed value at $x$ ($\sigma(x)=e^{\varphi_k(x)}$ is chosen if $\varphi_k$ denotes the psh local weight of $h_k$ near $x$) and with $L^2$-norm under control. After normalisation, $\sigma$ can be made to fit in the unit ball $\bar{B}_k(1)$ and $Const\, |\sigma(x)|^2$ provides an explicit lower bound for $a_k(x)$ in terms of $\varphi_k(x)$. By construction, this section $\sigma$ is ``not too small'' at $x$.

 Genuine {\it peak sections} in the sense of $L^2$-norms were constructed by Tian for high tensor powers of a positive holomorphic line bundle $L$ in [Tia90] (where the term {\it peak section} was used). H\"ormander's $L^2$-estimates were employed there to produce global holomorphic sections of $L^k$ whose $L^2$-norms get increasingly concentrated on increasingly smaller balls about a given point $x$ as $k\rightarrow +\infty$. 

 Both these methods completely break down in our case\!\!: no {\it non-integrable} analogues of the Ohsawa-Takegoshi and H\"ormander's theorems are known and no positivity assumption is made on the possibly {\it degenerate} form $\alpha$.

 As a substitute for these (by now) classical techniques, we propose a method of constructing global $C^{\infty}$ sections of $L_k$ belonging to the space ${\cal H}_k$ (hence {\it approximately holomorphic}) that {\it peak} at an arbitrary point $x\in X$ given beforehand and whose $L^2$-norms are under control. The starting point of the construction, consisting in the use of appropriate locally defined Gaussian sections, has been inspired by Donaldson's approach in [Don96].

\subsection{The local model}\label{subsection:local-model} 

 The notation is that of \ref{subsection:setting}. Fix an arbitrary point $x\in X$. Since $\alpha$ is a closed form, one can find a $C^{\infty}$ function $\varphi : U\rightarrow \R$ on an open neighbourhood $U$ of $x$ such that

$$\alpha=\frac{i}{2\pi}\partial\bar\partial\varphi \hspace{3ex}  \mbox{on}\,\, U.$$

\noindent It follows that

$$\alpha=d\bigg(\frac{i}{2\pi}\bar\partial\varphi\bigg) = d\bigg(-\frac{i}{2\pi}\partial\varphi\bigg),\hspace{3ex} \mbox{so} \hspace{1ex} \alpha= d\bigg(\frac{i}{4\pi}(\bar\partial\varphi - \partial\varphi)\bigg) = idA,$$

\noindent where we have denoted $A:= \frac{1}{4\pi}(\bar\partial\varphi - \partial\varphi)$, a $1$-form on $U$. Since the (possibly) non-rational class $\{\alpha\}$ need not correspond to a line bundle on $X$, $A$ need not be associated with a connection of a holomorphic line bundle with curvature form $\alpha$, but can be thought of as mimicking such a connection. We have

$$A^{1, \, 0} = - \frac{1}{4\pi}\partial\varphi \hspace{3ex} \mbox{and} \hspace{3ex} A^{0, \, 1}= \frac{1}{4\pi}\bar\partial\varphi   \hspace{3ex} \mbox{on}\,\, U.$$

 On the other hand, shrinking $U$ about $x$ if necessary, $L_k$ may be assumed trivial on $U$. Let $2\pi A_k$ be the $C^{\infty}$ $1$-form representing the connection $D_k$ of $L_k$ (known to have curvature form $\alpha_k$) in this local trivialisation $L_{k|U}\stackrel{\theta_k}{\simeq}U\times\C$\!\!:

$$D_k=d + 2\pi A_k \hspace{3ex} \mbox{on}\,\, U.$$ 

\noindent It follows that $\alpha_k = \frac{i}{2\pi}\,D_k^2=idA_k$ (hence $\alpha_k^{0, \, 2} = i\bar\partial A_k^{0, \, 1}$) on $U$ and  

\begin{equation}\label{eqn:dbar-k-local}\bar\partial_k=\bar\partial + 2\pi A_k^{0, \, 1} \hspace{3ex} \mbox{on}\,\, U.\end{equation}

 Since $||k\alpha - \alpha_k||_{C^{\infty}} = ||d(kA - A_k)||_{C^{\infty}}\leq\frac{C}{k^{1/b_2}}$ by (\ref{eqn:prop-def-alpha-k}), we can choose $\varphi$ and $\theta_k$ such that $||kA - A_k||_{C^{\infty}}\leq \frac{C}{k^{1/b_2}}$, which amounts to

\begin{equation}\label{eqn:A-Ak-types-close}||kA^{1, \, 0} - A_k^{1, \, 0}||_{C^{\infty}}\leq \frac{C}{k^{1/b_2}} \hspace{3ex} \mbox{and} \hspace{3ex} ||kA^{0, \, 1} - A_k^{0, \, 1}||_{C^{\infty}}\leq \frac{C}{k^{1/b_2}}.\end{equation}

 To exploit the proximity of $kA^{0, \, 1}$ to $A_k^{0, \, 1}$, we define the following coupled $\bar\partial$-operators on $U$ (which unlike $\bar\partial_k$ do not globalise to the whole of $X$ since the class $\{\alpha\}$ need not be rational).

\begin{Def}\label{Def:dbar-A-local} On the $L_k$-trivialising open subset $U\subset X$, set

\begin{equation}\label{eqn:dbar-A-local}\bar\partial_A=\bar\partial + 2\pi A^{0, \, 1}
 \hspace{3ex} \mbox{and} \hspace{3ex}  \bar\partial_{kA}=\bar\partial + 2\pi kA^{0, \, 1}, \hspace{3ex} k\in\N^{\star}.\end{equation}

\end{Def}

 Thanks to (\ref{eqn:A-Ak-types-close}), $\bar\partial_{kA}$ is {\it close} to $\bar\partial_k$ as will be seen shortly.

 Now choose local holomorphic coordinates $z_1, \dots , z_n$ centred at $x$ (and defined on $U$) such that

\begin{equation}\label{eqn:choice-local-coord}\omega(x)=\frac{i}{2\pi}\,\sum\limits_{j=1}^ndz_j\wedge d\bar{z}_j  \hspace{3ex} \mbox{and} \hspace{3ex} \alpha(x)=\frac{i}{2\pi}\,\sum\limits_{j=1}^n\lambda_j(x)\,dz_j\wedge d\bar{z}_j,\end{equation}

\noindent where $\lambda_1(x)\leq \ldots \leq \lambda_n(x)$ are the eigenvalues of $\alpha$ w.r.t. $\omega$ at $x$ (cf. notation in $\S.$\ref{subsection:Weitzenbock-special}). It is clear that the $C^{\infty}$ function $\varphi : U \rightarrow \R$ with the property $i/2\pi\,\partial\bar\partial\varphi=\alpha$ can be chosen such that

\begin{equation}\label{eqn:phi-choice}\varphi(z) = \sum\limits_{j=1}^n\lambda_j(x)\,|z_j|^2 + {\cal O}(|z|^3), \hspace{3ex} z\in U.\end{equation}

 Consider now the following $C^{\infty}$ function $u:U\rightarrow\R$,

\begin{equation}\label{eqn:loc-section-def}\nonumber u(z):=e^{-\frac{1}{2}\,\varphi(z)}= e^{-\frac{1}{2}\sum\limits_{j=1}^n\lambda_j(x)\,|z_j|^2 + {\cal O}(|z|^3)}, \hspace{3ex} z\in U.\end{equation}

\noindent We have\!\!: \,\,$\bar\partial u = -e^{-\frac{1}{2}\,\varphi(z)}\, (\frac{1}{2}\, \bar\partial\varphi) = - 2\pi A^{0, \, 1}\wedge u$, so

\begin{equation}\label{eqn:coupled-dbar-vanishes}\nonumber\bar\partial_Au=0  \hspace{3ex} \mbox{on}\,\,\, U.\end{equation}

\noindent Similarly, for 

\begin{equation}\label{eqn:loc-section-k-def}u^k(z)=e^{-\frac{k}{2}\,\varphi(z)} = e^{-\frac{k}{2}\sum\limits_{j=1}^n\lambda_j(x)\,|z_j|^2 + {\cal O}(|z|^3)}, \hspace{3ex} z\in U,k\in\N^{\star},\end{equation}

\noindent we have

\begin{equation}\label{eqn:coupled-dbar-k-vanishes}\nonumber\bar\partial_{kA}(u^k)=0  \hspace{3ex} \mbox{on}\,\,\, U.\end{equation}

\noindent We can now easily go from $\bar\partial_{kA}$ to $\bar\partial_k:$

\begin{eqnarray}\nonumber||\bar\partial_k(u^k)||_{C^{\infty}} & \leq & ||\bar\partial_k(u^k) -\bar\partial_{kA}(u^k)||_{C^{\infty}} + ||\bar\partial_{kA}(u^k)||_{C^{\infty}} = ||2\pi\,(A_k^{0, \, 1} -kA^{0, \, 1})u^k||_{C^{\infty}} \\
\nonumber  & \leq & 2\pi\,||A_k^{0, \, 1} -kA^{0, \, 1}||_{C^{\infty}}\, ||u^k||_{C^{\infty}} \leq \frac{C}{k^{1/b_2}}\, ||u^k||_{C^{\infty}},\end{eqnarray}

\noindent having used (\ref{eqn:dbar-k-local}), (\ref{eqn:dbar-A-local}) and (\ref{eqn:A-Ak-types-close}). We have thus obtained

\begin{equation}\label{eqn:dbar-k-uk-almost-Cinf}||\bar\partial_k(u^k)||_{C^{\infty}} \leq \frac{C}{k^{1/b_2}}\, ||u^k||_{C^{\infty}}.\end{equation}

\noindent It is clear that the same estimate also holds with $C^0$-norms and $L^2$-norms in place of $C^{\infty}$-norms. So we also have

\begin{equation}\label{eqn:dbar-k-uk-almost-C0-L2}||\bar\partial_k(u^k)||_{C^0} \leq \frac{C}{k^{1/b_2}}\, ||u^k||_{C^0} \hspace{2ex} \mbox{and} \hspace{2ex} ||\bar\partial_k(u^k)||_{L^2(U)} \leq \frac{C}{k^{1/b_2}}\, ||u^k||_{L^2(U)}.\end{equation}

 Now $u^k$ can be regarded as a $C^{\infty}$ section of $L_k$ over $U$. If $e^{(k)}$ denotes the $C^{\infty}$ local frame of $L_k$ corresponding to the trivialisation $\theta_k$ over $U$, we have

\begin{equation}\label{eqn:u^k-section}u^k(z)=e^{-\frac{k}{2}\,\varphi(z)}\stackrel{\theta_k}{\simeq}f_k(z)\otimes e^{(k)}(z), \hspace{3ex} z\in U,\end{equation}

\noindent where $f_k:=\theta_k(u^k)$ is the $C^{\infty}$ function on $U$ representing the section $u^k$ of $L_k$ in the trivialisation $\theta_k$. With respect to the fibre metric $h_k$ of $L_k$ we have\!\!:

\begin{equation}\label{eqn:av-u-hk-rep}|u^k(z)| = e^{-\frac{k}{2}\,\varphi(z)} = |f_k(z)\otimes e^{(k)}(z)|_{h_k},\end{equation}

\noindent where $|u^k|$ is the modulus of the function $u^k$, while $|f_k\otimes e^{(k)}|_{h_k}$ is the pointwise $h_k$-norm of the corresponding local section of $L_k$.

 The crucial estimate (\ref{eqn:dbar-k-uk-almost-Cinf}) shows that the Gaussian function $u^k$, viewed as a local $C^{\infty}$ section of $L_k$, is an {\it approximately holomorphic} section of $L_k$ over $U$ in the strong sense of the $C^{\infty}$-norm (cf. the much weaker $L^2$-norm inequality (\ref{eqn:L2-ah})). In the special case where $\lambda_1(x)>0$, the local section $u^k$ {\it peaks} at $x$ and does increasingly so as $k\rightarrow +\infty$.

 We can go further and define {\it jets} of {\it approximately holomorphic} sections of $L_k$ at $x$. Notice that for every $m_1, \dots , m_n\in\N$ we have on $U$\!\!:

$$\bar\partial_A(z_1^{m_1}\dots z_n^{m_n}\,e^{-\frac{1}{2}\,\varphi})=0 \hspace{2ex} \mbox{and} \hspace{2ex} \bar\partial_{kA}(z_1^{m_1}\dots z_n^{m_n}\,e^{-\frac{k}{2}\,\varphi})=0,$$

\noindent so, in particular, each $z_1^{m_1}\dots z_n^{m_n}\,e^{-\frac{k}{2}\,\varphi} = z_1^{m_1}\dots z_n^{m_n}\, u^k$ satisfies the same estimates (\ref{eqn:dbar-k-uk-almost-Cinf}) and (\ref{eqn:dbar-k-uk-almost-C0-L2}) as $u^k$ does. This motivates the following

\begin{Def}\label{Def:ah-jets} For all $m, k\in\N$ and every $x\in X$, the space of $m$-jets of approximately holomorphic sections of $L_k$ at $x$ is set to be

$$(J^mL_k)_x:=\bigg\{\sum\limits_{m_1+\dots + m_n\leq m}c_{(m_1, \dots , m_n)}\,z_1^{m_1}\dots z_n^{m_n}\,e^{-\frac{k}{2}\,\varphi} \,\,;\,\,c_{(m_1, \dots , m_n)}\in\C\bigg\},$$

\noindent where $z_1, \dots , z_n$ are local holomorphic coordinates of $X$ centred on $x$.

\end{Def}

 As with $u^k$, these jets can be regarded as $C^{\infty}$ sections of $L_k$ over $U$\!\!:\\

$z_1^{m_1}\dots z_n^{m_n}\,u^k(z)=z_1^{m_1}\dots z_n^{m_n}\,e^{-\frac{k}{2}\,\varphi(z)}\stackrel{\theta_k}{\simeq}f_k^{(m_1, \dots ,\, m_n)}(z)\otimes e^{(k)}(z), \hspace{2ex} z\in U,$\\
 
\noindent and the norms are given by

\begin{equation}\label{eqn:av-u-hk-rep_jets}|z_1^{m_1}\dots z_n^{m_n}\,u^k(z)| = z_1^{m_1}\dots z_n^{m_n}\,e^{-\frac{k}{2}\,\varphi(z)} = |f_k^{(m_1, \dots ,\, m_n)}(z)\otimes e^{(k)}(z)|_{h_k}.\end{equation}

 In particular, $(J^0L_k)_x$ consists of multiples of $u^k$. These jets will be used later on to show that a Kodaira-type map defined by {\it approximately holomorphic} sections of $L_k$ is an embedding when $\alpha>0$ and $k\gg 1$. (The reader may wish to compare the discussion of the {\it integrable} case treated in [Tia90] where (jets of) holomorphic peak sections (in the $L^2$-sense) are constructed in high tensor powers of an ample holomorphic line bundle -- a strategy that inspired in part our present treatment of the {\it non-integrable} case.)  

 The next step is to construct a global $C^{\infty}$ section of $L_k$ that belongs to ${\cal H}_k$ (so is {\it approximately holomorphic} in the $L^2$-norm sense) starting from the locally defined {\it peak section} $u^k$. We can define a global section $s$ by multiplying by a cut-off function $\theta$ with support in a neighbourhood of $x$. However, there is no reason that $s$ constructed in this fashion should belong to ${\cal H}_k$, so we need to correct it in the most economical way possible to bring it into ${\cal H}_k$. This will be done in the next subsection for an arbitrary $s\in C^{\infty}(X, \, L_k)$.

\subsection{Approximately holomorphic corrections of global $C^{\infty}$ sections}\label{subsection:corrections}

 For every $k\in\N^{\star}$, the non-negative, formally self-adjoint Laplace-Beltrami operator $\Delta_k'':C^{\infty}(X, \, L_k)\rightarrow C^{\infty}(X, \, L_k)$ is elliptic. So, since $X$ is compact, there is an orthonormal basis $(e_j)_{j\in\N}$ of $C^{\infty}(X, \, L_k)$ consisting of eigenvectors of $\Delta_k''$, while the spectrum of $\Delta_k''$ is discrete with $+\infty$ its only accumulation point. Fix any constant $\delta>0$ independent of $k$ and let

\begin{equation}\label{eqn:delta-eigenvalues}0\leq\mu_0\leq\ldots\leq\mu_{N_k} < \delta \leq \mu_{N_k+1}\leq\ldots\end{equation} 

\noindent denote the eigenvalues (ordered non-decreasingly) of $\Delta_k''$ acting in bi-degree $(0, \, 0)$. Thus $\Delta_k''e_j=\mu_je_j$ for every $j$. (Actually $e_j=e_{k, \, j}$ and $\mu_j=\mu_{k, \, j}$ depend on $k$ but we drop the $k$-index to lighten the notation.) The corresponding eigenspaces are denoted $E_{\Delta_k''}^{0, \, 0}(\mu_j)\subset C^{\infty}(X, \, L_k)$ (and similarly $E_{\Delta_k''}^{p, \, q}(\mu)\subset C^{\infty}_{p, \, q}(X, \, L_k)$ for $\Delta_k'':C^{\infty}_{p, \, q}(X, \, L_k)\rightarrow C^{\infty}_{p, \, q}(X, \, L_k)$). Set

\begin{equation}\label{eqn:delta-def-Hk}\widetilde{{\cal H}_k}:=\bigoplus\limits_{\mu<\delta}E_{\Delta_k''}^{0, \, 0}(\mu) \hspace{3ex} \mbox{and} \hspace{3ex} {\cal N}_k:= \bigoplus\limits_{\mu\geq\delta}E_{\Delta_k''}^{0, \, 0}(\mu)\end{equation}

\noindent (with the understanding that $E_{\Delta_k''}^{0, \, 0}(\mu)=\{0\}$ if $\mu$ is not an eigenvalue of $\Delta_k''$). We can regard $\widetilde{{\cal H}_k}$ as the space of {\it approximately holomorphic} sections of $L_k$ (in a sense less restrictive than for ${\cal H}_k$ of Definition \ref{Def:space-ah-sections}), while ${\cal N}_k$ is its orthogonal complement in $C^{\infty}(X, \, L_k)$.

\begin{Rem}\label{Rem:two-Hks} In the special case when $\alpha>0$, Corollary \ref{Cor:spectral-gap} shows that

\begin{equation}\label{eqn:precise-delta-eigenvalues}0\leq\mu_0\leq\ldots\leq\mu_{N_k} \leq \frac{C}{k^{1+\varepsilon}} < \delta < (\delta_0-\varepsilon_0)k \leq \mu_{N_k+1}\leq\ldots \hspace{2ex} \mbox{if}\,\, k\gg 1.\end{equation} 

 In particular, Laeng's space ${\cal H}_k$ of {\it approximately holomorphic} sections of $L_k$ introduced in Definition \ref{Def:space-ah-sections} coincides with the space $\widetilde{{\cal H}_k}$ defined in (\ref{eqn:delta-def-Hk}) if $k$ is large enough.

\end{Rem}

 In all cases we have an orthogonal splitting

\begin{equation}\label{eqn:orth-splitting}C^{\infty}(X, \, L_k) = \widetilde{{\cal H}_k} \oplus {\cal N}_k.\end{equation}

For every $j=0, \ldots, N_k$, let 

$$P_{k, \, j}:C^{\infty}(X, \, L_k)\longrightarrow \C e_j$$

\noindent be the orthogonal projection onto the $\C$-vector line of $C^{\infty}(X, \, L_k)$ generated by $e_j$. We introduce the following operator.

\begin{Def}\label{Def:Pk-def} For every $k\in\N^{\star}$, let

\begin{equation}\label{eqn:Pk-def}P_k:=\Delta_k'' - \sum\limits_{j=0}^{N_k}\mu_jP_{k, \, j}:C^{\infty}(X, \, L_k) \longrightarrow C^{\infty}(X, \, L_k).\end{equation}

\end{Def}

 It is clear that $\ker P_k=\widetilde{{\cal H}_k}$. Meanwhile, if $s\in C^{\infty}(X, \, L_k)$ is an arbitrary section, there is an orthogonal splitting of $s$ induced by (\ref{eqn:orth-splitting})\!\!:

\begin{equation}\label{eqn:orth-splitting-s}s = s_h + s_{nh}, \hspace{3ex} \mbox{with}\,\,\, s_h\in\widetilde{{\cal H}_k}, s_{nh}\in{\cal N}_k.\end{equation}

\noindent Moreover, if $s=\sum\limits_{j=0}^{+\infty}c_je_j$ (with $c_j\in\C$) is the decomposition of $s$ w.r.t. the orthonormal basis $(e_j)_{j\in\N}$, we have $s_h = \sum\limits_{j=0}^{N_k}c_je_j$, $s_{nh} = \sum\limits_{j\geq N_k+1}c_je_j$ and

\begin{equation}\label{eqn:Pk-s-formula}P_ks = \sum\limits_{j\geq N_k+1}^{+\infty}\mu_jc_je_j = \Delta_k''\bigg(\sum\limits_{j\geq N_k+1}c_je_j\bigg) = \Delta_k''s_{nh}\in{\cal N}_k.\end{equation}

\noindent In particular, $P_{k|{\cal N}_k} = \Delta''_{k|{\cal N}_k}$. If we denote by $P_k^{-1}:{\cal N}_k\rightarrow{\cal N}_k$ the Green operator of $P_k$ (i.e. the inverse of the restriction $P_{k|{\cal N}_k}:{\cal N}_k\rightarrow{\cal N}_k$), we have

\begin{equation}\label{eqn:s_nh-formula}P_k^{-1}P_ks=\sum\limits_{j\geq N_k+1}^{+\infty}c_je_j = s_{nh}\in{\cal N}_k.\end{equation}

 Our goal in this subsection is to estimate the $L^2$-norm of $s_{nh}$ in terms of the $L^2$-norm of $\bar\partial_k s$ for any section $s\in C^{\infty}(X, \, L_k)$. Being the orthogonal projection of $s$ onto ${\cal N}_k$, $s_{nh}$ has minimal $L^2$-norm among all sections $\xi\in C^{\infty}(X, \, L_k)$ for which $s-\xi\in\widetilde{{\cal H}_k}$. In other words, $s_{nh}$ is the minimal correction of an arbitrary $s\in C^{\infty}(X, \, L_k)$ to an {\it approximately holomorphic} section $s_h$.  

 Estimating the $L^2$-norm of $s_{nh}$ can be seen as a non-integrable analogue in this particular situation of H\"ormander's familiar $L^2$-estimates of the integrable case. Indeed, recall that the standard method of correcting an arbitrary global $C^{\infty}$ section $s$ of a {\it positive holomorphic} line bundle $L_k$ to a global {\it holomorphic} section $s_h$ of $L_k$ is to solve the $\bar\partial$-equation

$$\bar\partial_k\xi = \bar\partial_k s  \hspace{3ex} \mbox{on}\,\,\, X$$

\noindent by selecting the solution $\xi\in C^{\infty}(X, \, L_k)$ of minimal $L^2$-norm which is given explicitly by the familiar formula 

\begin{equation}\label{eqn:hol-xi-formula}\xi=G_k\bar\partial_k^{\star}(\bar\partial_k s)\end{equation} 

\noindent (where $G_k$ stands for the Green operator of $\Delta_k''$) and to set $s_h=s-\xi$. It is clear that $s_h$ and $\xi$ are nothing but the orthogonal projections of $s$ onto the subspace of global holomorphic sections and respectively its orthogonal complement in $C^{\infty}(X, \, L_k)$. In our non-integrable case, the roles of these two subspaces are played by $\widetilde{{\cal H}_k}$ and respectively ${\cal N}_k$, while formula (\ref{eqn:s_nh-formula}) is the analogue of (\ref{eqn:hol-xi-formula}) for $s_{nh}=\xi$.

\vspace{2ex}

 We shall now obtain the desired estimate in the particular case when the initial $(1, \, 1)$-form $\alpha$ is supposed to be {\it strictly positive} on $X$. (Recall that in the general case $\alpha$ is only to be assumed to satisfy Demailly's tremendously weaker hypothesis $\int_{X(\alpha, \, \leq 1)}\alpha^n>0$.) This assumption, which parallels the strict positivity curvature assumption in H\"ormander's $L^2$-estimates of the integrable case, will be relaxed in future work.

\begin{Prop}\label{Prop:s_nh-L2-est} Suppose $\alpha$ is a $C^{\infty}$ {\bf positive-definite} $d$-closed $(1, \, 1)$-form of {\bf possibly non-rational} De Rham cohomology class on a compact complex Hermitian manifold $(X, \, \omega)$. Fix an arbitrary constant $\delta>0$.

 Then, with the notation of subsections \ref{subsection:setting}, \ref{subsection:Lresults} and \ref{subsection:corrections} in place, the following property holds. For every $s\in C^{\infty}(X, \, L_k)$, the non-approximately-holomorphic component $s_{nh}$ of $s$ (cf. (\ref{eqn:orth-splitting-s}) and (\ref{eqn:delta-def-Hk})) satisfies the estimate\!\!:

\begin{equation}\label{eqn:s_nh-L2-est}||s_{nh}||^2\leq\frac{4}{\delta_0k}\bigg(1+ \frac{C}{\delta_0\delta^2}\,\frac{1}{k^{1+\frac{2}{b_2}}}\bigg)||\bar\partial_ks||^2, \hspace{3ex} k\geq k_{\delta},\end{equation}

\noindent where $||\,\,\,\,||$ stands for $L^2$-norm, $\delta_0>0$ is any constant for which $\alpha \geq 2\delta_0\omega$, while $C>0$ is a constant depending only on $(X, \, \omega)$ and $k_{\delta}\in\N^{\star}$ depends only on $\delta>0$ and $\alpha$.

\end{Prop}

\noindent {\it Proof.} By Remark \ref{Rem:two-Hks}, $\widetilde{{\cal H}_k}={\cal H}_k$ under the present assumptions. Note that by (\ref{eqn:Pk-s-formula}) and (\ref{eqn:delta-eigenvalues}) we have 

$$\langle\langle P_ks, \, s\rangle\rangle = \langle\langle\Delta_k''s_{nh}, \, s_{nh}\rangle\rangle\geq\mu_{N_k+1}||s_{nh}|^2\geq\delta ||s_{nh}|^2,  \hspace{3ex}  s\in C^{\infty}(X, \, L_k).$$

\noindent If we extend the Green operator $P_k^{-1}:{\cal N}_k\rightarrow {\cal N}_k$ to $P_k^{-1}:C^{\infty}(X, \, L_k)\rightarrow {\cal N}_k$ by letting $(P_k^{-1})_{|{\cal H}_k}=0$, we infer

\begin{equation}\label{eqn:P-inverse-00-ubound}\langle\langle P_k^{-1}s, \, s\rangle\rangle = \langle\langle P_k^{-1}s_{nh}, \, s_{nh}\rangle\rangle\leq\frac{1}{\mu_{N_k+1}}||s_{nh}||^2 \leq \frac{1}{\delta} ||s_{nh}||^2 \leq \frac{1}{\delta} ||s||^2\end{equation}

\noindent for all $s\in C^{\infty}(X, \, L_k)$, where the last inequality follows from $s_h$ and $s_{nh}$ being orthogonal (hence $||s||^2=||s_h||^2 + ||s_{nh}||^2\geq ||s_{nh}||^2$).

 On the other hand, definition (\ref{eqn:Pk-def}) of $P_k$ makes sense in any bi-degree $(p, \, q)$ and gives an operator $P_k:C^{\infty}_{p, \, q}(X, \, L_k)\rightarrow C^{\infty}_{p, \, q}(X, \, L_k)$ defined by the same formula (\ref{eqn:Pk-def}) as in bi-degree $(0, \, 0)$ if we make the convention that $P_{k, \, j}$ and $\mu_jP_{k, \, j}$ are the zero operator when $\mu_j$ is not an eigenvalue of $\Delta_k'': C^{\infty}_{p, \, q}(X, \, L_k) \rightarrow C^{\infty}_{p, \, q}(X, \, L_k)$. Indeed, due to the non-commutation of $\Delta_k''$ with $\bar\partial_k$ (because $\bar\partial_k^2\neq 0$), the eigenvalues $\mu_j$ of $\Delta_k''$ in bi-degree $(0, \, 0)$ need not be eigenvalues of $\Delta_k''$ in bi-degree $(p, \, q)\neq (0, \, 0)$. Furthermore, Corollary \ref{Cor:W-pos-cons} of the Weitzenb\"ock formula shows that in bi-degree $(0, \, 1)$ we have

$$P_k = \Delta_k''\geq \delta_0k>0  \hspace{3ex} \mbox{on} \hspace{1ex} C^{\infty}_{0, \, 1}(X, \, L_k) \hspace{3ex} \mbox{if} \hspace{1ex} k\gg 1 \,\,\,(\mbox{i.e. if}\,\, k>\frac{\delta}{\delta_0}).$$

\noindent (Implicitly $P_{k, \, j}=0$ in bi-degree $(0, \, 1)$ for all $j=0, \ldots , N_k$ and all large $k$.) So $P_k:C^{\infty}_{0, \, 1}(X, \, L_k)\rightarrow C^{\infty}_{0, \, 1}(X, \, L_k)$ is invertible for $k$ large enough and its inverse satisfies the estimate\!\!:

\begin{equation}\label{eqn:P-01-inverse-ubound}P_k^{-1} \leq \frac{1}{\delta_0k}  \hspace{3ex} \mbox{on} \hspace{1ex} C^{\infty}_{0, \, 1}(X, \, L_k) \hspace{3ex} \mbox{if} \hspace{1ex} k\gg 1.\end{equation}

 Let us introduce the operator

$$Q_k:=P_k^{-1}\bar\partial_k^{\star} - \bar\partial_k^{\star}P_k^{-1}: C^{\infty}_{0,\, 1}(X, \, L_k)\rightarrow C^{\infty}(X, \, L_k)$$

\noindent measuring the {\it commutation defect} of $P_k^{-1}$ with $\bar\partial_k^{\star}$. Similarly we set

$$S_k:=\bar\partial_k^{\star} P_k - P_k\bar\partial_k^{\star}: C^{\infty}_{0,\, 1}(X, \, L_k)\rightarrow C^{\infty}(X, \, L_k)$$

\noindent which measures the {\it commutation defect} of $P_k$ with $\bar\partial_k^{\star}$. We clearly have

\begin{equation}\label{eqn:Q_k-S_k-rel}Q_k=P_k^{-1}S_kP_k^{-1}.\end{equation}

\noindent In all these expressions, $P_k^{-1}$ and $P_k$ act on $C^{\infty}(X, \, L_k)$ or $C^{\infty}_{0, \, 1}(X, \, L_k)$ according to case.

\vspace{2ex}

 Now fix an arbitrary $s\in C^{\infty}(X, \, L_k)$. Using (\ref{eqn:s_nh-formula}) and (\ref{eqn:Pk-s-formula}), we get

$$s_{nh}=P_k^{-1}P_ks=P_k^{-1}\Delta_k''s_{nh}  = P_k^{-1}\bar\partial_k^{\star}\bar\partial_ks_{nh},$$

\noindent which, after writing $P_k^{-1}\bar\partial_k^{\star} = \bar\partial_k^{\star}P_k^{-1} + Q_k$, transforms to

\begin{equation}\label{eqn:s-nh-Qk}s_{nh}=\bar\partial_k^{\star}P_k^{-1}\bar\partial_ks_{nh} + Q_k\bar\partial_ks_{nh}.\end{equation}

\noindent We shall estimate separately the $L^2$-norms of $\bar\partial_k^{\star}P_k^{-1}\bar\partial_ks_{nh}$ and $Q_k\bar\partial_ks_{nh}$.

\vspace{2ex}

 In the case of $\bar\partial_k^{\star}P_k^{-1}\bar\partial_ks_{nh}$, we have\!\!:

\begin{eqnarray}\label{eqn:L2est-part1}\nonumber||\bar\partial_k^{\star}P_k^{-1}\bar\partial_ks_{nh}||^2 & = & \langle\langle \bar\partial_k^{\star}P_k^{-1}\bar\partial_ks_{nh}, \, \bar\partial_k^{\star}P_k^{-1}\bar\partial_ks_{nh} \rangle\rangle \\
\nonumber & = & \langle\langle \bar\partial_k\bar\partial_k^{\star}P_k^{-1}\bar\partial_ks_{nh}, \, P_k^{-1}\bar\partial_ks_{nh} \rangle\rangle \\
\nonumber & = & \langle\langle (\Delta_k'' - \bar\partial_k^{\star}\bar\partial_k)P_k^{-1}\bar\partial_ks_{nh}, \, P_k^{-1}\bar\partial_ks_{nh} \rangle\rangle \\
\nonumber & = & \langle\langle (P_k  - \bar\partial_k^{\star}\bar\partial_k)P_k^{-1}\bar\partial_ks_{nh}, \, P_k^{-1}\bar\partial_ks_{nh} \rangle\rangle \\
\nonumber & = & \langle\langle \bar\partial_ks_{nh}, \, P_k^{-1}\bar\partial_ks_{nh} \rangle\rangle - \langle\langle \bar\partial_k P_k^{-1}\bar\partial_ks_{nh}, \, \bar\partial_k P_k^{-1}\bar\partial_ks_{nh}\rangle\rangle \\
\nonumber & = & \langle\langle P_k^{-1}\bar\partial_ks_{nh}, \, \bar\partial_ks_{nh}\rangle\rangle - ||\bar\partial_k P_k^{-1}\bar\partial_ks_{nh}||^2 \\
          & \leq & \langle\langle P_k^{-1}\bar\partial_ks_{nh}, \, \bar\partial_ks_{nh}\rangle\rangle \leq \frac{1}{\delta_0k}\,||\bar\partial_ks_{nh}||^2  \hspace{2ex} \mbox{if} \,\, k\gg 1.\end{eqnarray}

\noindent In going from the third to the fourth line above, we have used the identity $\Delta_k''=P_k$ on $C^{\infty}_{0, \, 1}(X, \, L_k)$ (see the discussion preceding (\ref{eqn:P-01-inverse-ubound}) -- a consequence of the Weitzenb\"ock formula), while the passage from the fifth to the sixth line has used the fact that $\langle\langle \bar\partial_ks_{nh}, \, P_k^{-1}\bar\partial_ks_{nh} \rangle\rangle = \overline{\langle\langle P_k^{-1}\bar\partial_ks_{nh}, \, \bar\partial_ks_{nh}\rangle\rangle}$ is real. The last inequality follows from estimate (\ref{eqn:P-01-inverse-ubound}) 

\vspace{2ex}

 We shall now estimate the second term $Q_k\bar\partial_ks_{nh}$ in the expression (\ref{eqn:s-nh-Qk}) of $s_{nh}$. We shall actually estimate the $L^2$-norm of $S_k\bar\partial_ks_{nh}$ and then use (\ref{eqn:Q_k-S_k-rel}) to go from $S_k$ to $Q_k$. 

 Since $P_k = \Delta_k''$ on $C^{\infty}_{0, \, 1}(X, \, L_k)$ (hence also on $\bar\partial_ks_{nh}$) for $k\gg 1$, we get

\begin{eqnarray}\nonumber S_k\bar\partial_ks_{nh} & = & (\bar\partial_k^{\star}\Delta_k'' - \Delta_k''\bar\partial_k^{\star})(\bar\partial_ks_{nh})  + \sum\limits_{j=1}^{N_k}\mu_jP_{k, \, j}\bar\partial_k^{\star}(\bar\partial_ks_{nh}) \\
\nonumber  & = & \bar\partial_k^{\star 2}\bar\partial_k(\bar\partial_ks_{nh}) + \sum\limits_{j=1}^{N_k}\mu_jP_{k, \, j}\bar\partial_k^{\star}(\bar\partial_ks_{nh})  \hspace{15ex} \mbox{if} \hspace{1ex} k\gg 1.\end{eqnarray}

\noindent (Indeed, $\bar\partial_k^{\star}\Delta_k'' - \Delta_k''\bar\partial_k^{\star} = \bar\partial_k^{\star 2}\bar\partial_k - \bar\partial_k\bar\partial_k^{\star 2}$ but $\bar\partial_k^{\star 2}(\bar\partial_ks_{nh})=0$ for bi-degree reasons.) Since $\Delta_k'' = \bar\partial_k^{\star}\bar\partial_k$ in bi-degree $(0, \, 0)$, we see that $P_{k, \, j}\bar\partial_k^{\star}(\bar\partial_ks_{nh}) = P_{k, \, j}(\Delta_k''s_{nh}) =0$ for every $j\in\{1, \ldots , N_k\}$. Indeed, $s_{nh}\in{\cal N}_k$ by definition, so $\Delta_k''s_{nh}\in{\cal N}_k$, while $P_{k, \, j}$ is the orthogonal projection onto a subspace of ${\cal H}_k = {\cal N}_k^{\perp}.$ Thus

\begin{equation}\label{eqn:Sk-s_nh-expression}S_k\bar\partial_ks_{nh} = \bar\partial_k^{\star 2}\bar\partial_k(\bar\partial_ks_{nh})  \hspace{3ex} \mbox{if} \hspace{1ex} k\gg 1.\end{equation}

 We pause briefly to prove in full generality (i.e. without using the positivity assumption made on $\alpha$ in Proposition \ref{Prop:s_nh-L2-est}) the following estimate reminiscent of Laeng's estimate (\ref{eqn:comm-defect-Laurent}).

\begin{Lem}\label{Lem:comm-defect-est} For every section $s\in C^{\infty}(X, \, L_k)$ we have

$$||\bar\partial_k^{\star 2}\bar\partial_k(\bar\partial_ks)||^2 \leq \frac{C}{k^{2/b_2}}\,\bigg(k||s||^2 + ||\bar\partial_ks||^2\bigg), \hspace{3ex} k\in\N^{\star},$$

\noindent where $C>0$ is a constant independent of $k$.

\end{Lem}

\noindent {\it Proof of Lemma \ref{Lem:comm-defect-est}.}  Recall that the fundamental commutation relations for non-K\"ahler metrics (that are common to the integrable and non-integrable cases -- see e.g. [Dem85b] or [Lae02] or [Don96]) give (cf. notation in subsection \ref{subsection:Weitzenbock})\!\!:

\begin{equation}\label{eqn:comm-rel}i(\bar\partial_k^{\star} + \bar\tau^{\star})=[\Lambda, \, \partial_k]  \hspace{3ex}  \mbox{or equivalently} \hspace{3ex} \bar\partial_k^{\star} = -i[\Lambda, \, \partial_k] - \bar\tau^{\star}.\end{equation}

 Thus for every $\sigma\in C^{\infty}_{0, \, 1}(X, \, L_k)$ we get from (\ref{eqn:comm-rel})\!\!:

\begin{eqnarray}\nonumber \bar\partial_k^{\star 2}\bar\partial_k\sigma & = & (i[\Lambda, \, \partial_k] + \bar\tau^{\star})(i[\Lambda, \, \partial_k] + \bar\tau^{\star})\bar\partial_k\sigma\\
\nonumber  & = & (i[\Lambda, \, \partial_k] + \bar\tau^{\star})(i\Lambda\partial_k\bar\partial_k\sigma + \bar\tau^{\star}(\bar\partial_k\sigma))\\
\nonumber  & = & -\Lambda\partial_k\Lambda\partial_k\bar\partial_k\sigma + i\Lambda\partial_k\bar\tau^{\star}(\bar\partial_k\sigma) + i\bar\tau^{\star}\Lambda\partial_k\bar\partial_k\sigma + \bar\tau^{\star 2}(\bar\partial_k\sigma).\end{eqnarray}

\noindent Since the latter half of (\ref{eqn:del-squared}) amounts to $\partial_k\bar\partial_k = -2\pi i\,\alpha_k^{1, \, 1}- \bar\partial_k\partial_k$, we get\!\!:

\begin{eqnarray}\label{eqn:dbar-square-dbar-first}\nonumber \bar\partial_k^{\star 2}\bar\partial_k\sigma = \Lambda\partial_k\Lambda\,(2\pi i\,\alpha_k^{1, \, 1}\wedge\sigma & + & \bar\partial_k\partial_k\sigma) + i\bar\tau^{\star}\Lambda\,(-2\pi i\alpha_k^{1, \, 1}\wedge\sigma- \bar\partial_k\partial_k\sigma)\\
    & + & i\Lambda\partial_k\bar\tau^{\star}(\bar\partial_k\sigma) + \bar\tau^{\star 2}(\bar\partial_k\sigma).\end{eqnarray}

  Now suppose that $\sigma=\bar\partial_ks\in C^{\infty}_{0, \, 1}(X, \, L_k)$ for some $s\in C^{\infty}(X, \, L_k)$. Then\!\!:

\begin{eqnarray}\nonumber\bar\partial_k\partial_k\sigma & = & \bar\partial_k(\partial_k\bar\partial_ks) = -2\pi i\,\bar\partial_k(\alpha_k^{1, \, 1}\wedge s) - \bar\partial_k^2\partial_ks \\
\nonumber  & = & -2\pi i\,\bar\partial\alpha_k^{1, \, 1}\wedge s - 2\pi i \,\alpha_k^{1, \, 1}\wedge\bar\partial_k s + 2\pi i \,\alpha_k^{0, \, 2}\wedge\partial_ks.\end{eqnarray}

\noindent having used the former half of (\ref{eqn:del-squared}) to obtain the last term. Now $d\alpha_k=0$, so passing to the $(1, \, 2)$-component we see that $\bar\partial\alpha_k^{1, \, 1} = -\partial\alpha_k^{0, \, 2}$. Hence we get\!\!:

$$\bar\partial_k\partial_k\sigma = 2\pi i\,(\partial\alpha_k^{0, \, 2}\wedge s - \alpha_k^{1, \, 1}\wedge\bar\partial_ks + \alpha_k^{0, \, 2}\wedge\partial_ks),$$

\noindent from which, since $\bar\partial_ks=\sigma$, we further get\!\!:

\begin{equation}\label{eqn:dbar-square-dbar-second}2\pi i\,\alpha_k^{1, \, 1}\wedge\sigma + \bar\partial_k\partial_k\sigma = 2\pi i\, (\partial\alpha_k^{0, \, 2}\wedge s + \alpha_k^{0, \, 2}\wedge\partial_ks).\end{equation}

\noindent Combining (\ref{eqn:dbar-square-dbar-first}) and (\ref{eqn:dbar-square-dbar-second}), we get \!\!:

\begin{eqnarray}\label{eqn:dbar-square-dbar-third} \bar\partial_k^{\star 2}\bar\partial_k\sigma & = & 2\pi i\, \Lambda\partial_k\Lambda\,(\partial\alpha_k^{0, \, 2}\wedge s + \alpha_k^{0, \, 2}\wedge\partial_ks) + \\
\nonumber  &  & 2\pi\, \bar\tau^{\star}\Lambda\,(\partial\alpha_k^{0, \, 2}\wedge s + \alpha_k^{0, \, 2}\wedge\partial_ks) + i\Lambda\partial_k\bar\tau^{\star}(\bar\partial_k^2s) + \bar\tau^{\star 2}(\bar\partial_k^2s)\end{eqnarray}

\noindent for all $\sigma = \bar\partial_ks\in C^{\infty}_{0, \, 1}(X,\, L_k)$. Recall that $\bar\partial_k^2s = -2\pi i\,\alpha_k^{0, \, 2}\wedge s$ (cf. (\ref{eqn:del-squared})) and that $||\alpha_k^{0, \, 2}||_{C^{\infty}}\leq\frac{C}{k^{1/b_2}}$ (cf. (\ref{eqn:prop-def-alpha-k-bis})). Consequently, in the above identity (\ref{eqn:dbar-square-dbar-third}), the $L^2$-norms of the expressions $\bar\partial_k^2s$ (a zero-order operator acting on $s$) and $\partial\alpha_k^{0, \, 2}\wedge s + \alpha_k^{0, \, 2}\wedge\partial_ks$ (a first-order operator acting on $s$) can be controlled in terms of the $L^2$-norms of $s$ and $\partial_ks$. Meanwhile, $\Lambda$, $\bar\tau^{\star}$ are zero-order operators independent of $k$, hence bounded independently of $k$. Thus so are $\bar\tau^{\star}\Lambda$ and $\bar\tau^{\star 2}$, too.

 Putting these facts together, we see that the $L^2$-norms of the terms featuring on the right of (\ref{eqn:dbar-square-dbar-third}) are estimated as follows (where the constant $C>0$ is independent of $k$ and is allowed to vary from line to line)\!\!:

\begin{equation}\label{eqn:dbar-square-dbar-fourth}\nonumber ||\partial\alpha_k^{0, \, 2}\wedge s + \alpha_k^{0, \, 2}\wedge\partial_ks|| \leq \frac{C}{k^{1/b_2}}\,(||s|| + ||\partial_ks||),\end{equation}

\noindent and similarly

\begin{equation}\label{eqn:dbar-square-dbar-fifth}\nonumber ||\Lambda\partial_k\Lambda\,(\partial\alpha_k^{0, \, 2}\wedge s + \alpha_k^{0, \, 2}\wedge\partial_ks)|| \leq \frac{C}{k^{1/b_2}}\,(||s|| + ||\partial_ks||),\end{equation}

\begin{equation}\label{eqn:dbar-square-dbar-sixth}\nonumber ||\bar\tau^{\star}\Lambda\,(\partial\alpha_k^{0, \, 2}\wedge s + \alpha_k^{0, \, 2}\wedge\partial_ks)|| \leq \frac{C}{k^{1/b_2}}\,(||s|| + ||\partial_ks||),\end{equation}

\begin{equation}\label{eqn:dbar-square-dbar-seventh}\nonumber ||\Lambda\partial_k\bar\tau^{\star}(\alpha_k^{0, \, 2}\wedge s)|| \leq \frac{C}{k^{1/b_2}}\,(||s|| + ||\partial_ks||),\end{equation}

\begin{equation}\label{eqn:dbar-square-dbar-eightth}\nonumber ||\bar\tau^{\star 2}\,(\alpha_k^{0, \, 2}\wedge s)|| \leq \frac{C}{k^{1/b_2}}\,||s||.\end{equation}

\noindent If we now take the squared $L^2$-norm on either side of (\ref{eqn:dbar-square-dbar-third}) (with $\sigma=\bar\partial_ks\in C^{\infty}_{0, \, 1}(X, \, L_k)$), the above estimates add up to

\begin{equation}\label{eqn:dbar-square-dbar-seventhth}||\bar\partial_k^{\star 2}\bar\partial_k(\bar\partial_ks)||^2 \leq \frac{C}{k^{2/b_2}}\,(||s||^2 + ||\partial_ks||^2)  \hspace{3ex} \mbox{for all} \hspace{1ex} s\in C^{\infty}(X, \, L_k).\end{equation}

 Since $\bar\partial_ks$ (measuring how far short a section $s$ of $L_k$ falls from being holomorphic) is better adapted to our purposes than $\partial_ks$, we wish to replace $\partial_ks$ by $\bar\partial_ks$ in the right-hand side of the above estimate (\ref{eqn:dbar-square-dbar-seventhth}). The transition from $\partial_ks$ to $\bar\partial_ks$ was done in a natural way by Laeng in [Lae02, p. 89] using the Bochner-Kodaira-Nakano identity (\ref{eqn:BKN-nK-nH}) which, when specialised to the case $E=(L_k, \, h_k, \, D_k)\rightarrow (X, \, \omega)$, reads

$$\Delta_k'' = \Delta_{k, \, \tau}' + 2\pi\,[\alpha_k^{1, \, 1}, \, \Lambda] + T_{\omega}.$$

\noindent This allows one to express $||\bar\partial_ks||^2 = \langle\langle\Delta_k''s, \, s\rangle\rangle$ in terms of $||\partial_ks + \tau s||^2 = \langle\langle\Delta_{k, \, \tau}'s, \, s\rangle\rangle$. The straightforward calculation performed in [Lae02, p. 89] gave the estimate\!\!:

\begin{equation}\label{eqn:dbar_k-del_k-passage}||\partial_ks||^2 \leq C(k||s||^2 + ||\bar\partial_ks||^2),  \hspace{3ex}  s\in C^{\infty}(X, \, L_k).\end{equation}

\noindent Note that the factor $k$ of $||s||^2$ comes from the curvature term $\alpha_k^{1, \, 1}$ which is close (in $C^{\infty}$-norm) to $k\alpha$.

 Using (\ref{eqn:dbar_k-del_k-passage}), (\ref{eqn:dbar-square-dbar-seventhth}) transforms to

$$||\bar\partial_k^{\star 2}\bar\partial_k(\bar\partial_ks)||^2 \leq \frac{C}{k^{2/b_2}}\,(k||s||^2 + ||\bar\partial_ks||^2)  \hspace{3ex} \mbox{for all} \hspace{1ex} s\in C^{\infty}(X, \, L_k)$$

\noindent which is precisely the estimate claimed in the statement. The proof of Lemma \ref{Lem:comm-defect-est} is complete.  \hfill $\Box$

\vspace{2ex}

 Thanks to (\ref{eqn:Sk-s_nh-expression}), Lemma \ref{Lem:comm-defect-est} immediately implies the following

\begin{Cor}\label{Cor:comm-defect-s_nh} Under the hypotheses of Proposition \ref{Prop:s_nh-L2-est}, every $s\in C^{\infty}(X, \, L_k)$ satisfies the following estimate\!\!:

\begin{equation}\label{eqn:cor-Sk-est}||S_k\bar\partial_ks_{nh}||^2 \leq \frac{C}{k^{2/b_2}}\,\bigg(k||s||^2 + ||\bar\partial_ks||^2\bigg), \hspace{3ex} k\gg 1,\end{equation}

\noindent where $C>0$ is a constant independent of $k$.

\end{Cor}

\noindent {\it Proof.} Applying Lemma \ref{Lem:comm-defect-est} to $s_{nh}$ and using (\ref{eqn:Sk-s_nh-expression}), we get

\begin{equation}\label{eqn:cor-Sk-est-aux}||S_k\bar\partial_ks_{nh}||^2 \leq \frac{C}{k^{2/b_2}}\,\bigg(k||s_{nh}||^2 + ||\bar\partial_ks_{nh}||^2\bigg), \hspace{3ex} k\gg 1.\end{equation}

 As already noticed, $||s_{nh}||^2\leq ||s||^2$ by virtue of $s_h$ and $s_{nh}$ being orthogonal in the splitting $s=s_h + s_{nh}$.

 On the other hand, taking $\bar\partial_k$ in this splitting, we get $\bar\partial_k s = \bar\partial_k s_h + \bar\partial_ks_{nh}.$ We claim that the $L_k$-valued $(0, \, 1)$-forms $\bar\partial_k s_h$ and $\bar\partial_ks_{nh}$ are orthogonal. Indeed, we see that

$$\langle\langle \bar\partial_ks_h, \, \bar\partial_ks_{nh}\rangle\rangle = \langle\langle \bar\partial_k^{\star}\bar\partial_ks_h, \, s_{nh}\rangle\rangle = \langle\langle \Delta_k''s_h, \, s_{nh}\rangle\rangle = 0.$$

\noindent The reason for the last equality is that $s_h\in{\cal H}_k$ (by construction), hence $\Delta_k''s_h\in{\cal H}_k$, while $s_{nh}\in{\cal N}_k$ (again by construction) and ${\cal H}_k\perp{\cal N}_k$. Thus $\bar\partial_ks_h\perp\bar\partial_ks_{nh}$ and it follows that $||\bar\partial_k s||^2 = ||\bar\partial_k s_h||^2 + ||\bar\partial_ks_{nh}||^2\geq ||\bar\partial_ks_{nh}||^2$.

 It is now clear that the right-hand side of (\ref{eqn:cor-Sk-est-aux}) is $\leq$ than the right-hand side of (\ref{eqn:cor-Sk-est}). This completes the proof.  \hfill $\Box$

\vspace{3ex}

\noindent {\it End of proof of Proposition \ref{Prop:s_nh-L2-est}.} By (\ref{eqn:Q_k-S_k-rel}) we have

$$Q_k\bar\partial_ks_{nh} = P_k^{-1}S_kP_k^{-1}(\bar\partial_ks_{nh}).$$

\noindent Using (\ref{eqn:cor-Sk-est-aux}), (\ref{eqn:P-inverse-00-ubound}) and (\ref{eqn:P-01-inverse-ubound}) we get for every $s\in C^{\infty}(X, \, L_k)$ the estimate

\begin{equation}\label{eqn:Qk-est}||Q_k\bar\partial_ks_{nh}||^2\leq\frac{C}{\delta^2(\delta_0k)^2k^{2/b_2}}\,\bigg(k||s_{nh}||^2 + ||\bar\partial_ks_{nh}||^2\bigg), \hspace{3ex} k\gg 1,\end{equation}

\noindent where $C>0$ is a constant independent of $k$. By (\ref{eqn:P-inverse-00-ubound}), the $\delta^2$ of the above denominator can be improved to $\mu_{N_k+1}^2$, hence also to $(\delta_0-\varepsilon_0)^2k^2$ by (\ref{eqn:precise-delta-eigenvalues}), but this will be of no consequence in the sequel.

 Using now the splitting (\ref{eqn:s-nh-Qk}) and the estimates (\ref{eqn:L2est-part1}) and (\ref{eqn:Qk-est}) of its two terms, we get for every $s\in C^{\infty}(X, \, L_k)$ the estimate

\begin{eqnarray}\nonumber ||s_{nh}||^2 & \leq & 2\,(||\bar\partial_k^{\star}P_k^{-1}\bar\partial_ks_{nh}||^2 + ||Q_k\bar\partial_ks_{nh}||^2)\\
 \nonumber & \leq & \frac{2}{\delta_0k}\,\bigg(1+\frac{C}{\delta_0\delta^2}\,\frac{1}{k^{1+\frac{2}{b_2}}}\bigg)\,||\bar\partial_ks_{nh}||^2 + \frac{2C}{(\delta_0\delta)^2}\,\frac{1}{k^{1+\frac{2}{b_2}}}\,||s_{nh}||^2, \hspace{3ex} k\gg 1,\end{eqnarray}

\noindent which is equivalent to

$$\bigg(1 - \frac{2C}{(\delta_0\delta)^2}\,\frac{1}{k^{1+\frac{2}{b_2}}}\bigg)\,||s_{nh}||^2 \leq \frac{2}{\delta_0k}\,\bigg(1+\frac{C}{\delta_0\delta^2}\,\frac{1}{k^{1+\frac{2}{b_2}}}\bigg)\,||\bar\partial_ks_{nh}||^2, \hspace{3ex} k\gg 1.$$

\noindent Now it is clear that the coefficient on the left-hand side above satisfies

$$\frac{1}{2}\leq 1 - \frac{2C}{(\delta_0\delta)^2}\,\frac{1}{k^{1+\frac{2}{b_2}}} <1  \hspace{3ex} \mbox{for} \hspace{1ex} k\gg 1,$$

\noindent so we get

$$||s_{nh}||^2 \leq \frac{4}{\delta_0k}\,\bigg(1+\frac{C}{\delta_0\delta^2}\,\frac{1}{k^{1+\frac{2}{b_2}}}\bigg)\,||\bar\partial_ks_{nh}||^2 \hspace{3ex} \mbox{for} \hspace{1ex} k\gg 1.$$

\noindent Since $||\bar\partial_ks_{nh}|| \leq ||\bar\partial_ks||$ (as explained in the proof of Corollary \ref{Cor:comm-defect-s_nh}), the above estimate implies estimate (\ref{eqn:s_nh-L2-est}). The proof of Proposition \ref{Prop:s_nh-L2-est} is complete.   \hfill  $\Box$

\subsection{Global approximately holomorphic peak sections}\label{subsection:global-peak} We now bring the discussions of subsections \ref{subsection:local-model} and \ref{subsection:corrections} together. We suppose that $\alpha >0$ on $X$ as in Proposition \ref{Prop:s_nh-L2-est}. As in the previous subsections, the symbol $||\,\,\,||$ will stand for the global $L^2$-norm on $X$ when it has no index, while an index will change its meaning to the norm it indicates. 

 Let $x\in X$ be an arbitrary point and let $U\subset X$ be an open neighbourhood of $x$ as in subsection \ref{subsection:local-model} with local holomorphic coordinates as in (\ref{eqn:choice-local-coord}). Consider, for every $k\in\N^{\star}$, the Gaussian section $u^k$ of $L_k$ over $U$ defined in (\ref{eqn:loc-section-k-def}). Choose an open neighbourhood $V$ of $x$ such that $V\Subset U$ and a $C^{\infty}$ cut-off function $\theta:X\rightarrow \R$ such that 

$$\theta\equiv 1 \hspace{1ex} \mbox{on} \hspace{1ex} V \hspace{1ex} \mbox{and} \hspace{1ex} \mbox{Supp}\,\theta\Subset U.$$

 We can apply the results of subsection \ref{subsection:corrections} to the global section 

$$s:=\theta u^k\in C^{\infty}(X, \, L_k)$$

\noindent whose $s_{nh}$ component satisfies thus the $L^2$-estimate (\ref{eqn:s_nh-L2-est}). This estimate can be refined in the special case of $s=\theta u^k$ using the $C^{\infty}$-estimate (\ref{eqn:dbar-k-uk-almost-Cinf}) satisfied by $u^k$. Indeed, applying $\bar\partial_k$ we get\!\!:

$$\bar\partial_ks = \bar\partial_k(\theta u^k) = \theta\,\bar\partial_ku^k + (\bar\partial\theta)\,u^k,$$

\noindent hence $\bar\partial_ks = 0$ on $X\setminus U$ and $\bar\partial_ks = \bar\partial_ku^k$ on $V$. Thus

$$||\bar\partial_ks||_{C^{\infty}}\leq ||\bar\partial_ku^k||_{C^{\infty}} + C\,||u^k||_{C^{\infty}} \leq C\bigg(1 + \frac{1}{k^{1/b_2}}\bigg)\,||u^k||_{C^{\infty}},$$

\noindent having used (\ref{eqn:dbar-k-uk-almost-Cinf}) to get the last estimate. The analogous estimate holds for $C^0$-norms by (\ref{eqn:dbar-k-uk-almost-C0-L2}), so

\begin{equation}\label{eqn:d-bar_k-local-est}||\bar\partial_ks||_{C^{\infty}}\leq C\,||u^k||_{C^{\infty}}  \hspace{3ex} \mbox{and} \hspace{3ex} ||\bar\partial_ks||_{C^0}\leq C\,||u^k||_{C^0} = C \hspace{2ex} \mbox{for} \,\, k\gg 1.\end{equation}

\noindent (The last identity holds whenever $\alpha\geq 0$ since $u^k(0)=1$ and $u^k(z)\leq u^k(0)$ in this case for all $z\in U\setminus\{0\}$ if $U$ is small enough.) Thus estimate (\ref{eqn:s_nh-L2-est}) yields

\begin{equation}\label{eqn:_nh-L2-est-improved}||s_{nh}||^2\leq C(X, \, \omega)\,\delta_k, \hspace{3ex} k\gg 1, \hspace{3ex} \mbox{where}\hspace{1ex}\delta_k:=\frac{4}{\delta_0k}\,\bigg(1+ \frac{C}{\delta_0\delta^2}\, \frac{1}{k^{1+\frac{2}{b_2}}}\bigg),\end{equation}

\noindent since $||\bar\partial_ks||^2\leq ||\bar\partial_ks||^2_{C^0}\,\mbox{Vol}_{\omega}(X)\leq C^2\,\mbox{Vol}_{\omega}(X)$. We have set $C(X, \, \omega): =C\,\mbox{Vol}_{\omega}(X)>0$.

 We can now go from this global $L^2$-estimate to a local $L^{\infty}$-estimate, but we first need to rescale the coordinates in a way similar to [Don96, $\S.2$]. If $B(0, 1)\subset\C^n$ denotes the unit ball in $\C^n$ and $\chi=(z_1, \dots , z_n):U\rightarrow B(0, \, 1)$ is the chart of coordinates $z_1, \dots , z_n$ centred at $x$ already used above, let $\alpha_0$ stand for the closed $(1, \, 1)$-form on $B(0, \, 1)$ such that $\chi^{\star}\alpha_0=\alpha_{|U}$. If $\widetilde{\chi} : 1/\sqrt{k}\,U\rightarrow B(0, \, 1)$ is the chart of rescaled coordinates $w_j:=\sqrt{k}\,z_j$ ($j=1, \dots , n$), then $(k\alpha)_{|1/\sqrt{k}\,U}=\widetilde{\chi}^{\star}\alpha_0$. Since the curvature form $\alpha_k$ of $L_k$ is close to $k\alpha$ (see (\ref{eqn:prop-def-alpha-k})), if we denote by $\alpha_0^{(k)}$ the closed $2$-form on $B(0, \, 1)$ for which $\alpha_{k|1/\sqrt{k}\,U}=\widetilde{\chi}^{\star}\alpha_0^{(k)}$, we have $||\alpha_0^{(k)}-\alpha_0||_{C^{\infty}}\leq C/k^{1/b_2}$ in consequence of (\ref{eqn:prop-def-alpha-k}). Now $\widetilde{\chi}$ lifts to a connection-preserving bundle map

$$\widetilde{\chi} : L_{k|\frac{1}{\sqrt{k}}U} \rightarrow \xi_k$$

\noindent where $\xi_k\rightarrow B(0, \, 1)$ is the algebraically trivial complex line bundle endowed with the connection of matrix $A_k$ (see \ref{subsection:local-model}). Thus local sections $s\in C^{\infty}(\frac{1}{\sqrt{k}}U, \, L_k)$ identify with sections of $\xi_k$ over $B(0, \, 1)$.

Applying the {\it a priori} estimate to the elliptic operator $\Delta_k''$ on the interior of $V$ (say on some open subset $V'\Subset V$), we get the following Sobolev $W^2(V')$-norm estimate\!\!:

\begin{equation}\label{eqn:elliptic-est}||s_{nh}||^2_{W^2(\frac{1}{\sqrt{k}}\,V')} \leq C\, \bigg(||\Delta_k''s_{nh}||^2_{L^2(\frac{1}{\sqrt{k}}\,V)} + ||s_{nh}||^2_{L^2(\frac{1}{\sqrt{k}}\,V)}\bigg).\end{equation}

\noindent The constant $C>0$ depends only on the ellipticity constant of $\Delta_k''$, hence only on the principal part of $\Delta_k''$ and this is independent of $k$. Indeed, the operators $\bar\partial_k$ have the same principal part ($=\bar\partial$, see (\ref{eqn:dbar-k-local})) for all $k\in\N^{\star}$, hence the Laplacians $\Delta_k''$ have the same principal part for all $k\in\N^{\star}$. Thus $C>0$ is independent of $k$. We have, using (\ref{eqn:_nh-L2-est-improved}) for the second inequality below, that

\begin{equation}\label{eqn:s_nh-L2-L2}||s_{nh}||^2_{L^2(\frac{1}{\sqrt{k}}\,V)}\leq ||s_{nh}||^2\leq C(X, \, \omega)\,\delta_k,   \hspace{2ex} k\gg 1.\end{equation}

\noindent On the other hand, for $s=\theta u^k=s_h + s_{nh}$, we have (since $\theta = 1$ on $V$)

$$s_{nh} = u^k - s_h  \hspace{3ex}  \mbox{on} \hspace{1ex} V,$$

\noindent so

\begin{eqnarray}\label{eqn:Delta_k''s_nh}\nonumber ||\Delta_k''s_{nh}||^2_{L^2(\frac{1}{\sqrt{k}}\,V)} & = & ||\Delta_k''u^k - \Delta_k''s_h||^2_{L^2(\frac{1}{\sqrt{k}}\,V)} \leq 2\,||\Delta_k''u^k||^2_{L^2(\frac{1}{\sqrt{k}}\,V)} + 2\,||\Delta_k''s_h||^2 \\
            & \leq & 2\,||\Delta_k''u^k||^2_{L^2(\frac{1}{\sqrt{k}}\,V)} + 2\frac{C}{k^{2 + 2\varepsilon}}||s||^2,\end{eqnarray}

\noindent having used the fact that $s_h\in{\cal H}_k$ (see Definition \ref{Def:space-ah-sections}) and that $||s_h|| \leq ||s||$. To estimate $||\Delta_k''u^k||^2_{L^2(1/\sqrt{k}\,V)}$ from above, we see that $||\Delta_k''u^k||^2_{L^2(1/\sqrt{k}\,V)}\leq ||\Delta_k''u^k||^2$ (if we extend $u^k$ to $X$ by setting $(u^k)_{|X\setminus U}\equiv 0$) and

\begin{equation}\nonumber\langle\langle\Delta_k''u^k, \, u^k\rangle\rangle = ||\bar\partial_ku^k||^2 \leq \frac{C}{k^{2/b_2}}\, ||u^k||^2,\end{equation}

\noindent having used (\ref{eqn:dbar-k-uk-almost-C0-L2}). Hence

\begin{eqnarray}\label{eqn:Delta_k''-u_k}\nonumber||\Delta_k''u^k||^2_{L^2(\frac{1}{\sqrt{k}}\,V)} & \leq & ||\Delta_k''u^k||^2 \leq \frac{C^2}{k^{4/b_2}}\,||u^k||^2\\
  & \leq & \mbox{Vol}_{\omega}(X)\,\frac{C^2}{k^{4/b_2}}\,||u^k||^2_{C^0} = \mbox{Vol}_{\omega}(X)\,\frac{C^2}{k^{4/b_2}},\end{eqnarray}

\noindent where the last identity holds since $||u^k||^2_{C^0} = 1$ whenever $\alpha\geq 0$.

 Putting (\ref{eqn:Delta_k''s_nh}) and (\ref{eqn:Delta_k''-u_k}) together, we get

$$||\Delta_k''s_{nh}||^2_{L^2(\frac{1}{\sqrt{k}}\,V)}\leq \frac{C(X, \, \omega)}{k^{4/b_2}} + \frac{C}{k^{2 + 2\varepsilon}}||s||^2$$

\noindent and combining this with (\ref{eqn:elliptic-est}) and (\ref{eqn:s_nh-L2-L2}) we finally get\!\!:

$$||s_{nh}||^2_{W^2(\frac{1}{\sqrt{k}}\,V')}\leq \frac{C(X, \, \omega)}{k^{4/b_2}} + \frac{C}{k^{2 + 2\varepsilon}}||s||^2 + C(X, \, \omega)\,\delta_k,  \hspace{3ex} k\gg 1.$$

\noindent Since $s=\theta u^k$, we have $||s||_{C^0}\leq ||u^k||_{C^0}=1$, hence $||s||\leq \mbox{Vol}_{\omega}(X)\,||s||_{C^0}\leq\mbox{Vol}_{\omega}(X)$. Thus the above estimate gives

\begin{equation}\label{eqn:s_nh-W2}||s_{nh}||^2_{W^2(\frac{1}{\sqrt{k}}\,V')}\leq C(X, \, \omega)\, \bigg(\frac{1}{k^{4/b_2}} + \delta_k\bigg),  \hspace{3ex} k\gg 1\end{equation}

\noindent because $2+2\varepsilon > 4/b_2$ (since $0<2\varepsilon < 4/b_2$ with $2\varepsilon$ very close to $4/b_2$).

 We can now go from this $W^2$-norm estimate to a $W^{2p}$-estimate for any $p\in\N^{\star}$. Indeed, the {\it a priori} estimate applied to the elliptic operator $(\Delta_k'')^p$ on $V'\Subset V$ gives

\begin{equation}\label{eqn:elliptic-est-p}||s_{nh}||^2_{W^{2p}(\frac{1}{\sqrt{k}}\,V')} \leq C\, \bigg(||(\Delta_k'')^ps_{nh}||^2_{L^2(\frac{1}{\sqrt{k}}\,V)} + ||s_{nh}||^2_{L^2(\frac{1}{\sqrt{k}}\,V)}\bigg).\end{equation}

\noindent Repeating the above arguments for every $p\in\N^{\star}$ and using the Sobolev embedding theorem, we finally get

\begin{equation}\label{eqn:s_nh-Cinfty}||s_{nh}||^2_{C^{\infty}(\frac{1}{\sqrt{k}}\,V')}\leq C(X, \, \omega)\, \bigg(\frac{1}{k^{4/b_2}} + \delta_k\bigg),   \hspace{3ex} k\gg 1.\end{equation}

 We have thus constructed a global {\it approximately holomorphic} section $s_h=\theta u^k - s_{nh}\in {\cal H}_k$ of $L_k$ that peaks at an arbitrary point $x\in X$ given beforehand.

\begin{Prop}\label{Prop:peak-section} Suppose $\alpha>0$ is a $d$-closed $(1, \, 1)$-form on a compact Hermitian manifold $(X, \, \omega)$. Then, for every $x\in X$ and every $k\in\N^{\star}$, there exists a global section $s_h=s_h^{(k)}\in{\cal H}_k$ of $L_k$ such that, for all $k\gg 1$, we have\!\!: \\

\noindent $(i)$\, $\displaystyle 1 - C(X, \, \omega)\, \bigg(\frac{1}{k^{4/b_2}} + \delta_k\bigg) \leq |s_h(x)|_{h^k}\leq 1 + C(X, \, \omega)\, \bigg(\frac{1}{k^{4/b_2}} + \delta_k\bigg).$ \\

\noindent In particular, $s_h(x)\neq 0$ if $k$ is large enough. \\

\noindent $(ii)$\,$C(X, \, \omega)\, (1-\delta_k^{1/2}) \leq ||s_h|| \leq C(X, \, \omega)\,(1+\delta_k^{1/2}),$ \\

\noindent $(iii)$\,$- C(X, \, \omega)\bigg(\frac{1}{k^{4/b_2}} + \delta_k\bigg)^{1/2} \leq ||s_h||_{C^0(\frac{1}{\sqrt{k}}\,V')}-C \leq C(X, \, \omega)\bigg(\frac{1}{k^{4/b_2}} + \delta_k\bigg)^{1/2}.$

\end{Prop}

\noindent {\it Proof.} The statement follows immediately from the above considerations. For $s= \theta u^k$, we have $s(x)=u^k(0)=1$, so $(i)$ follows from (\ref{eqn:s_nh-Cinfty}) (in which the $C^0(1/\sqrt{k}\,V')$-norm on the left suffices). 

 To get $(ii)$, recall that $||s_h||^2 = ||s||^2 - ||s_{nh}||^2$, use (\ref{eqn:_nh-L2-est-improved}) and notice that $||s||\leq \mbox{Vol}_{\omega}(X)||s||_{C^0}$ while $||s||_{C^0} = ||u^k||_{C^0(U)}$ is bounded independently of $k$ since $u^k(0)=1$ and $u^k$ is non-increasing on a neighbourhood of $0$. Finally $(iii)$ follows from the same arguments as $(i)$ and $(ii)$.  \hfill $\Box$

\vspace{2ex}

 We end this subsection by noticing that all the above estimates still hold if $u^k$ is replaced by any $m$-jet of {\it approximately holomorphic} sections of $L_k$ at $x$ (cf. Definition \ref{Def:ah-jets}). Indeed, as pointed out before that definition, any linear combination of expressions of the form $z_1^{m_1}\cdots z_n^{m_n}\, u^k$ (with $m_1, \dots , m_n\in\N$) satisfies the same estimates (\ref{eqn:dbar-k-uk-almost-Cinf}) and (\ref{eqn:dbar-k-uk-almost-C0-L2}) as $u^k$ does. It follows that if we start off with a global section $s\in C^{\infty}(X, \, L_k)$ obtained by multiplying a local jet (with coefficients, say, $c_{(m_1, \dots , m_n)}\in\C$) by a cut-off function $\theta$\!\!:

$$s(z):=\theta(z)\, c_{(m_1, \dots , m_n)}z_1^{m_1}\cdots z_n^{m_n}\, u^k(z),  \hspace{3ex}  z\in X,$$

\noindent the non-anti-holomorphic component $s_{nh}$ satisfies the $C^{\infty}$-estimate (\ref{eqn:s_nh-Cinfty}) for $k$ large enough. Moreover, setting $m:=m_1 + \cdots + m_n$, we have

$$\frac{1}{m_1!\cdots m_n!}\,\frac{\partial^m(c_{(m_1, \dots , m_n)}z_1^{m_1}\cdots z_n^{m_n}u^k)}{\partial z_1^{m_1}\cdots\partial z_n^{m_n}}(0) = c_{(m_1, \dots , m_n)}e^{-\frac{k}{2}\varphi(0)}=c_{(m_1, \dots , m_n)},$$

\noindent so $|m_1!\cdots m_n!\,c_{(m_1, \dots , m_n)}| - \varepsilon_k\leq |\frac{\partial^m s_h}{\partial z_1^{m_1}\cdots\partial z_n^{m_n}}(0)| \leq |m_1!\cdots m_n!\,c_{(m_1, \dots , m_n)}| + \varepsilon_k$ on a neighbourhood of $x$, where we have denoted $\varepsilon_k:= C(X, \, \omega)\,(\frac{1}{k^{4/b_2}}+\delta_k)\downarrow 0$ (when $k\rightarrow +\infty$) the right-hand term of (\ref{eqn:s_nh-Cinfty}). Therefore $|\frac{\partial^m s_h}{\partial z_1^{m_1}\cdots\partial z_n^{m_n}}(0)|\neq 0$ for $k\gg 1$ if $c_{(m_1, \dots , m_n)}\neq 0$. Thus we obtain the following addition to Proposition \ref{Prop:peak-section}.

\begin{Prop}\label{Prop:peak-section-addition} The assumptions are the same as in Proposition \ref{Prop:peak-section}. Fix $x\in X$ and local holomorphic coordinates $z_1, \dots , z_n$ centred on $x$. Then, for every $k\in\N^{\star}$ and every $m_1, \dots , m_n\in\N$, there exists a global approximately holomorphic section $s_h=s_h^{(k, \, (m_1, \dots , m_n))}\in{\cal H}_k$ of $L_k$ such that 

$$\frac{\partial^{m_1+\cdots + m_n}s_h}{\partial z_1^{m_1}\cdots\partial z_n^{m_n}}(x)\neq 0  \hspace{3ex}  \mbox{if}\,\,\, k\,\,\,\mbox{is large enough}.$$

\noindent (We have denoted by $s_h$ both the global section of $L_k$ and the function that represents it in a local trivialisation of $L_k$ on a neighbourhood of $x=z(0)$.)

\end{Prop}

 This means that ${\cal H}_k$ generates all $m$-jets of {\it approximately holomorphic} sections of $L_k$ at any $x\in X$ for any $m\in\N$.

\section{Approximately holomorphic projective embeddings}\label{section:proj-embed}

 In this section we prove the {\it transcendental} analogue of the Kodaira Embedding Theorem\!\!: any $C^{\infty}$ $d$-closed {\it positive definite} $(1, \, 1)$-form $\alpha >0$ on a compact complex manifold $X$ with a (possibly) {\it non-rational} De Rham cohomology class\footnote{Giving such an $\alpha$ is, of course, equivalent to giving a K\"ahler metric of (possibly) non-rational class on $X$.} $\{\alpha\}\in H_{DR}^2(X, \, \R)$ defines, by means of the spaces ${\cal H}_k$ of Definition \ref{Def:space-ah-sections}, {\it approximately holomorphic} embeddings of $X$ into complex projective spaces $\Proj^{N_k}$ for $k$ large enough. The classical Kodaira Embedding Theorem corresponds to the case when $\{\alpha\}$ is {\it integral} (or merely {\it rational})\!\!: one gets then genuine {\it holomorphic} projective embeddings.

 As noticed in $(i)$ of Proposition \ref{Prop:peak-section}, for every point $x\in X$ one can find a global section of $L_k$ belonging to ${\cal H}_k$ that does not vanish at $x$ if $k$ is large enough. In other words, the sections in ${\cal H}_k$ have no common zeroes for $k\gg 1$. Hence the {\it approximately holomorphic} Kodaira maps

\begin{equation}\label{eqn:Kod-maps1}\Phi_k:X\longrightarrow P{\cal H}_k \simeq\Proj^{N_k},  \hspace{3ex} \Phi_k(z)=H_z:=\{s\in{\cal H}_k \,\,;\,\, s(z)=0\}\end{equation}

\noindent (where $P{\cal H}_k$ stands for the complex projective space whose points are hyperplanes of ${\cal H}_k$), defined equivalently by choosing any orthonormal basis $(\sigma_{k, \, l})_{0\leq l\leq N_k}$ of ${\cal H}_k$ and putting

\begin{equation}\label{eqn:Kod-maps2}\Phi_k:X\longrightarrow\Proj^{N_k},  \hspace{3ex} \Phi_k(z):=[\sigma_{k, \, 0}(z) : \dots : \sigma_{k, \, N_k}(z)],\end{equation}

\noindent are everywhere defined on $X$.

\begin{The}\label{The:Kod-embed} Suppose there exists a $C^{\infty}$ $d$-closed positive definite $(1, \, 1)$-form $\alpha >0$ (i.e. a K\"ahler metric) on a compact complex manifold $X$. Then, for every $k$ large enough, the map $\Phi_k:X\longrightarrow\Proj^{N_k}$ is an embedding.

\end{The}

 We will adapt to our non-integrable context the classical strategy of proof of the Kodaira Embedding Theorem. It remains to prove that, for $k\gg 1$, the sections in ${\cal H}_k$ separate points on $X$ (i.e. for any distinct points $x, y\in X$, there exists a section $s\in{\cal H}_k$ such that $s(x)\neq 0\in (L_k)_x$ but $s(y)=0\in (L_k)_y$) and generate $1$-jets of sections of $L_k$ at every point $x\in X$. We will work our new arguments peculiar to the present context and the necessary modifications of the classical integrable case into the presentation of Demailly's book [Dem97, chapter VII, $\S.13$]. 

\vspace{2ex}

 We begin by analysing a situation to which our case will be reduced.

\begin{Lem}\label{Lem:point-div-separation} Suppose there exists an effective divisor $E$ on $X$ and let $x\in X\setminus\mbox{Supp}\, E$. Fix any $m\in\N$. Then, for every $k$ large enough, there exists an approximately holomorphic section $\tau\in{\cal H}_k$ of $L_k$ such that

$$\tau(x)\neq 0  \hspace{3ex} \mbox{and} \hspace{3ex} \tau \,\,\, \mbox{vanishes on}\,\, E \,\, \mbox{to order}\,\, \geq m+1.$$

\end{Lem}

\noindent {\it Proof.} Consider open subsets $U, V$ such that $x\in V\Subset U\Subset X$, $L_k$ is trivial on $U$ and $U\cap\mbox{Supp}\,E=\emptyset$. Let $u^k\in C^{\infty}(U, \, L_k)$ be the local Gaussian section of $L_k$ {\it peaking} at $x$ constructed in subsection \ref{subsection:local-model}. For a cut-off function $\theta : X \rightarrow \R$ such that $\theta \equiv 1$ on $V$ and $\mbox{Supp}\,\theta\Subset U$ set, as in $\S.$\ref{subsection:global-peak},

$$s:=\theta u^k\in C^{\infty}(X, \, L_k).$$ 

\noindent We know by $(i)$ of Proposition \ref{Prop:peak-section} that, in the splitting $s=s_h + s_{nh}$ with $s_h\in{\cal H}_k$ and $s_{nh}\in{\cal N}_k$, we have $s_h(x)\neq 0$. However, there is a priori no reason for $s_h$ to vanish on $E$. 

 Let $h\in H^0(X, \, {\cal O}(E))$ be the canonical holomorphic section of the holomorphic line bundle associated with $E$. Thus $\mbox{div}\,(h)=E$. Since $s$ vanishes identically on a neighbourhood of $\mbox{Supp}\, E$, we get a smooth section of the $C^{\infty}$ complex line bundle $F_k:={\cal O}(-(m+1)E)\otimes L_k$ by setting\!\!:

$$\sigma:=h^{-(m+1)}\otimes s\in C^{\infty}(X, \, F_k).$$

\noindent Put any $C^{\infty}$ Hermitian metric $h_E$ on ${\cal O}(E)$ and endow the holomorphic line bundle ${\cal O}(-(m+1)E)$ with the Chern connection associated with the induced metric $h_E^{-(m+1)}$. Together with the connection $D_k=\partial_k + \bar\partial_k$ of $L_k$ (cf. $\S.$ \ref{subsection:setting}) this induces a connection

$$D_{F_k} = \partial_{F_k} + \bar\partial_{F_k}$$

\noindent on $F_k$ that is compatible with the metric $h_{F_k}$ induced on $F_k$ by $h_E^{-(m+1)}$ and the metric $h_k$ of $L_k$. Since ${\cal O}(-(m+1)E)$ is holomorphic, we actually have $\bar\partial_{F_k} = \bar\partial_k$ in the following sense\!\!: the $(0, \, 1)$-type connection $\bar\partial_k$ of $L_k = {\cal O}((m+1)E)\otimes F_k$ splits as\!\!:

\begin{equation}\label{eqn:dbar_k-splitting}\bar\partial_k = \bar\partial\otimes\mbox{Id}_{F_k} + \mbox{Id}_{(m+1)E}\otimes\bar\partial_{F_k}.\end{equation}

\noindent It follows that the formal adjoint of $\bar\partial_k$ w.r.t. $h_k=h_E^{m+1}\otimes h_{F_k}$ splits as\!\!:

\begin{equation}\label{eqn:dbar-star_k-splitting}\bar\partial_k^{\star} = \bar\partial^{\star}\otimes\mbox{Id}_{F_k} + \mbox{Id}_{(m+1)E}\otimes\bar\partial_{F_k}^{\star}.\end{equation}

\noindent The corresponding curvature form of $F_k$ reads\!\!:

$$\frac{i}{2\pi}\Theta(F_k) = \gamma_m + \alpha_k,$$

\noindent where $\gamma_m:=\frac{i}{2\pi}\Theta({\cal O}(-(m+1)E)$ is of type $(1, \, 1)$. Hence

$$\frac{i}{2\pi}\Theta(F_k)^{0, \, 2} = \alpha_k^{0, \, 2} \hspace{3ex}  \mbox{and}  \hspace{3ex} \frac{i}{2\pi}\Theta(F_k)^{1, \, 1} = \gamma_m + \alpha_k^{1, \, 1},$$

\noindent so using (\ref{eqn:prop-def-alpha-k-bis}) we see that

\begin{equation}\label{eqn:prop-curv-Fk}\bigg|\bigg|\frac{i}{2\pi}\Theta(F_k)^{1, \, 1} - k(\alpha + \frac{1}{k}\gamma_m)\bigg|\bigg|_{C^{\infty}} \leq\frac{C}{k^{\frac{1}{b_2}}} \hspace{3ex}  \mbox{and}  \hspace{3ex} \bigg|\bigg|\frac{i}{2\pi}\Theta(F_k)^{0, \, 2}\bigg|\bigg|_{C^{\infty}} \leq\frac{C}{k^{\frac{1}{b_2}}},\end{equation}

\noindent where $\alpha + \frac{1}{k}\gamma_m > 0$ for all $k$ large enough since $\alpha >0$ by assumption.

 This shows that the sequence $(F_k)_{k\geq 1}$ of $C^{\infty}$ line bundles on $X$ is {\it asymptotically holomorphic} in the same way as the sequence $(L_k)_{k\geq 1}$ (introduced in $\S.$ \ref{subsection:setting}) is. We can thus apply to the bundles $F_k$ the results obtained for $L_k$ in the previous sections. In particular, we can define anti-holomorphic Laplace-Beltrami operators 

$$\Delta_{F_k}'':=\bar\partial_{F_k}\bar\partial_{F_k}^{\star} + \bar\partial_{F_k}^{\star}\bar\partial_{F_k}:C^{\infty}_{p, \, q}(X, \, F_k)\rightarrow C^{\infty}_{p, \, q}(X, \, F_k)$$

\noindent and spaces of {\it approximately holomorphic} sections of $F_k$ analogous to those of Laeng (cf. Definition \ref{Def:space-ah-sections})\!\!:

$${\cal H}_{F_k}:=\bigoplus\limits_{\mu\leq\frac{C}{k^{1+\varepsilon}}}E^{0, \, 0}_{\Delta_{F_k}''}(\mu)\subset C^{\infty}(X, \, F_k)$$

\noindent which induce orthogonal splittings

$$C^{\infty}(X, \, F_k) = {\cal H}_{F_k} \oplus {\cal N}_{F_k},$$

\noindent where ${\cal N}_{F_k}:=({\cal H}_{F_k})^{\perp}$ (cf. (\ref{eqn:delta-def-Hk})). Accordingly $\sigma$ splits as

$$C^{\infty}(X, \, F_k)\ni\sigma=h^{-(m+1)}\otimes s = \sigma_h + \sigma_{nh}, \hspace{3ex} \mbox{with} \hspace{3ex} \sigma_h\in{\cal H}_{F_k},\,\, \sigma_{nh}\in{\cal N}_{F_k}.$$

\noindent By Proposition \ref{Prop:s_nh-L2-est}, $\sigma_{nh}$ satisfies (in the same notation) the $L^2$-estimate

\begin{equation}\label{eqn:sigma_nh-L2-est}||\sigma_{nh}||^2\leq\frac{4}{\delta_0k}\bigg(1+ \frac{C}{\delta_0\delta^2}\,\frac{1}{k^{1+\frac{2}{b_2}}}\bigg)||\bar\partial_{F_k}\sigma||^2, \hspace{3ex} k\geq k_{\delta},\end{equation}

\noindent while $\sigma_{nh}$ also satisfies by $\S.$\ref{subsection:global-peak} the $C^{\infty}$-estimate (\ref{eqn:s_nh-Cinfty}).

 Now set $\xi:=h^{m+1}\otimes\sigma_{nh}\in C^{\infty}(X, \, L_k)$ and 

$$\tau: = s - \xi = h^{m+1}\otimes \sigma_h\in C^{\infty}(X, \, L_k).$$

\noindent It is clear that $\tau$ vanishes to order $\geq m+1$ on $E$, by construction. On the other hand, estimate (\ref{eqn:sigma_nh-L2-est}) reads

\begin{equation}\label{eqn:h-xi-L2-est}||h^{-(m+1)}\otimes\xi||^2\leq\frac{4}{\delta_0k}\bigg(1+ \frac{C}{\delta_0\delta^2}\,\frac{1}{k^{1+\frac{2}{b_2}}}\bigg)||h^{-(m+1)}\otimes\bar\partial_ks||^2, \hspace{3ex} k\gg 1.\end{equation}

\noindent Now $||h^{-(m+1)}\otimes\xi||^2\geq C_1\,||\xi||^2$ for a constant $C_1>0$ independent of $k$ (depending only on $||h^{m+1}||_{C^0}$) because the holomorphic section $h^{m+1}$ is bounded above on the compact manifold $X$. Meanwhile $s$ vanishes identically, hence so does $\bar\partial_ks$, on a neighbourhood $W$ of $\mbox{Supp}\, E$ in $X$. So $||h^{-(m+1)}\otimes\bar\partial_ks||^2\leq C_2\, ||\bar\partial_ks||^2$ for a constant $C_2>0$ independent of $k$ (depending only on $\inf\limits_{X\setminus W}|h^{m+1}|$). Therefore, (\ref{eqn:h-xi-L2-est}) yields\!\!:

$$||\xi||^2\leq\frac{4\,C_2}{C_1\delta_0k}\bigg(1+ \frac{C}{\delta_0\delta^2}\,\frac{1}{k^{1+\frac{2}{b_2}}}\bigg)||\bar\partial_ks||^2, \hspace{3ex} k\gg 1.$$

\noindent This estimate is the analogue of (\ref{eqn:s_nh-L2-est}) for $\xi$ in place of $s_{nh}$. It leads, by a repetition of the arguments of subsection \ref{subsection:global-peak}, to the $C^{\infty}$-estimate for $\xi$ analogous to (\ref{eqn:s_nh-Cinfty}) which, in turn, leads to the analogue for $\xi$ of $(i)$ of Proposition \ref{Prop:peak-section}\!\!:

$$\tau(x) = (s - \xi)(x)\neq 0.$$

 It remains to check that $\tau\in{\cal H}_k$. If we choose an orthonormal basis $(f_{k, \, l})_{l\in\N}$ of $C^{\infty}(X, \, F_k)$ consisting of eigenvectors of $\Delta_{F_k}''$ and denote by $(\mu_{k, \, l})_{l\in\N}$ the corresponding eigenvalues, we have

$$\sigma_h = \sum\limits_{l=0}^{N_k'}c_l\, f_{k, \, l}, \hspace{3ex} \mbox{hence} \hspace{3ex} \tau = h^{m+1}\otimes\sigma_h = \sum\limits_{l=0}^{N_k'}c_l\, h^{m+1}\otimes f_{k, \, l},$$

\noindent where $N_k' + 1 :=\mbox{dim}{\cal H}_{F_k}$ and $\mu_{k, \, l}\leq \frac{C}{k^{1+\varepsilon}}$ for $l=0, \dots , N_k'$. Then, using formulae (\ref{eqn:dbar_k-splitting}) and (\ref{eqn:dbar-star_k-splitting}), we get\!\!:

\begin{eqnarray}\nonumber\Delta_k''(h^{m+1}\otimes f_{k, \, l}) & = & \bar\partial_k^{\star}\bar\partial_k(h^{m+1}\otimes f_{k, \, l}) = \bar\partial_k^{\star}(h^{m+1}\otimes\bar\partial_kf_{k, \, l})\\
\nonumber  & = & h^{m+1}\otimes\Delta_k''(f_{k, \, l}) = h^{m+1}\otimes (\mu_{k, \, l}f_{k, \, l}),\end{eqnarray}

\noindent which shows that $h^{m+1}\otimes f_{k, \, l}$ is an eigenvector for $\Delta_k''$ corresponding to the eigenvalue $\mu_{k, \, l}\leq \frac{C}{k^{1+\varepsilon}}$ for all $l=0, \dots , N_k'$. Hence $\tau\in{\cal H}_k$.  \hfill $\Box$

\vspace{2ex}

 We can now show that ${\cal H}_k$ separates points on $X$ (and even more).

\begin{Lem}\label{Lem:point-separation} Let $x, y\in X$ such that $x\neq y$. Fix any $m\in\N$. Then, for every $k$ large enough, there exists an approximately holomorphic section $\tau\in{\cal H}_k$ of $L_k$ such that

$$\tau(x)\neq 0  \hspace{3ex} \mbox{and} \hspace{3ex} \tau \hspace{1ex} \mbox{vanishes at}\,\, y\,\, \mbox{to order}\,\, \geq m+1.$$

\end{Lem}

\noindent {\it Proof.} Let $\pi:\widetilde{X}\longrightarrow X$ be the blow-up of $y$ in $X$ and let $E$ be the exceptional divisor. Then $\pi^{\star}\alpha$ is a $C^{\infty}$ $(1, \, 1)$-form on $\widetilde{X}$ (since $\pi$ is holomorphic) satisfying

$$\pi^{\star}\alpha\geq 0 \hspace{2ex} \mbox{on}\,\,\widetilde{X}  \hspace{3ex} \mbox{and} \hspace{3ex} \pi^{\star}\alpha > 0 \hspace{2ex} \mbox{on}\,\,\widetilde{X}\setminus\mbox{Supp}\,E.$$

\noindent Since ${\cal O}(E)_{|E}\simeq{\cal O}_{\Proj(T_xX)}(-1)$, we can equip ${\cal O}(E)_{|E}$ with the smooth metric coming from ${\cal O}_{\Proj^{n-1}}(-1)$ ($n:=\mbox{dim}_{\C}X$) and then extend it in an arbitrary way to a smooth metric of ${\cal O}(E)$. Thus there exists $k_0\in\N^{\star}$ such that $\frac{i}{2\pi}\Theta({\cal O}(-E)) + k_0\pi^{\star}\alpha > 0$ on $\widetilde{X}$. It follows that

\begin{equation}\label{eqn:strict-pos-up}(m+1)\frac{i}{2\pi}\Theta({\cal O}(-E)) + k\pi^{\star}\alpha > 0 \hspace{2ex}\mbox{on}\,\,\widetilde{X} \hspace{2ex}\mbox{for all}\,\, k\geq k_0(m+1).\end{equation}

 If we equip the $C^{\infty}$ complex line bundle $F_k:= {\cal O}(-(m+1)E)\otimes \pi^{\star}L_k$ with the smooth metric induced from the metrics of ${\cal O}(E)$ and $L_k$, the $(1, \, 1)$-component of the associated curvature $2$-form reads

$$\frac{i}{2\pi}\Theta(F_k)^{1, \, 1} = k\bigg(\frac{m+1}{k}\,\frac{i}{2\pi}\Theta({\cal O}(-E)) + \pi^{\star}\alpha\bigg) + \pi^{\star}(\alpha_k^{1, \, 1}-k\alpha) \hspace{3ex} \mbox{on}\,\,\widetilde{X},$$

\noindent while $\frac{i}{2\pi}\Theta(F_k)^{0, \, 2}=\pi^{\star}\alpha_k^{0, \, 2}$, where, thanks to (\ref{eqn:prop-def-alpha-k-bis}), we have

$$||\pi^{\star}(\alpha_k^{1, \, 1}-k\alpha)||_{C^{\infty}}\leq\frac{C}{k^{\frac{1}{b_2}}} \hspace{2ex} \mbox{and} \hspace{2ex} ||\pi^{\star}\alpha_k^{0, \, 2}||_{C^{\infty}}\leq\frac{C}{k^{\frac{1}{b_2}}} \hspace{3ex} \mbox{for}\,\, k\geq 1.$$

\noindent These relations compare to (\ref{eqn:prop-curv-Fk}). Given the strict positivity property (\ref{eqn:strict-pos-up}), the $C^{\infty}$ {\it approximately holomorphic} line bundles $F_k\rightarrow\widetilde{X}$ are analogous to the line bundles $F_k\rightarrow X$ of the proof of Lemma \ref{Lem:point-div-separation}. Thus if we take open neighbourhoods $V\Subset U\Subset X$ of $x$ in $X$ such that $y\notin \overline{U}$, a cut-off function $\theta$ and the section $s=\theta u^k\in C^{\infty}(X, \, L_k)$ peaking at $x$ as in the proof of Lemma \ref{Lem:point-div-separation}, we can run the argument of that proof for $\pi^{\star}s\in C^{\infty}(\widetilde{X}, \, \pi^{\star}L_k)$ on $\widetilde{X}$ in place of $s$ on $X$. Keeping the notation of the proof of Lemma \ref{Lem:point-div-separation} (possibly up to a $\widetilde{\,}\,$), we get sections $\tilde\xi=h^{m+1}\otimes\sigma_{nh}\in C^{\infty}(\widetilde{X}, \, \pi^{\star}L_k)$ (so $\tilde\xi$ vanishes to order $\geq m+1$ on $E$) and

$$\tilde\tau:=\pi^{\star}s - \tilde\xi = h^{m+1}\otimes\sigma_h\in C^{\infty}(\widetilde{X}, \,\pi^{\star}L_k)$$

\noindent with $\tilde\tau(\tilde{x})\neq 0$ (where $\tilde{x}:=\pi^{-1}(x)$) and $\tilde\tau$ vanishing to order $\geq m+1$ on $E$. Since $\pi^{\star}L_k$ is trivial on a neighbourhood of $E$ (because $L_k$ is trivial near $y$), there exists a section $\tau\in{\cal H}_k\subset C^{\infty}(X, \, L_k)$ such that $\tilde\tau=\pi^{\star}s - \tilde\xi=\pi^{\star}\tau$. Since $\tilde\xi$ vanishes to order $\geq m+1$ on $E$ and $\pi_{|\widetilde{X}\setminus\mbox{Supp}\,E}: \widetilde{X}\setminus\mbox{Supp}\,E \rightarrow X\setminus\{y\}$ is a biholomorphism, we see that

$$s - \xi=\tau\in{\cal H}_k\subset C^{\infty}(X, \, L_k) \hspace{2ex} \mbox{on}\,\, X,$$

\noindent where $\xi:=\pi_{\star}\tilde\xi\in C^{\infty}(X, \, L_k)$ vanishes to order $\geq m+1$ at $y$. Since $s$ vanishes identically on a neighbourhood of $y$, $\tau\in{\cal H}_k\subset C^{\infty}(X, \, L_k)$ is the desired section.   \hfill $\Box$

\vspace{3ex}

\noindent {\it End of proof of Theorem \ref{The:Kod-embed}.} The space ${\cal H}_k$ separating points on $X$ (the case $m=0$ in Lemma \ref{Lem:point-separation}) amounts to the Kodaira-type map $\Phi_k$ (cf. (\ref{eqn:Kod-maps1}) or (\ref{eqn:Kod-maps2})) being injective for $k$ large. On the other hand, by the case $m_1 + \cdots + m_n=1$ of Proposition \ref{Prop:peak-section-addition}, the sections in ${\cal H}_k$ generate all $1$-jets of approximately holomorphic sections of $L_k$ at any point $x$. This amounts to $\Phi_k$ being an immersion if $k$ is large enough. The proof of Theorem \ref{The:Kod-embed} is complete.  \hfill  $\Box$

\section{The original form $\alpha$ as a limit}\label{section:a-isometry}

 In this section we prove the analogue for transcendental classes of Tian's almost isometry theorem [Tia90, Theorem A]. We assume throughout that $\alpha>0$ on $X$ but its class $\{\alpha\}\in H^2(X, \, \R)$ need not be rational.

 Let $\Phi_k: X\rightarrow\Proj^{N_k}$ be the approximately holomorphic embedding of the previous section defined by the subspace ${\cal H}_k\subset C^{\infty}(X, \, L_k)$ and let $\omega_{FS}^{(k)}$ denote the Fubini-Study metric of $\Proj^{N_k}$. Then $\frac{1}{k}\Phi_k^{\star}\omega_{FS}^{(k)}$ is again a $d$-closed $C^{\infty}$ $2$-form on $X$ but in general not of type $(1, \, 1)$ (since pull-backs under non-holomorphic maps need not preserve bi-degrees). On the other hand, the current $T_k$ introduced in (\ref{eqn:def1-currents}) is now a genuine $C^{\infty}$ $(1, \, 1)$-form on $X$ (since the sections in ${\cal H}_k$ do not have common zeroes when $\alpha>0$ -- see Proposition \ref{Prop:peak-section}) if $k$ is large enough. In the classical case when the class $\{\alpha\}$ is integral, $T_k$ coincides with $\frac{1}{k}\Phi_k^{\star}\omega_{FS}^{(k)}$ (since $\Phi_k$ is holomorphic in that case, hence it commutes with $\partial\bar\partial$) and is termed the $k^{th}$ Bergman metric on $X$. However, in our case the class $\{\alpha\}$ is non-rational and the above $2$-forms are different for bi-degree reasons. Since $\Phi_k$ is {\it approximately (or asymptotically) holomorphic}, the $(2, \, 0)$ and $(0, \, 2)$-components of $\frac{1}{k}\Phi_k^{\star}\omega_{FS}^{(k)}$ are intuitively expected to converge to zero when $k\rightarrow +\infty$. This fact will be borne out by a calculation below. On the other hand, the $(1, \, 1)$-component of $\frac{1}{k}\Phi_k^{\star}\omega_{FS}^{(k)}$ need not be closed, hence it need not coincide with $T_k$ but we will show that $\frac{1}{k}(\Phi_k^{\star}\omega_{FS}^{(k)})^{1, \, 1}$ and $T_k$ converge to the same limit, so they are in a sense {\it asymptotically} equal. Moreover, we will prove that this limit is the original K\"ahler form $\alpha$ as was the case in [Tia90] when $\{\alpha\}$ was integral.

\begin{The}\label{The:a-isometry} Suppose there exists a K\"ahler metric $\alpha >0$ on a compact complex manifold $X$. For an arbitrary orthonormal basis $(\sigma_{k, \, l})_{l\in\N}$ of ${\cal H}_k$, set 

\begin{equation}\label{eqn:def1bis-currents}T_k:=\alpha + \frac{i}{2\pi k}\partial\bar\partial\log\sum\limits_{l=0}^{N_k}|\sigma_{k, \, l}|^2_{h_k}\end{equation}

\noindent (and note that $T_k$ is independent of the choice of orthonormal basis). Then\!\!: \\

\noindent $(a)$\, $||T_k-\alpha||_{C^2} = O(\frac{1}{\sqrt{k}})$ as $k\rightarrow +\infty$. 

\noindent $(b)$\, $||\frac{1}{k}\Phi_k^{\star}\omega_{FS}^{(k)} - T_k||_{C^2}= O(\frac{1}{\sqrt{k}})$ as $k\rightarrow +\infty$. \\

\noindent In particular, $T_k$ and $\frac{1}{k}\Phi_k^{\star}\omega_{FS}^{(k)}$ converge to $\alpha$, while $\frac{1}{k}(\Phi_k^{\star}\omega_{FS}^{(k)})^{1, \, 1} - T_k$, $\frac{1}{k}(\Phi_k^{\star}\omega_{FS}^{(k)})^{2, \, 0}$ and $\frac{1}{k}(\Phi_k^{\star}\omega_{FS}^{(k)})^{0, \, 2}$ converge to zero in the $C^2$-topology as $k\rightarrow +\infty$.

\end{The}

 The rest of this section will be devoted to proving these statements. Notice that it suffices to prove the estimates locally with constants independent of the open subset chosen. So fix a point $x\in X$, local holomorphic coordinates $z_1, \dots , z_n$ centred at $x$ and a local trivialisation $\theta_k$ of $L_k$ over a neighbourhood $U$ of $x$ as in $\S.$\ref{subsection:local-model}. For any global section $\sigma$ of $L_k$, denote by $f$ the $C^{\infty}$ function on $U$ that represents $\sigma$ w.r.t. $\theta_k$. In particular, the functions $f_{k, \, l}$ represent on $U$ the sections $\sigma_{k, \, l}$ forming an orthonormal basis of ${\cal H}_k$\!\!:

$$\sigma_{k, \, l}\stackrel{\theta_k}{\simeq}f_{k, \, l}\otimes e^{(k)} \hspace{2ex} \mbox{on}\hspace{1ex} U, \hspace{2ex} \mbox{for} \hspace{1ex} l=0, \dots , N_k.$$

 We will choose an orthonormal basis $(\sigma_{k, \, l})_{0\leq l\leq N_k}$ of ${\cal H}_k$ that will enable us to compute derivatives of $f_{k, \, l}$ at $x$ in the same way as Tian chose his basis in [Tia90, (3.7)]. Since the evaluation linear map

$$ev_x : {\cal H}_k \longrightarrow \C,  \hspace{3ex}  \sigma \longmapsto f(x),$$

\noindent does not vanish identically (cf. $(i)$ of Proposition \ref{Prop:peak-section}), its kernel is a hyperplane in ${\cal H}_k$ and we can choose $\sigma_{k, \, 0}\in{\cal H}_k\setminus\ker\,(ev_x)$ such that $\sigma_{k, \, 0}\perp\ker\,(ev_x)$. Thus $f_{k, \, 0}(x)\neq 0$. Since the evaluation linear map

$$ev_x\frac{\partial}{\partial z_1} : \ker\,(ev_x)\longrightarrow \C, \hspace{3ex} \sigma \longmapsto \frac{\partial f}{\partial z_1}(x),$$

\noindent does not vanish identically (otherwise ${\cal H}_k$ would not generate the {\it approximately holomorphic} $1$-jet $z_1\, e^{-k/2\,\varphi}$ at $x$ -- see Proposition \ref{Prop:peak-section-addition}), its kernel is a hyperplane in $\ker\,(ev_x)$ and we can choose $\sigma_{k, \, 1}\in\ker\,(ev_x)\setminus\ker(ev_x\frac{\partial}{\partial z_1})$ such that $\sigma_{k,\, 1}\perp\ker(ev_x\frac{\partial}{\partial z_1})$. Thus $f_{k,\, 1}(x)=0$ but $\frac{\partial f_{k,\, 1}}{\partial z_1}(x)\neq 0$. We can thus construct inductively a decreasing sequence of subspaces

$${\cal H}_k\supset\ker(ev_x)\supset\ker\bigg(ev_x\frac{\partial}{\partial z_1}\bigg)\supset\dots\supset\ker\bigg(ev_x\frac{\partial}{\partial z_n}\bigg)\supset\ker\bigg(ev_x\frac{\partial^2}{\partial z_1^2}\bigg)\supset\dots,$$

\noindent each containing the next as a hyperplane (since ${\cal H}_k$ generates {\it approximately holomorphic} jets at $x$), choose $\sigma_{k,\, l}\in\ker(ev_x\frac{\partial}{\partial z_l})^{\perp}\subset\ker\,(ev_x\frac{\partial}{\partial z_{l-1}})$ for $l=1, \dots , n$ and $\sigma_{k,\, n+1}, \dots , \sigma_{k,\, N_k}$ in the analogous way. Normalising each $\sigma_{k,\, l}$ to norm $1$, we get an orthonormal basis of ${\cal H}_k$ such that (cf. [Tia90, (3.7)])\!\!:

\begin{eqnarray}\label{eqn:OBchoice}\nonumber f_{k,\, 0}(0) & \neq & 0 \hspace{2ex} \mbox{and} \hspace{1ex} f_{k,\, l}(0)=0 \hspace{2ex} \mbox{for all} \hspace{1ex} l\geq 1,\\
\nonumber \frac{\partial f_{k,\, l}}{\partial z_1}(x) & = & \dots = \frac{\partial f_{k,\, l}}{\partial z_{l-1}}(x) = 0  \hspace{2ex} \mbox{but} \hspace{1ex} \frac{\partial f_{k,\, l}}{\partial z_l}(x) \neq 0 \hspace{2ex} \mbox{for all} \hspace{1ex} 1\leq l\leq n,\\
\nonumber \frac{\partial f_{k,\, n+1}}{\partial z_1}(x) & = & \dots = \frac{\partial f_{k,\, n+1}}{\partial z_n}(x) = 0  \hspace{2ex} \mbox{but} \hspace{1ex} \frac{\partial^2 f_{k,\, n+1}}{\partial z_1^2}(x) \neq 0,\\
 \frac{\partial^2 f_{k,\, l}}{\partial z_1^2}(x) & = & 0 \hspace{2ex} \mbox{for all} \hspace{1ex} l\geq n+2.\end{eqnarray}

 To streamline the calculations, we may assume that the local holomorphic coordinates $z_1, \dots , z_n$ about $x$, chosen originally as in (\ref{eqn:choice-local-coord}) (where $\lambda_j(x)>0$ for all $j$ since $\alpha>0$), have been further rescaled such that

\begin{equation}\label{eqn:rescaled-coord}\alpha(x)=\frac{i}{2\pi}\,\partial\bar\partial\varphi(x)=\frac{i}{2\pi}\,\sum\limits_{j=1}^ndz_j\wedge d\bar{z_j}.\end{equation}

\begin{Lem}\label{Lem:hol-dir-deriv} Suppose the local coordinates $z_1, \dots , z_n$ about $x$ have been rescaled as in (\ref{eqn:rescaled-coord}) and the local potential $\varphi$ of $\alpha$ has been chosen as in (\ref{eqn:phi-choice}) (with each $\lambda_j(x)$ replaced by $1$). For $u^k$ defined in $\S.$\ref{subsection:local-model} write $u^k=f_k\otimes e^{(k)}$ on $U$ as in (\ref{eqn:u^k-section}), where $e^{(k)}$ denotes the local frame of $L_k$ w.r.t. $\theta_k$. Then, for all $j=1, \dots , n$, we have\!\!:

$$\frac{\partial f_k}{\partial\bar{z_j}}(z) = -\frac{k}{2}\,f_k(z)\,\frac{\partial\varphi}{\partial\bar{z_j}}(z), \hspace{2ex} z\in U.$$

\noindent In particular, $\frac{\partial f_k}{\partial\bar{z_j}}(0)=0$ and $\frac{\partial^2 f_k}{\partial z_j\partial\bar{z_j}}(0)=-\frac{k}{2}\,f_k(0)$ for all $j=1, \dots , n$.

\noindent The statements still hold if we replace $u^k$ by any jet $z_1^{m_1}\cdots z_n^{m_n}\, u^k$.

\end{Lem}

\noindent {\it Proof.} Since $\bar\partial_{kA}(u^k)=0$ on $U$ (see $\S.$\ref{subsection:local-model}), we have

$$0=\bar\partial_{kA}(u^k) = (\bar\partial f_k + 2\pi k f_kA^{0, \, 1})\otimes e_k \hspace{2ex} \mbox{hence} \hspace{2ex} \bar\partial f_k = -\frac{k}{2}f_k\,\bar\partial\varphi  \hspace{2ex} \mbox{on}\hspace{1ex} U,$$

\noindent having used the identity $A^{0, \, 1}=1/4\pi\,\bar\partial\varphi$ on $U$. The first statement of the Lemma follows. The first part of the second statement follows by taking $z=0$ in the first one and using (\ref{eqn:phi-choice}) to see that $\partial\varphi/\partial\bar{z_j}(0)=0$ for every $j$. The second part follows by applying $\partial/\partial z_j$ in the first statement and using (\ref{eqn:phi-choice}). These properties still hold for jets since $\bar\partial_{kA}(z_1^{m_1}\cdots z_n^{m_n}\,u^k)=0$ on $U$.  \hfill $\Box$

\vspace{2ex}

 We need one more preliminary observation in the spirit of [Tia90, Lemmas 2.1, 2.2, 2.3] before performing the actual calculations. We can apply to every {\it approximately holomorphic} jet $z_1^{m_1}\cdots z_n^{m_n}\,u^k$ at $x$ the procedure described in $\S.$\ref{subsection:corrections}: multiply it by a cut-off function $\theta$ (with $\theta\equiv 1$ on $V$ and $\mbox{Supp}\,\theta\Subset U$) and then take the orthogonal projections $s_h:=s^{(k),\, h}_{(m_1, \dots , m_n)}$, resp. $s_{nh}:=s^{(k),\, nh}_{(m_1, \dots , m_n)}$, of $s:=\theta\,z_1^{m_1}\cdots z_n^{m_n}\,u^k\in C^{\infty}(X, \, L_k)$ onto ${\cal H}_k$, resp. ${\cal N}_k$. So 

$$s^{(k),\, h}_{(m_1, \dots , m_n)} = \theta\,z_1^{m_1}\cdots z_n^{m_n}\,u^k - s^{(k),\, nh}_{(m_1, \dots , m_n)}   \hspace{3ex} \mbox{on}\hspace{1ex} X.$$

\noindent The results obtained in $\S.$\ref{subsection:corrections} and $\S.$\ref{subsection:global-peak} starting from $u^k$ still apply if we start off with $z_1^{m_1}\cdots z_n^{m_n}\,u^k$ instead. Therefore $s^{(k),\, nh}_{(m_1, \dots , m_n)}$ satisfies the $L^2$-estimate (\ref{eqn:s_nh-L2-est}) on $X$ and the $C^{\infty}$-estimate (\ref{eqn:s_nh-Cinfty}) on $1/\sqrt{k}\,V'\Subset 1/\sqrt{k}\,V\Subset 1/\sqrt{k}\,U$.

 On the other hand, given the properties (\ref{eqn:OBchoice}) satisfied by the orthonormal basis $(\sigma_{k, \, l})_{l\in\N}$ of ${\cal H}_k$, it is not hard to see, after normalising each peak section $s^{(k),\, h}_{(m_1, \dots , m_n)}$ to $\widetilde{s}^{(k)}_{(m_1, \dots , m_n)}$ of $L^2$-norm $1$, that $\sigma_{k, \,  0}$ is {\it close} to $\widetilde{s}^{(k)}_{(0, \dots , 0)}$, $\sigma_{k, \,  l}$ is {\it close} to $\widetilde{s}^{(k)}_{(0, \dots, 1, \dots, 0)}$ (with $1$ in the $l^{th}$ spot) for every $l=1, \dots , n$, $\sigma_{k, \,  n+1}$ is {\it close} to $\widetilde{s}^{(k)}_{(2, 0, \dots , 0)}$, etc. In Tian's holomorphic case, this was an $L^2$-norm proximity (cf. [Tia90, Lemma 3.1]). We can show furthermore that the forms in each of these pairs are close to each other in the $C^{\infty}$-norm on a neighbourhood of $x$.

\begin{Lem}\label{Lem:s-sigma-closeness} We have\!\!:   $||\sigma_{k, \,  0} - \widetilde{s}^{(k)}_{(0, \dots , 0)}||_{C^{\infty}(\frac{1}{\sqrt{k}}\,V')}\leq C(X, \, \omega)\,(\frac{1}{k^{4/b_2}} + \delta_k)$ and

$$||\sigma_{k, \,  l} - \widetilde{s}^{(k)}_{(0, \dots, 1, \dots, 0)}||_{C^{\infty}(\frac{1}{\sqrt{k}}\,V')}\leq C(X, \, \omega)\,(\frac{1}{k^{4/b_2}} + \delta_k), \hspace{2ex} \mbox{for all} \hspace{1ex} l=1, \dots , n,$$

$$||\sigma_{k, \,  n+1} - \widetilde{s}^{(k)}_{(2, 0, \dots , 0)}||_{C^{\infty}(\frac{1}{\sqrt{k}}\,V')}\leq C(X, \, \omega)\,(\frac{1}{k^{4/b_2}} + \delta_k).$$

\end{Lem}

\noindent {\it Proof.} The proof is essentially contained in [Tia90], we only re-interpret it in the light of our estimate (\ref{eqn:s_nh-Cinfty}). It is well known that, for every $l\in\N$, any homogeneous polynomial $P(X_1, \dots , X_d)\in\R[X_1, \dots , X_d]$ of degree $l$ for which $\Delta_{\R^d}P=0$ (where $\Delta_{\R^d}$ is the usual Laplacian of $\R^d$) restricts to an eigenvector $P_{|S^{d-1}}$ of the (non-positive) Laplace-Beltrami operator $\Delta_{S^{d-1}}$ of the unit sphere $S^{d-1}\subset\R^d$ of eigenvalue $-l(l+d-2)$. Furthermore, all {\it spherical harmonics} arise as such restrictions. For $z\in\C^n$, it follows that 

$$\int\limits_{S^{2n-1}}z_1^{m_1}\cdots z_n^{m_n}\,\bar{z}_1^{p_1}\cdots \bar{z}_n^{p_n}\,d\sigma(z)=0  \hspace{2ex}  \mbox{for all}\hspace{1ex} (m_1, \dots , m_n)\neq (p_1, \dots , p_n),$$

\noindent hence, after integrating by parts, we get for all $(m_1, \dots , m_n)\neq (p_1, \dots , p_n)$\!\!:

$$i^n\,\int\limits_{|z|\leq R}z_1^{m_1}\cdots z_n^{m_n}\,\bar{z}_1^{p_1}\cdots \bar{z}_n^{p_n}\,\rho(|z|^2)\,dz_1\wedge d\bar{z}_1\wedge \cdots \wedge dz_n\wedge d\bar{z}_n=0$$

\noindent for any function $\rho(z)=\rho(|z|)$ depending only on $|z|$. This means that different {\it approximately holomorphic} jets $z_1^{m_1}\cdots z_n^{m_n}\,u^k$ and $z_1^{p_1}\cdots z_n^{p_n}\,u^k$ at $x$ (with $(m_1, \dots , m_n)\neq (p_1, \dots , p_n)$) are mutually orthogonal on a neighbourhood of $x$. The orthogonality defect on $X$ between the global sections $s^{(k), \, h}_{(m_1, \dots , m_n)}$ and $s^{(k),\, h}_{(p_1, \dots , p_n)}$ is due solely to the distorsion introduced by the correcting sections $s^{(k),\,nh}_{(m_1, \dots , m_n)}$ and $s^{(k), nh}_{(p_1, \dots , p_n)}$. However, these correcting sections satisfy estimate (\ref{eqn:s_nh-Cinfty}), so they are {\it small} in $C^{\infty}$-norm on $1/\sqrt{k}\,V'$. The lemma follows.  \hfill $\Box$

\vspace{2ex}

 With the choice (\ref{eqn:OBchoice}) of an orthonormal basis of ${\cal H}_k$ we shall now estimate the latter term in the right-hand side of (\ref{eqn:def1bis-currents}) in the same way as Tian did in the {\it holomorphic} case. Extra terms containing $\bar{z_j}$-derivatives of $f_{k, \, l}$ appear in our {\it approximately holomorphic} context compared to Tian's case; Lemma \ref{Lem:hol-dir-deriv} will contribute to their estimates. The terms containing only $z_j$-derivatives of $f_{k, \, l}$ can be estimated as in [Tia90] and we will be rather brief on details. However, we will spell out in detail the estimates of the new terms peculiar to our non-holomorphic case. It suffices to obtain $C^2$-estimates at the fixed point $x$ that are independent of $x$. 

  Thanks to Lemma \ref{Lem:s-sigma-closeness}, we can compute the derivatives of (\ref{eqn:def1bis-currents}) at $x$ as if the $f_{k, \, l}$ were the {\it approximately holomorphic} jets $u^k$ (for $l=0$), $z_l\,u^k$ (for $1\leq l\leq n$), $z_1^2\,u^k$ (for $l=n+1$), etc. Indeed, only distortions very small in $C^{\infty}$-norm on $1/\sqrt{k}\, V'$ are introduced by this substitution, hence the derivatives at $x$ are only distorted by $O(1/k^{4/b_2} + \delta_k)$. Thus by (\ref{eqn:av-u-hk-rep}), (\ref{eqn:av-u-hk-rep_jets}) we have\!\!:

\begin{equation}\label{eqn:def2-currents}T_{k|U}:=\alpha_{|U} + \frac{i}{2\pi k}\partial\bar\partial\log\sum\limits_{l=0}^{N_k}|f_{k, \, l}|^2,\end{equation}

\noindent where $|f_{k, \, l}|$ denotes the modulus of $f_{k, \, l}$.

 Proving that the latter term in the right-hand side of (\ref{eqn:def2-currents}) converges to {\it zero} amounts to proving that, for every $r, s=1, \dots , n$, the $\partial^2/\partial z_r\partial\bar{z_s}$-derivative of $1/k\cdot\log\sum|f_{k, \, l}|^2$ converges to {\it zero}. Using an orthogonal transformation, it suffices to prove this fact for $r=s=1$.

\vspace{2ex} 

 We begin with the $C^0$-estimate. 

\begin{Lem}\label{Lem:C0-est} With the above choices we have\!\!:

$$\bigg|\frac{1}{k}\,\frac{\partial^2\log\sum\limits_{l=0}^{N_k}|f_{k, \, l}|^2}{\partial z_1\partial\bar{z_1}}(x)\bigg|\leq \frac{C}{k} \hspace{3ex} \mbox{for all} \hspace{1ex} k\gg 1,$$

\noindent where $C>0$ is a constant independent of $x$.

\end{Lem}

\noindent {\it Proof.} Straightforward calculations give\!\!:

\begin{eqnarray}\label{eqn:C0-terms}\nonumber & & \frac{1}{k}\frac{\partial^2\log\sum\limits_{l=0}^{N_k}|f_{k, \, l}|^2}{\partial z_1\partial\bar{z}_1}(x) = \frac{1}{k}\,\frac{\partial}{\partial z_1}\bigg(\frac{\sum\limits_{l=0}^{N_k}f_{k, \, l}\,\frac{\partial\bar{f}_{k,\, l}}{\partial\bar{z}_1} + \sum\limits_{l=0}^{N_k}\bar{f}_{k,\, l}\,\frac{\partial f_{k, \, l}}{\partial\bar{z}_1}}{\sum\limits_{l=0}^{N_k}|f_{k, \, l}|^2}\bigg)(x)\\
\nonumber & = &  \frac{1}{k}\,\frac{\sum\limits_{l=0}^{N_k}\bigg|\frac{\partial f_{k, \, l}}{\partial z_1}\bigg|^2 + \sum\limits_{l=0}^{N_k}\bigg|\frac{\partial f_{k,\, l}}{\partial\bar{z_1}}\bigg|^2 + \sum\limits_{l=0}^{N_k} f_{k,\, l}\,\frac{\partial^2\bar{f}_{k,\, l}}{\partial z_1\partial\bar{z}_1} + \sum\limits_{l=0}^{N_k} \bar{f}_{k,\, l}\,\frac{\partial^2f_{k, \, l}}{\partial z_1\partial\bar{z}_1}}{\sum\limits_{l=0}^{N_k}|f_{k, \, l}|^2}(x)\\
\nonumber &  & \hspace{23ex} - \frac{1}{k}\,\frac{\bigg|\sum\limits_{l=0}^{N_k}f_{k, \, l}\,\frac{\partial\bar{f}_{k,\, l}}{\partial\bar{z}_1} + \sum\limits_{l=0}^{N_k}\bar{f}_{k,\, l}\,\frac{\partial f_{k, \, l}}{\partial\bar{z}_1}\bigg|^2}{(\sum\limits_{l=0}^{N_k}|f_{k, \, l}|^2)^2}(x)\\
\nonumber & = & \frac{1}{k}\,\frac{\bigg|\frac{\partial f_{k, \, 1}}{\partial z_1}(x)\bigg|^2}{|f_{k,\, 0}(x)|^2} + \frac{1}{k}\,\frac{\sum\limits_{l=1}^{N_k}\bigg|\frac{\partial f_{k, \, l}}{\partial\bar{z}_1}(x)\bigg|^2 - \bigg|\frac{\partial f_{k, \, 0}}{\partial\bar{z}_1}(x)\bigg|^2}{|f_{k, \, 0}(x)|^2}\\
\nonumber &  & \hspace{21ex} + \frac{1}{k}\,\frac{f_{k, \, 0}(x)\,\frac{\partial^2\bar{f}_{k,\, 0}}{\partial z_1\partial\bar{z}_1}(x) + \bar{f}_{k, \, 0}(x)\,\frac{\partial^2f_{k, \, 0}}{\partial z_1\partial\bar{z}_1}(x)}{|f_{k, \, 0}(x)|^2}\\
 &  & \hspace{3ex} -\frac{1}{k}\,\frac{f_{k, \, 0}(x)^2\,\frac{\partial\bar{f}_{k,\, 0}}{\partial z_1}(x)\,\frac{\partial\bar{f}_{k,\, 0}}{\partial\bar{z}_1}(x) + \bar{f}_{k, \, 0}(x)^2\,\frac{\partial f_{k,\, 0}}{\partial z_1}(x)\,\frac{\partial f_{k,\, 0}}{\partial\bar{z}_1}(x)}{|f_{k, \, 0}(x)|^4}.\end{eqnarray}

\noindent We have used (\ref{eqn:OBchoice}) in $f_{k, \, l}(x)=0$ for all $l\geq 1$ and in $\frac{\partial f_{k, \, l}}{\partial z_1}(x)=0$ for all $l\geq 2$.

 By Lemma \ref{Lem:hol-dir-deriv}, the anti-holomorphic first order derivatives $\partial/\partial\bar{z}_j$ vanish at $0=z(x)$ in the case of jets, so we are left with calculating the first term and the terms containing second order derivatives $\partial^2/\partial z_1\partial\bar{z}_1$ at $0$ in (\ref{eqn:C0-terms}). 

 By Lemma \ref{Lem:hol-dir-deriv}, we further have

$$\frac{\partial^2f_{k, \, 0}}{\partial z_1\partial\bar{z}_1}(0)=-\frac{k}{2}\,f_{k, \, 0}(0), \hspace{2ex} \mbox{hence} \hspace{1ex} \bar{f}_{k, \, 0}(0)\,\frac{\partial^2f_{k, \, 0}}{\partial z_1\partial\bar{z}_1}(0) = -\frac{k}{2}\,|f_{k, \, 0}(0)|^2.$$

\noindent Thus $\bar{f}_{k, \, 0}(0)\,\frac{\partial^2f_{k, \, 0}}{\partial z_1\partial\bar{z}_1}(0)$ is real and therefore equals its conjugate $f_{k, \, 0}(0)\,\frac{\partial^2\bar{f}_{k, \, 0}}{\partial z_1\partial\bar{z}_1}(0)$. It follows that at $x=0$ we have

\begin{equation}\label{eqn:2nd-order-terms_C0}\frac{1}{k}\,\frac{f_{k, \, 0}(x)\,\frac{\partial^2\bar{f}_{k,\, 0}}{\partial z_1\partial\bar{z}_1}(x) + \bar{f}_{k, \, 0}(x)\,\frac{\partial^2f_{k, \, 0}}{\partial z_1\partial\bar{z}_1}(x)}{|f_{k, \, 0}(x)|^2} = - 1.\end{equation}

\noindent As for the first term in the last sum on the right-hand side of (\ref{eqn:C0-terms}), it can be estimated in the same way as the analogous term was estimated in [Tia90]. Indeed, if $\sim$ stands for equality up to terms $O(1/k)$ involving constants independent of $x$, we have as in [Tia90, $\S.2$, $\S.3$]\!\!:

$$\frac{1}{k}\,\frac{\bigg|\frac{\partial f_{k, \, 1}}{\partial z_1}(x)\bigg|^2}{|f_{k,\, 0}(x)|^2}\sim \frac{C^2_{(1, 0, \dots , 0)}}{k\,C^2_{(0, \dots , 0)}},$$

\noindent where $C_{(1, 0, \dots , 0)}$, $C_{(0, \dots , 0)}$ are the coefficients of the normalised peak sections\!\!:

$$\widetilde{s}^{(k)}_{(1, \dots , 0)}(z)=C_{(1, 0, \dots , 0)}\,(z_1\,u^k(z) - s^{(k),\, nh}_{(1, \dots , 0)}(z)),  \hspace{2ex} \widetilde{s}_{(0, \dots , 0)}(z)=C_{(0, 0, \dots , 0)}\,(u^k(z) - s^{(k),\, nh}_{(0, \dots , 0)}(z))$$

\noindent given by the formulae\!\!:

$$C^2_{(1, 0, \dots , 0)}\sim\frac{1}{\int\limits_{|z|\leq\frac{\log k}{\sqrt{k}}}|z_1|^2\,e^{-k\varphi(z)}\,dV_{\omega}(z)} \hspace{2ex}  \mbox{and} \hspace{2ex} C^2_{(0, 0, \dots , 0)}\sim\frac{1}{\int\limits_{|z|\leq\frac{\log k}{\sqrt{k}}}e^{-k\varphi(z)}\,dV_{\omega}(z)}.$$

\noindent In the rescaled coordinates of (\ref{eqn:rescaled-coord}), $\lambda_j(x)$ becomes $1$ in (\ref{eqn:phi-choice}), so 

$$e^{-\varphi(z)}=1-|z|^2 + O(|z|^3), \hspace{2ex} z\in U,$$

\noindent hence

\begin{equation}\nonumber C^2_{(1, \,0, \dots , 0)}\sim\frac{1}{\int\limits_{|z|\leq\frac{\log k}{\sqrt{k}}}|z_1|^2\,(1-|z|^2)^k \,dz_1\wedge d\bar{z}_1\wedge\dots\wedge dz_n\wedge d\bar{z}_n},\end{equation}

\begin{equation}\nonumber C^2_{(0, \,0, \dots , 0)} \sim \frac{1}{\int\limits_{|z|\leq\frac{\log k}{\sqrt{k}}}(1-|z|^2)^k \,dz_1\wedge d\bar{z}_1\wedge\dots\wedge dz_n\wedge d\bar{z}_n},\end{equation}

\noindent from which it follows that

\begin{eqnarray}\label{eqn:z-1-deriv_frac-f0}\nonumber \frac{1}{k}\,\frac{\bigg|\frac{\partial f_{k, \, 1}}{\partial z_1}(x)\bigg|^2}{|f_{k,\, 0}(x)|^2} & \sim & \frac{\int\limits_{|z|\leq\frac{\log k}{\sqrt{k}}}(1-|z|^2)^k \,dz_1\wedge d\bar{z}_1\wedge\dots\wedge dz_n\wedge d\bar{z}_n}{k\,\int\limits_{|z|\leq\frac{\log k}{\sqrt{k}}}|z_1|^2\,(1-|z|^2)^k \,dz_1\wedge d\bar{z}_1\wedge\dots\wedge dz_n\wedge d\bar{z}_n}\\
 & \sim & 1+O(\frac{1}{k}),\end{eqnarray}

\noindent the last estimate appearing in [Tia90, proof of Lemma 3.3].

 Putting together (\ref{eqn:2nd-order-terms_C0}) and (\ref{eqn:z-1-deriv_frac-f0}) we see that

$$\frac{1}{k}\,\frac{\bigg|\frac{\partial f_{k, \, 1}}{\partial z_1}(x)\bigg|^2}{|f_{k,\, 0}(x)|^2} + \frac{1}{k}\,\frac{f_{k, \, 0}(x)\,\frac{\partial^2\bar{f}_{k,\, 0}}{\partial z_1\partial\bar{z}_1}(x) + \bar{f}_{k, \, 0}(x)\,\frac{\partial^2f_{k, \, 0}}{\partial z_1\partial\bar{z}_1}(x)}{|f_{k, \, 0}(x)|^2} \sim O(\frac{1}{k}).$$

 This completes the proof of Lemma \ref{Lem:C0-est}.  \hfill $\Box$

\vspace{2ex}

 It follows from the uniformity w.r.t. $x$ of the estimate of Lemma \ref{Lem:C0-est} that $\frac{i}{2\pi k}\,\partial\bar\partial\log\sum\limits_{l=0}^{N_k}\log|f_{k, \, l}|^2$ converges uniformly to {\it zero} on $U$. Thus we have\!\!

\begin{Cor}\label{Cor:Tk-C0-conv} As $k\rightarrow +\infty$, $T_k$ converges to $\alpha$ in the $C^0$-topology.

\end{Cor}

 The $C^1$-estimate is handled in a similar way. It suffices to estimate uniformly the $\partial^3/\partial z_1^2\partial\bar{z}_1$-derivative of $1/k\cdot\log\sum|f_{k, \, l}|^2$ at $x$.

\begin{Lem}\label{Lem:C1-est} With the above choices we have\!\!:

$$\bigg|\frac{1}{k}\,\frac{\partial^3\log\sum\limits_{l=0}^{N_k}|f_{k, \, l}|^2}{\partial z_1^2\partial\bar{z_1}}(x)\bigg|\leq \frac{C}{k} \hspace{3ex} \mbox{for all} \hspace{1ex} k\gg 1,$$

\noindent where $C>0$ is a constant independent of $x$.

\end{Lem}

\noindent {\it Proof.} The third-order $\partial^3/\partial z_1^2\partial\bar{z}_1$-derivatives of the $f_{k, \, l}$ and $\bar{f}_{k, \, l}$ (which would vanish if the $f_{k, \, l}$ were holomorphic, but do not in our case) are handled as follows. Applying $\partial^2/\partial z_j^2$ in the first conclusion of Lemma \ref{Lem:hol-dir-deriv}, we get

$$\frac{\partial^3f_{k, \, l}}{\partial z_j^2\partial\bar{z}_j}(z) = -\frac{k}{2}\,f_{k, \, l}(z)\frac{\partial^3\varphi}{\partial z_j^2\partial\bar{z}_j}(z) - k\,\frac{\partial f_{k, \, l}}{\partial z_j}(z)\,\frac{\partial^2\varphi}{\partial z_j\partial\bar{z}_j}(z) - \frac{k}{2}\,\frac{\partial^2f_{k, \, l}}{\partial z_j^2}(z)\,\frac{\partial\varphi}{\partial\bar{z}_j}(z).$$

\noindent A similar formula is obtained for $\partial^3\bar{f}_{k, \, l}/\partial z_j^2\partial\bar{z}_j$. Taking now $z=0$ and using the facts that $\partial\varphi/\partial z_j(0) = \partial\varphi/\partial\bar{z}_j(0)=0$ (by (\ref{eqn:phi-choice})), $\partial^2\varphi/\partial z_j\partial\bar{z}_j(0)=1$ (again by (\ref{eqn:phi-choice}) in which each $\lambda_j(x)$ becomes $1$ in the rescaled coordinates of (\ref{eqn:rescaled-coord})) and $\partial\bar{f}_{k, \, l}/\partial z_j(0)=0$ (by Lemma \ref{Lem:hol-dir-deriv}), we find for $j=1$ that

$$\frac{\partial^3f_{k, \, l}}{\partial z_1^2\partial\bar{z}_1}(0) = -\frac{k}{2}\,f_{k, \, l}(0)\frac{\partial^3\varphi}{\partial z_1^2\partial\bar{z}_1}(0) - k\,\frac{\partial f_{k, \, l}}{\partial z_1}(0)$$

\noindent and

$$\frac{\partial^3\bar{f}_{k, \, l}}{\partial z_1^2\partial\bar{z}_1}(0) = -\frac{k}{2}\,\bar{f}_{k, \, l}(0)\frac{\partial^3\varphi}{\partial z_1^2\partial\bar{z}_1}(0).$$

 Using these facts and (\ref{eqn:OBchoice}), straightforward but tedious calculations give

\begin{eqnarray}\nonumber \frac{1}{k}\,\frac{\partial^3\log\sum\limits_{l=0}^{N_k}|f_{k, \, l}|^2}{\partial z_1^2\partial\bar{z_1}}(0) & = & \frac{\frac{\partial^2f_{k, \, 1}}{\partial z_1^2}(0)\,\frac{\partial\bar{f}_{k, \, 1}}{\partial\bar{z}_1}(0)}{k\,|f_{k, \, 0}(0)|^2} - \frac{|\frac{\partial f_{k, \, 1}}{\partial z_1}(0)|^2\,\bar{f}_{k, \, 0}(0)\,\frac{\partial f_{k, \, 0}}{\partial z_1}(0)}{k\,|f_{k, \, 0}(0)|^4}\\
\nonumber  & + & 2\,\frac{\bar{f}_{k, \, 0}(0)\,\frac{\partial f_{k, \, 0}}{\partial z_1}(0)}{|f_{k, \, 0}(0)|^2}.\end{eqnarray}

\noindent The two terms on the right-hand side of the first line above also appear in the holomorphic case. The estimates of [Tia90, Lemma 3.2] apply to the sections involved in all of the above expressions thanks to arguments very similar to those recalled in the proof of Lemma \ref{Lem:C0-est}. Using those estimates, we get the following uniform growth rates for the three right-hand terms above\!\!: \\

\noindent $\cdot$\, $O(\frac{k^{(n-1)/2}\,k^{n+1}}{k\,k^{2n}}) = O(\frac{1}{k^{(n+1)/2}})$ (for the first term)\!; \\

\noindent $\cdot$\, $O(\frac{k^{n+1}\,k^{n/2}\,k^{(n/2)-1}}{k\,k^{4n}}) = O(\frac{1}{k^{2n+1}})$ (for the second term)\!; \\

\noindent $\cdot$\, $O(\frac{k^{n/2}\,k^{(n/2)-1}}{k^n}) = O(\frac{1}{k})$ (for the third term). \\

 The contention follows.   \hfill $\Box$

\vspace{2ex}

 The $C^2$ estimate can be proved in the same way and is left to the reader. We have thus proved part $(a)$ of Theorem \ref{The:a-isometry}. We now prove part $(b)$.

 Since for any system of homogeneous coordinates $[w_0:\cdots :w_{N_k}]$ of $\Proj^{N_k}$ the Fubini-Study metric reads $\omega_{FS}^{(k)}=\frac{i}{2\pi}\,\partial\bar\partial\log\sum\limits_{l=0}^{N_k}|w_l|^2$, we get\!\!:

$$\omega_{FS}^{(k)}=\frac{1}{\sum\limits_{l=0}^{N_k}|w_l|^2}\bigg(\sum\limits_{l=0}^{N_k}\frac{i}{2\pi}\,dw_l\wedge d\bar{w}_l - \frac{1}{\sum\limits_{l=0}^{N_k}|w_l|^2}\,\sum\limits_{l,\, r=0}^{N_k}w_r\bar{w}_l\,\frac{i}{2\pi}\,dw_l\wedge d\bar{w}_r\bigg).$$

\noindent Hence

\begin{eqnarray}\nonumber \frac{1}{k}\,\Phi_k^{\star}\omega_{FS}^{(k)} & = & \frac{1}{k\,\sum\limits_{l=0}^{N_k}|f_{k,\, l}|^2}\bigg(\sum\limits_{l=0}^{N_k}\frac{i}{2\pi}\,(\partial f_{k,\, l} + \bar\partial f_{k,\, l})\wedge (\partial\bar{f}_{k,\, l} + \bar\partial\bar{f}_{k,\, l})\\
\nonumber  & - & \frac{1}{\sum\limits_{l=0}^{N_k}|f_{k,\, l}|^2}\,\sum\limits_{l,\, r=0}^{N_k}f_{k,\, r}\,\bar{f}_{k,\, l}\,\frac{i}{2\pi}\,(\partial f_{k,\, l} + \bar\partial f_{k,\, l})\wedge (\partial\bar{f}_{k,\, r} + \bar\partial\bar{f}_{k,\, r})\bigg),\end{eqnarray}

\noindent from which it follows that

\begin{eqnarray}\label{eqn:Phi_k_11}(\frac{1}{k}\,\Phi_k^{\star}\omega_{FS}^{(k)})^{1, \, 1} & = & \frac{1}{k\,\sum\limits_{l=0}^{N_k}|f_{k, \, l}|^2}\bigg(\sum\limits_{l=0}^{N_k}(\frac{i}{2\pi}\,\partial f_{k, \, l}\wedge\bar\partial\bar f_{k, \, l} - \frac{i}{2\pi}\,\partial\bar f_{k, \, l}\wedge\bar\partial f_{k, \, l}) \\
\nonumber & - & \frac{1}{\sum\limits_{l=0}^{N_k}|f_{k, \, l}|^2}\,\,\sum\limits_{l, \, r=0}^{N_k}f_{k, \, r}\,\bar f_{k, \, l}(\frac{i}{2\pi}\,\partial f_{k, \, l}\wedge\bar\partial\bar f_{k, \, r} - \frac{i}{2\pi}\,\partial\bar f_{k, \, r}\wedge\bar\partial f_{k, \, l})\bigg),\end{eqnarray}

\begin{eqnarray}\label{eqn:Phi_k_20}(\frac{1}{k}\,\Phi_k^{\star}\omega_{FS}^{(k)})^{2, \, 0}  =  \frac{1}{k\,\sum\limits_{l=0}^{N_k}|f_{k, \, l}|^2} & \bigg( & \sum\limits_{l=0}^{N_k}\frac{i}{2\pi}\,\partial f_{k, \, l}\wedge\partial\bar f_{k, \, l}\\
 \nonumber & - & \frac{1}{\sum\limits_{l=0}^{N_k}|f_{k, \, l}|^2}\sum\limits_{l, \, r=0}^{N_k}f_{k, \, r}\,\bar f_{k, \, l}\,\frac{i}{2\pi}\,\partial f_{k, \, l}\wedge\partial\bar f_{k, \, r}\bigg),\end{eqnarray}

\begin{eqnarray}\label{eqn:Phi_k_02}(\frac{1}{k}\,\Phi_k^{\star}\omega_{FS}^{(k)})^{0, \, 2}  =  \frac{1}{k\,\sum\limits_{l=0}^{N_k}|f_{k, \, l}|^2} & \bigg( & \sum\limits_{l=0}^{N_k}\frac{i}{2\pi}\,\bar\partial f_{k, \, l}\wedge\bar\partial\bar f_{k, \, l}\\
\nonumber  & - & \frac{1}{\sum\limits_{l=0}^{N_k}|f_{k, \, l}|^2}\sum\limits_{l, \, r=0}^{N_k}f_{k, \, r}\,\bar f_{k, \, l}\,\frac{i}{2\pi}\,\bar\partial f_{k, \, l}\wedge\bar\partial\bar f_{k, \, r}\bigg),\end{eqnarray}

\noindent where $|f_{k, \, l}|$ stands for the modulus of the function $f_{k, \, l}$ that represents $\sigma_{k, \, l}$ in a local trivialisation of $L_k$. Thanks to Lemma \ref{Lem:hol-dir-deriv}, if the $f_{k, \, l}$'s were the actual {\it approximately holomorphic} jets at $x=z(0)$, we would have

\begin{equation}\label{eqn:dbar-point-vanishing}\bar\partial f_{k, \, l}(0) = 0,  \hspace{3ex} l=0, \dots, N_k,\end{equation}

\noindent and implicitly $(\Phi_k^{\star}\omega_{FS}^{(k)})^{2, \, 0}$ and $(\Phi_k^{\star}\omega_{FS}^{(k)})^{0, \, 2}$ would vanish at $x$. Now, as already noticed earlier, the $f_{k, \, l}$'s need not be the jets at $x$ but they lie within a $C^{\infty}$-norm distance $O(1/k^{4/b_2} + \delta_k)$ of the jets by Lemma \ref{Lem:s-sigma-closeness}. This estimate being uniform w.r.t. $x$, we infer that

$$||(\frac{1}{k}\,\Phi_k^{\star}\omega_{FS}^{(k)})^{2, \, 0}||_{C^2}\leq O(\frac{1}{\sqrt{k}})  \hspace{2ex}  \mbox{and} \hspace{2ex} ||(\frac{1}{k}\,\Phi_k^{\star}\omega_{FS}^{(k)})^{0, \, 2}||_{C^2}\leq O(\frac{1}{\sqrt{k}}).$$

On the other hand, straightforward calculations show that

\begin{eqnarray}\label{eqn:T_k_calc1}\nonumber \frac{i}{2\pi k}\,\partial\bar\partial\log\sum\limits_{l=0}^{N_k}|f_{k, \, l}|^2 & \hspace{-2ex} = & \frac{1}{k\,\sum\limits_{l=0}^{N_k}|f_{k, \, l}|^2}\sum\limits_{l=0}^{N_k}(\frac{i}{2\pi}\,\partial f_{k, \, l}\wedge\bar\partial\bar f_{k, \, l} + \frac{i}{2\pi}\,\partial\bar f_{k, \, l}\wedge\bar\partial f_{k, \, l}\\
\nonumber & \hspace{-2ex} &  \hspace{15ex} + f_{k, \, l}\,\frac{i}{2\pi}\,\partial\bar\partial\bar f_{k, \, l} + \bar f_{k, \, l}\,\frac{i}{2\pi}\,\partial\bar\partial f_{k, \, l})\\
\nonumber  & \hspace{-22ex}  & \hspace{-22ex} -\frac{1}{k\,\bigg(\sum\limits_{l=0}^{N_k}|f_{k, \, l}|^2\bigg)^2}\sum\limits_{l, \, r=0}^{N_k}(f_{k, \, l}\,f_{k, \, r}\, \frac{i}{2\pi}\,\partial\bar f_{k, \, l} \wedge\bar\partial\bar f_{k, \, r} + f_{k, \, l}\,\bar f_{k, \, r}\, \frac{i}{2\pi}\,\partial\bar f_{k, \, l}\wedge\bar\partial f_{k, \, r}\\
  & \hspace{-12ex} + & \hspace{-6ex} f_{k, \, r}\,\bar f_{k, \, l}\, \frac{i}{2\pi}\,\partial f_{k, \, l}\wedge\bar\partial\bar f_{k, \, r} + \bar f_{k, \, l}\,\bar f_{k, \, r}\, \frac{i}{2\pi}\,\partial f_{k, \, l}\wedge\bar\partial f_{k, \, r}).\end{eqnarray}

\noindent Notice that the right-hand sides of the expressions (\ref{eqn:Phi_k_11}) and (\ref{eqn:T_k_calc1}) contain precisely the same terms featuring products $\partial f_{k, \, l}\wedge\bar\partial\bar{f}_{k, \, r}$, while all the products containing a factor $\partial\bar{f}_{k, \, l}$ or $\bar\partial f_{k, \, l}$ would vanish at $0=z(x)$ if the $f_{k, \, l}$'s were the actual {\it approximately holomorphic} jets at $0$ by (\ref{eqn:dbar-point-vanishing}). Thus the terms in this latter group are negligible in the $C^{\infty}$-norm by Lemma \ref{Lem:s-sigma-closeness}. The only two terms featuring on the right of (\ref{eqn:T_k_calc1}) but not on the right of (\ref{eqn:Phi_k_11}) are those containing second-order derivatives $i\partial\bar\partial\bar{f}_{k, \, l}$ and $i\partial\bar\partial f_{k, \, l}$. Thus we have\!\!:

\begin{eqnarray}\label{eqn:ddbar-11-rel}\nonumber\frac{i}{2\pi k}\,\partial\bar\partial\log\sum\limits_{l=0}^{N_k}|f_{k, \, l}|^2(x) & \hspace{-3ex} = & \hspace{-3ex} (\frac{1}{k}\,\Phi_k^{\star}\omega_{FS}^{(k)})^{1, \, 1}(x)\\
\nonumber & \hspace{-3ex} + & \hspace{-3ex} \frac{1}{k\,\sum\limits_{l=0}^{N_k}|f_{k, \, l}|^2}\sum\limits_{l=0}^{N_k}(f_{k, \, l}\,\frac{i}{2\pi}\,\partial\bar\partial\bar f_{k, \, l} + \bar f_{k, \, l}\,\frac{i}{2\pi}\,\partial\bar\partial f_{k, \, l})(x)\\
 &  + &  O(\frac{1}{\sqrt{k}}).\end{eqnarray}

\noindent Now by Lemma \ref{Lem:hol-dir-deriv}, if the $f_{k, \, l}$'s were the actual {\it approximately holomorphic} jets at $0=z(x)$, for every $l=0, \dots , N_k$ the following would hold on $U$\!\!:  

$$\bar\partial f_{k, \, l} = -\frac{k}{2}\,f_{k, \, l}\,\bar\partial\varphi  \hspace{2ex} \mbox{hence} \hspace{1ex} \frac{i}{2\pi}\,\partial\bar\partial f_{k, \, l} = -\frac{k}{2}\,f_{k, \, l}\, \frac{i}{2\pi}\,\partial\bar\partial\varphi - \frac{k}{2}\,\frac{i}{2\pi}\,\partial f_{k, \, l}\, \wedge \bar\partial\varphi.$$

\noindent Since $i\partial\bar\partial\varphi = 2\pi\,\alpha$ on $U$ (see $\S.$\ref{subsection:local-model}) and $\bar\partial\varphi(0)=0$ (see (\ref{eqn:phi-choice})), the last identity applied at $z=0$ reads\!\!:

$$\frac{i}{2\pi}\,\partial\bar\partial f_{k, \, l}(0) = -\frac{k}{2}\,f_{k, \, l}(0)\,\alpha(0),$$

\noindent from which we get (at $x=z(0)$)\!\!:

$$\bigg(f_{k, \, l}\,\frac{i}{2\pi}\,\partial\bar\partial\bar f_{k, \, l} + \bar f_{k, \, l}\,\frac{i}{2\pi}\,\partial\bar\partial f_{k, \, l}\bigg)(0) = -k\,|f_{k, \, l}(0)|^2\,\alpha(0),   \hspace{2ex}   l=0, \dots , N_k.$$

\noindent Since $f_{k, \, l}(0)=0$ for all $l\geq 1$ (cf. (\ref{eqn:OBchoice})),
we infer that

$$\frac{1}{k\,\sum\limits_{l=0}^{N_k}|f_{k, \, l}(0)|^2}\,\,\sum\limits_{l=0}^{N_k}(f_{k, \, l}\,\frac{i}{2\pi}\,\partial\bar\partial\bar f_{k, \, l} + \bar f_{k, \, l}\,\frac{i}{2\pi}\,\partial\bar\partial f_{k, \, l})(0) = -\alpha(0).$$

\noindent Since $x=0$ in the chosen coordinates, from (\ref{eqn:ddbar-11-rel}) we get\!\!:

$$\frac{i}{2\pi k}\,\partial\bar\partial\log\sum\limits_{l=0}^{N_k}|f_{k, \, l}|^2(x) =  (\frac{1}{k}\,\Phi_k^{\star}\omega_{FS}^{(k)})^{1, \, 1}(x) - \alpha(x) + O(\frac{1}{\sqrt{k}}).$$

\noindent On the other hand, the left-hand term above equals $T_k(x) - \alpha(x)$ thanks to (\ref{eqn:def2-currents}), so we get \!\!:

$$T_k(x) =  (\frac{1}{k}\,\Phi_k^{\star}\omega_{FS}^{(k)})^{1, \, 1}(x) + O(\frac{1}{\sqrt{k}}).$$

\noindent Since the constant implicit in the quantity $O(\frac{1}{\sqrt{k}})$ is independent of $x$ and corresponds to a $C^{\infty}$ estimate, we have obtained the uniform estimate proving part $(b)$ of Theorem \ref{The:a-isometry}.

 This completes the proof of Theorem \ref{The:a-isometry}.  \hfill $\Box$

\vspace{3ex}

\noindent {\bf References.}\\

\noindent [BU00]\, D. Borthwick, A. Uribe --- {\it Nearly K\"ahlerian Embeddings of Symplectic Manifolds} --- Asian J. Math. {\bf 4} (2000), no. 3, 599–620.

\vspace{1ex}

\noindent [BG81]\, L. Boutet de Monvel, V. Guillemin --- {\it The Spectral Theory of Toeplitz Operators} --- Ann. Math. Stud. {\bf 99}, Princeton University Press, Princeton 1981.

\vspace{1ex}

\noindent [BS76]\, L. Boutet de Monvel, J. Sj\"ostrand --- {Sur la singularit\'e des noyaux de Bergman et de Szeg\"o} --- Ast\'erisque {\bf 34-35} (1976), 123-164.

\vspace{1ex}

\noindent [Dem85a]\, J.-P. Demailly --- {\it Champs magn\'etiques et in\'egalit\'es de Morse pour la $d''$-cohomologie} --- Ann. Inst. Fourier (Grenoble) {\bf 35} (1985), no. 4, 189–229.

\vspace{1ex}

\noindent [Dem85b]\, J.-P. Demailly --- {\it Sur l'identit\'e de Bochner-Kodaira-Nakano en g\'eom\'etrie hermitienne} --- S\'eminaire P. Lelong, P. Dolbeault, H. Skoda (Analyse) 1983/1984, LNM {\bf 1198}, Springer Verlag (1985) 88-97.

\vspace{1ex}

\noindent [Dem92]\, J.-P. Demailly --- {\it Regularization of Closed Positive Currents and Intersection Theory} --- J. Alg. Geom. {\bf 1} (1992) 361-409.

\vspace{1ex}

\noindent [Dem 97] \, J.-P. Demailly --- {\it Complex Analytic and Differential Geometry}---http://www-fourier.ujf-grenoble.fr/~demailly/books.html



\vspace{1ex}

\noindent [DP04]\, J.-P. Demailly, M. Paun --- {\it Numerical Charaterization of the K\"ahler Cone of a Compact K\"ahler Manifold} --- Ann. of Math. (2) {\bf 159(3)} (2004) 1247-1274.

\vspace{1ex}

\noindent [Don90]\, S.K. Donaldson --- {\it Yang-Mills Invariants of Four-manifolds} --- Geometry of low-dimensional manifolds, 1 (Durham, 1989), 5–40, London Math. Soc. Lecture Note Ser. {\bf 150}, Cambridge Univ. Press, Cambridge, 1990.

\vspace{1ex}

\noindent [Don96]\, S.K. Donaldson --- {\it Symplectic Submanifolds and Almost-Complex Geometry} --- J. Diff. Geom. {\bf 44} (1996) 666-705.

\vspace{1ex}

\noindent [Gri66]\, Ph. Griffiths --- {\it The Extension Problem in Complex Analysis II: Embeddings with Positive Normal Bundle} --- Amer. J. of Math. {\bf 88} (1966), 366-446.

\vspace{1ex}

\noindent [Lae02]\, L. Laeng --- {\it Estimations spectrales asymptotiques en g\'eom\'etrie hermitienne} --- Th\`ese de doctorat de l'Universit\'e Joseph Fourier, Grenoble (octobre 2002), http://tel.archives-ouvertes.fr/tel-00002098/en/.

\vspace{1ex}

\noindent [MM02]\, X. Ma, G. Marinescu --- {\it The $\mbox{Spin}^c$ Dirac Operator on High Tensor Powers of a Line Bundle} --- Math. Z. {\bf 240}  651-664 (2002).

\vspace{1ex}

\noindent [Pop08]\, D. Popovici --- {\it Regularisation of Currents with Mass Control and Singular Morse Inequalities} --- J. Diff. Geom. {\bf 80} (2008) 281-326.

\vspace{1ex}

\noindent [Pop09]\, D. Popovici --- {\it Limits of Projective Manifolds under Holomorphic Deformations} --- arXiv e-print math.AG/0910.2032v1.

\vspace{1ex}

\noindent [SZ02]\, B. Shiffman, S. Zelditch --- {\it Asymptotics of almost holomorphic sections of ample line bundles on symplectic manifolds} --- J. Reine Angew. Math. {\bf 544} (2002) 181-222.

\vspace{1ex}

\noindent [Tia90]\, G. Tian --- {\it On a Set of Polarized K\"ahler Metrics on Algebraic Manifolds} --- J. Diff. Geom. {\bf 32} (1990) 99-130.

\vspace{1ex}

\noindent [Zel98]\, S. Zelditch --- {\it Szeg\"o Kernels and a Theorem of Tian} --- Int. Math. Res. Notices {\bf 6} (1998) 317-331.

\vspace{3ex}

\noindent Institut de Math\'ematiques de Toulouse, Universit\'e Paul Sabatier

\noindent 118, route de Narbonne, 31062, Toulouse, France

\noindent Email: popovici@math.univ-toulouse.fr

\end{document}